\pdfoutput=1
\RequirePackage{silence}
\PassOptionsToPackage{table}{xcolor}
\documentclass[12pt,a4paper,svgnames]{amsart}
\usepackage[hmarginratio={1:1},vmarginratio={1:1},lmargin=60.0pt,tmargin=55.0pt]{geometry}

\usepackage{booktabs}
\allowdisplaybreaks
\usepackage{amsmath,amssymb,amsthm}
\usepackage{mathtools,xcolor,colortbl}
\usepackage{enumitem}
\setlist[enumerate]{itemsep=0.15cm,label=\emph{\upshape(\alph*)}}
\usepackage{multicol}
\usepackage{bbm}
\usepackage{seqsplit}

\newcommand{\Tilt}{\mathrm{Tilt}}
\newcommand{\bk}{\mathbf{k}}

\newcommand{\R}{\mathbb{R}}
\newcommand{\C}{\mathbb{C}}
\newcommand{\Z}{\mathbb{Z}}
\newcommand{\N}{\Z_{\geq 0}}
\newcommand{\F}{\mathbb{F}}
\newcommand{\Q}{\mathbb{Q}}

\numberwithin{equation}{section} 
\usepackage{color}

\newcommand{\per}{\nu}
\newcommand{\alnum}{\bar{t}}

\newcommand{\placeholder}{{}_{-}}

\newcommand{\mP}{\mathbb{P}}
\newcommand{\mA}{\mathbb{A}}

\newcommand{\mC}{\mathcal{C}}

\newcommand{\rp}[1][p]{d_{#1}}

\newcommand{\period}{P}

\newcommand{\avalue}{\tau}
\newcommand{\alvalue}{t}

\newcommand{\hone}{h}
\newcommand{\htwo}{\tilde{h}}

\usepackage{tikz}
\tikzset{anchorbase/.style={baseline={([yshift=-0.5ex]current bounding box.center)}},}

\usepackage{colortbl}
\usepackage{xcolor}

\definecolor{orchid}{RGB}{143,40,194}
\definecolor{lava}{RGB}{207,16,32}
\definecolor{mydarkblue}{RGB}{10,10,170}

\usepackage{aliascnt}

\def\NewTheorem#1{%
\newaliascnt{#1}{equation}%
\newtheorem{#1}[#1]{#1}%
\aliascntresetthe{#1}%
\expandafter\def\csname #1autorefname\endcsname{#1}%
}
\def\equationautorefname~#1\null{(#1)\null}

\numberwithin{equation}{subsection}

\NewTheorem{Proposition}
\NewTheorem{Theorem}
\NewTheorem{Main Theorem}
\NewTheorem{Corollary}
\AtEndEnvironment{Corollary}{\null\hfill$\square$}%
\NewTheorem{Lemma}
\theoremstyle{definition}
\NewTheorem{Definition}
\NewTheorem{Notation}
\NewTheorem{Question}
\NewTheorem{Example}
\AtEndEnvironment{Example}{\null\hfill$\Diamond$}%
\NewTheorem{Examples}
\AtEndEnvironment{Examples}{\vskip-10mm\null\hfill$\Diamond$}%

\theoremstyle{remark}
\NewTheorem{Remark}
\AtEndEnvironment{Remark}{\null\hfill$\Diamond$}%
\NewTheorem{Assumption}

\usepackage[hypertexnames=false]{hyperref}
\usepackage{bookmark}
\hypersetup{
pdftoolbar=true,
pdfmenubar=true,
pdffitwindow=false,
pdfstartview={FitH},
pdftitle={Fractal behavior of tensor powers of the two dimensional space in prime characteristic},
pdfauthor={Kevin Coulembier, Pavel Etingof, Victor Ostrik and Daniel Tubbenhauer},
pdfsubject={},
pdfcreator={Kevin Coulembier, Pavel Etingof, Victor Ostrik and Daniel Tubbenhauer},
pdfproducer={Kevin Coulembier, Pavel Etingof, Victor Ostrik and Daniel Tubbenhauer},
pdfkeywords={},
pdfnewwindow=true,
colorlinks=true,
linkcolor=mydarkblue,
citecolor=teal,
filecolor=magenta,
urlcolor=orchid,
linkbordercolor=lava,
citebordercolor=teal,
urlbordercolor=orchid,
linktocpage=true
}

\usepackage{etoolbox}
\def\makeautorefname#1#2{\csdef{#1autorefname}{#2}}
\makeautorefname{section}{Section}%
\makeautorefname{subsection}{Section}%
\makeautorefname{subsubsection}{Section}%

\begin{document}
\title[Fractal behavior...in prime characteristic]{Fractal behavior of tensor powers of the two dimensional space in prime characteristic}
\author[K. Coulembier, P. Etingof, V. Ostrik and D. Tubbenhauer]{Kevin Coulembier, Pavel Etingof, Victor Ostrik and Daniel Tubbenhauer}

\address{K.C.: The University of Sydney, School of Mathematics and Statistics F07, Office Carslaw 717, NSW 2006, Australia, \href{https://www.maths.usyd.edu.au/u/kevinc/}{www.maths.usyd.edu.au/u/kevinc}, \href{https://orcid.org/0000-0003-0996-3965}{ORCID 0000-0003-0996-3965}}
\email{kevin.coulembier@sydney.edu.au}

\address{P.E.: Department of Mathematics, Room 2-282, Massachusetts Institute of Technology, 77 Massachusetts Ave., Cambridge, MA 02139, USA, \href{https://math.mit.edu/~etingof/}{math.mit.edu/~etingof/}, \href{https://orcid.org/0000-0002-0710-1416}{ORCID 0000-0002-0710-1416}}
\email{etingof@math.mit.edu}

\address{V.O.: University of Oregon, Department of Mathematics, Eugene, OR 97403, USA, \newline \href{https://pages.uoregon.edu/vostrik/}{pages.uoregon.edu/vostrik}, \href{https://orcid.org/0009-0002-2264-669X}{ORCID 0009-0002-2264-669X}}
\email{vostrik@math.uoregon.edu}

\address{D.T.: The University of Sydney, School of Mathematics and Statistics F07, Office Carslaw 827, NSW 2006, Australia, \href{http://www.dtubbenhauer.com}{www.dtubbenhauer.com}, \href{https://orcid.org/0000-0001-7265-5047}{ORCID 0000-0001-7265-5047}}
\email{daniel.tubbenhauer@sydney.edu.au}

\begin{abstract}
We study the number of 
indecomposable summands in 
tensor powers of the vector representation of SL2. Our main focus is on positive characteristic where this sequence of 
numbers and its generating function show fractal behavior 
akin to Mahler functions.
\end{abstract}


\addtocontents{toc}{\protect\setcounter{tocdepth}{1}}

\maketitle

\tableofcontents

\section{Introduction}

Let $\mathbf{k}$ be a field and let $\Gamma$ be a group, possibly an algebraic group or affine group scheme over $\mathbf{k}$. For any finite dimensional $\Gamma$-representation $W$ 
over $\mathbf{k}$ let $\nu(W)\in\Z_{\geq 0}$ be an integer such that
\begin{gather*}
W\cong{\textstyle\bigoplus_{i=1}^{\nu(W)}}W_{i},
\end{gather*} 
where $W_{i}$ are indecomposable $\Gamma$-representations. This number is well-defined
by the Krull--Schmidt theorem. Now let $V$ be a finite dimensional $\Gamma$-representation and define
the integer sequence
\begin{gather*}
b_n=b_{n}(V):=\nu(V^{\otimes n}),\; n\in\Z_{\geq 0}.
\end{gather*}
The study of the growth of $b_{n}$ is the main motivation of this paper.

The related question to study the growth of $\ell(V^{\otimes n})$, where $\ell(W)$ is the \emph{length} of $W$, is usually much easier 
and we discuss this along the way, e.g. in 
\autoref{P:FiniteGroups} and \autoref{SS:Length} below.

We write $\sim$ for `asymptotically equal'. One could hope that
\begin{gather}\label{Eq:Ansatz}
b_{n}\sim
\hone(n)\cdot n^{\avalue}\cdot\beta^{n},
\quad
\begin{aligned}
&\hone\colon\N\to\R_{>0}\text{ is a function \emph{bounded away from $0,\infty$}},
\\[-0.1cm]
&n^{\avalue}\text{ is the \emph{subexponential factor}, $\avalue\in\R$},
\\[-0.1cm]
&\beta^{n}\text{ is the \emph{exponential factor}, 
$\beta\in\R_{\geq 1}$}.
\end{aligned}
\end{gather}
In practice, $\hone(n)$ is often a constant or alternates between finitely many constants, but sometimes $\hone(n)$ is more complicated.

We do not know in what generality \autoref{Eq:Ansatz} holds, 
but expressions of this form are very common in the theory of asymptotics of generating functions, see for example 
\cite{FaSe-analytic-combinatorics} or \cite{Mi-analytic-combinatorics} for the relation between counting problems and 
the analysis of generating functions.

\subsection{The exponential factor}

Even without assuming \eqref{Eq:Ansatz} one can define the value $\beta:=\lim_{n\to\infty}\sqrt[n]{b_{n}(V)}$, and it is proven in \cite[Theorem 1.4]{CoOsTu-growth} that
\begin{gather}\label{Eq:ExponentialFactor}
\beta=\dim_{\mathbf{k}}V.
\end{gather}
Thus, \autoref{Eq:ExponentialFactor} gives $(\dim_{\mathbf{k}}V)^{n}$ as the exponential factor 
of the growth rate of $b_{n}$.

The goal of this paper is to provide a more precise asymptotics, that is, one of the form \autoref{Eq:Ansatz}, for the sequence $b_{n}$ in 
the case $\dim_{\mathbf{k}}V=2$ and $\Gamma=GL(V)$ or, equivalently, $\Gamma=SL(V)$. This is a nontrivial case where one can hope to understand the asymptotics explicitly.
In this setting the subexponential factor is determined by \autoref{T:MainTheorem}, as well as \autoref{P:ReductiveGroup} and \autoref{P:FiniteGroups}. 

We do not calculate 
$\hone(n)$, but it appears to be oscillating, cf. \autoref{E:Hfactor}. See 
also \autoref{S:Fourier} where 
we determine $\hone(n)$ for $p=2$.

\subsection{The main theorem -- the subexponential factor}

If $\mathbf{k}$ is finite, then $\Gamma=SL_{2}(\mathbf{k})$ is a finite group 
and the asymptotic behavior of $b_{n}=b_{n}(\mathbf{k}^{2})$ can be easily obtained, see 
for example \autoref{P:FiniteGroups} below or
\cite[Example 8]{LaTuVa-growth-pfdim} which \emph{settle the case of a finite field}. We give an example which is prototypical in this situation:

\begin{Example}\label{E:IntroFiniteGroup}
For $\mathbf{k}=\F_{3}$ we get
\begin{gather*}
(b_{n})_{n\in\N}=(1,1,2,2,6,6,22,22,86,\dots)
\quad\text{and}\quad
b_{n}\sim
\left(\frac{1}{6}+\frac{(-1)^{n}}{12}\right)
\cdot n^{0}\cdot 2^{n}.
\end{gather*}
(We highlight that the subexponential factor is $n^{0}=1$.)
Note that $\lim_{n\to\infty}\frac{b_{n}}{2^{n}}$ does not exist, but 
the ratio $\frac{b_{n}}{2^{n}}$ is rather governed by the periodic function
$\hone(n)=\frac{1}{6}+\frac{(-1)^{n}}{12}$.
\end{Example}

From now on let $V$ be a two dimensional vector space over an infinite field. Let $\Gamma$ be the algebraic group $SL(V)\simeq SL_{2}(\mathbf{k})$. Let 
$b_{n}=b_{n}(V)$ be the sequence as above. As a matter of fact, the sequence $b_{n}$ depends only on the characteristic of $\mathbf{k}$ (which is only implicit in our notation), so we only need to study it for one infinite field for each $p\geq 0$. That is, we have
\begin{gather}\label{Eq:MainSequence}
b_{n}=
\text{number of indecomposable $SL_{2}(\mathbf{k})$-representations in 
$(\mathbf{k}^{2})^{\otimes n}$},
\end{gather}
where $\mathbf{k}$ is either $\C$ for $p=0$ or $\bar{\F}_{p}$ otherwise.

The following example is classical, and \emph{settles the case $p=0$}:

\begin{Example}\label{E:IntroCharZero}
Assume $p=0$. Then $b_{n}$ is the sequence of middle binomial coefficients
\begin{gather*}
(b_{n})_{n\in\N}=\left(\binom{n}{\lfloor n/2\rfloor}\right)_{n\in\N}
=(1,1,2,3,6,10,20,35,70,\dots)
\end{gather*}
Then Stirling's formula implies
\begin{gather*}
b_{n}\sim\sqrt{2/\pi}\cdot n^{-1/2}\cdot 2^{n},
\end{gather*}
as one easily checks. In this case the periodic function is constant, i.e. $\hone(n)=\sqrt{2/\pi}=
\lim_{n\to\infty}\frac{b_{n}}{n^{-1/2}\cdot 2^{n}}$. Since this situation is semisimple, the same formula works for the length instead of the number of indecomposable summands.
\end{Example}

For the remainder of the paper let $p>0$.
Then the asymptotic behavior of the sequence $b_{n}$ turns out to be more complicated.
Namely let us define
\begin{gather}\label{Eq:Alpha}
\alvalue_{p}:=-\frac{1}{2}\log_{p}\frac{2p^{2}}{p+1}=
\frac{1}{2}
\frac{\ln(p+1)-\ln(2p^{2})}{\ln p}=\frac{1}{2}\left(\log_p\frac{p+1}{2}\right)-1.
\end{gather}
As we will see, $\alvalue_{p}$ plays the role of $\avalue$ in \autoref{Eq:Ansatz}. To give some numerical expressions, we have for example:
\begin{gather*}
\alvalue_{2}=-0.7075187496\ldots\,,
\quad
\alvalue_{3}=-0.6845351232\ldots\,,
\quad
\alvalue_{5}=-0.6586969028\ldots\,,
\\
\alvalue_{7}=-0.6437928129\ldots\,,
\quad
\alvalue_{101}=-0.5740278484\ldots\,.
\end{gather*}
It is easy to see that all $\alvalue_{p}$ are irrational. Moreover, all these numbers are transcendental, as follows from the
(logarithm form of the) Gelfond--Schneider theorem which implies that $\log_{p}r$, for $p$ a prime number and $r\in\Q$, is transcendental unless $r$ is a power of $p$.

\begin{Remark}
By very similar arguments, all $\avalue$ that we will see below that are not manifestly rational will be irrational and even transcendental.
\end{Remark}

\begin{Notation}
Throughout, we will use \emph{big O notation} (or rather its generalization 
often called \emph{Bachmann--Landau notation}). Most important for us 
from this family of notations are $\sim$ as above, small $o$, big $O$, big $\Theta$ and $\asymp$. That is, given two real-valued functions $f$, $g$ which are both defined on $\N$ or on $D\subset\R$ for the final point, we will write
\begin{gather}\label{Eq:Asymp}
\begin{aligned}
f\sim g
&\;\Leftrightarrow\;
\forall\varepsilon>0,\,\exists n_{0}
\;\text{ such that }\;
\fbox{$|\tfrac{f(n)}{g(n)}-1|<\varepsilon$},
\;\forall n>n_{0}
,
\\
f\in o(g)
&\;\Leftrightarrow\;
\forall C>0,\,\exists n_{0}
\;\text{ such that }\;
\fbox{$|f(n)|\leq C\cdot g(n)$},
\;\forall n>n_{0}
,
\\
f\in O(g)
&\;\Leftrightarrow\;
\exists C>0,\,\exists n_{0}
\;\text{ such that }\;
\fbox{$|f(n)|\leq C\cdot g(n)$},
\;\forall n>n_{0}
,
\\
f\in\Theta(g)
&\;\Leftrightarrow\;
\exists C_{1},C_{2}>0,\,\exists n_{0}
\;\text{ such that }\;
\fbox{$C_{1}\cdot g(n)\leq f(n)\leq C_{2}\cdot g(n)$},
\;\forall n>n_{0}
,
\\
f\asymp_{D}g
&\;\Leftrightarrow\;
\exists
C>1
\;\text{ such that }\;
\fbox{$C^{-1}\cdot g(x)\leq f(x)\leq C\cdot g(x)$},
\;\forall x\in D.
\end{aligned}
\end{gather}
We sometimes omit the subscript $D$ if no confusion can arise.
\end{Notation}

Our main result is:

\begin{Main Theorem}\label{T:MainTheorem}
For any $p>0$ there exist $C_{1}=C_{1}(p),C_{2}=C_{2}(p)\in\R_{>0}$
such that we have
\begin{gather*}
C_{1}\cdot n^{\alvalue_{p}}\cdot 2^{n}\leq b_{n}\leq C_{2}\cdot n^{\alvalue_{p}}\cdot 2^{n},\qquad n\geq 1.
\end{gather*}
Thus, $(n\mapsto b_{n})\in\Theta(n^{\alvalue_{p}}\cdot 2^{n})$.
\end{Main Theorem}

Moreover, for $p=2$ we have a quite complete picture and we prove stronger results: we show that \autoref{Eq:Ansatz} holds and we determine $\hone(n)$, see \autoref{S:PisTwo} and \autoref{S:Fourier}.
We expect that similar methods can be used to handle the case $p>2$.

\begin{Remark}\label{R:IntroCharZero}
As usual in modular representation theory, the case $p=0$ behaves like 
$p=\infty$. Indeed, we have
\begin{gather*}
\alvalue_{\infty}:=\lim_{p\to\infty}\alvalue_{p}=-\frac{1}{2},
\end{gather*}
so we can compare \autoref{E:IntroCharZero} with \autoref{T:MainTheorem}.
\end{Remark}

\begin{Remark}\label{R:Conjecture}
In contrast to characteristic zero as in \autoref{E:IntroCharZero}, the limit
\begin{gather*}
\lim_{n\to\infty}\frac{b_{n}}{n^{\alvalue_{p}}\cdot 2^{n}}
\end{gather*}
does not exist, i.e. for $p>0$ there is no constant $c$ such that $b_{n}\sim c\cdot n^{\alvalue_{p}}\cdot 2^{n}$.
\end{Remark}

\begin{Example}\label{E:Hfactor}
For $p=2$ we have $b_{2n-1}=b_{2n}$ and the logplot of the even values of $b_{n}/(n^{-0.708}\cdot 2^{n})$ (splitting the sequence $b_{n}/(n^{-0.708}\cdot 2^{n})$ into even and odd makes it more regular, see for example \autoref{S:Monoton}) gives
\begin{gather*}
\begin{tikzpicture}[anchorbase]
\node at (0,0) {\includegraphics[height=5cm]{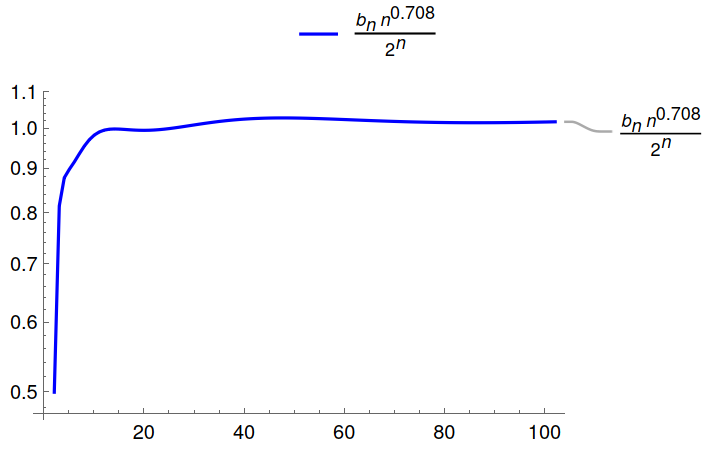}};
\node at (-2,1.7) {even $n$};
\end{tikzpicture}
\hspace{-0.5cm},
\begin{tikzpicture}[anchorbase]
\node at (0,0) {\includegraphics[height=5cm]{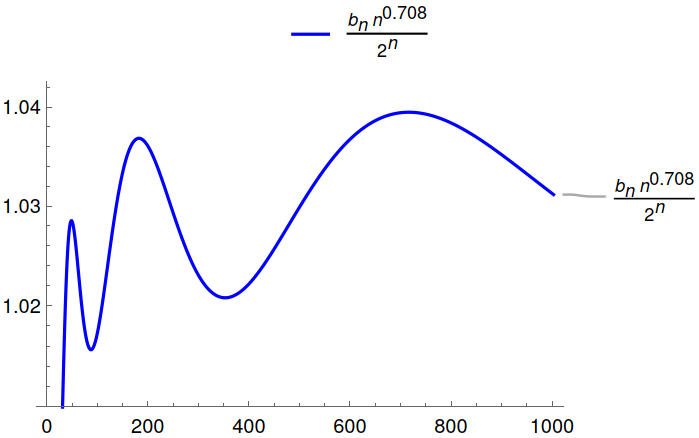}};
\node at (-2,1.7) {even $n$ zoom};
\end{tikzpicture}
.
\end{gather*}
These illustrations show the graph of $b_{n}/(n^{-0.708}\cdot 2^{n})$ for $n\in\{0,2,\dots,198,200\}$ and $n\in\{0,2,\dots,1998,2000\}$ (the fluctuations in the beginning vanish quickly) and we see the tiny oscillation, which is proven to be true in \autoref{S:Fourier}.
\end{Example}

\begin{Remark}\label{R:Thanks}
For $p=2$ \autoref{T:MainTheorem}, or rather the stronger version in \autoref{P:MainTheorem}, was independently proven 
in \cite{La-char2-story} using quite different methods. Moreover, and similar but different to \autoref{S:Fourier}, the paper \cite{La-char2-story} also identifies $h(n)$ from \autoref{Eq:Ansatz}. For $p=2$ \cite{La-char2-story} even implies some results for other representations as well. 
\end{Remark}

\begin{Remark}
Let us comment on the quantum case; the analog of \autoref{T:MainTheorem} 
for quantum $SL_{2}$. If the quantum parameter is not a root of unity, then 
the same discussion as in \autoref{E:IntroCharZero} works. The next case, where the quantum parameter is a root of unity and the underlying field is of characteristic zero can be deduced from \cite{LaPoReSo-growth-qgroup}. For the final possibility, the so-called mixed case as e.g. in \cite{SuTuWeZh-mixed-tilting}, we expect a result very similar to \autoref{T:MainTheorem}.
\end{Remark}

\subsection{Proof outline of \autoref{T:MainTheorem}}

It is a classical fact that all direct summands
of $V^{\otimes n}$ are \emph{tilting representations} over $\Gamma=SL(V)$. In the special case
of $\Gamma=SL_{2}(\mathbf{k})$ the characters of indecomposable tilting representations are explicitly known; using this we derive a recursion for the sequence $b_{n}$, see \autoref{SS:Recursion}. 

Then we translate this recursion into a functional equation for the generating function $F(z)$ of the sequence $b_{n}$, see \autoref{th2.4}. Then we analyze the behavior of the function $F(z)$ in the vicinity of its radius of convergence and deduce \autoref{T:MainTheorem} by employing suitable Tauberian theorems in \autoref{S:MainTheorem}. 

One of the steps is the following elementary looking inequality from \autoref{T:Mono}:
\begin{gather*}
b_{n+2}\leq 4b_{n},\; n\geq 0.
\end{gather*}

\subsection{Odds and ends}

Here are a few open questions that might be interesting to explore. (i) is proven for $p=2$, see \autoref{S:PisTwo} and \autoref{S:Fourier}, and similar methods apply in general.

\begin{enumerate}[label=(\roman*)]

\item In \autoref{T:Mono} we prove that
$b_{n+2}\leq 4b_{n}$.
A natural question is whether we have $\lim_{n\to\infty}b_{n+2}/b_{n}=4$.
Moreover, recall from \autoref{R:Conjecture} that there is no constant such that
$b_{n}\sim c\cdot n^{\alvalue_{p}}\cdot 2^{n}$.
However, one could conjecture that there exists a continuous real function $\htwo(x)$ of period one such that
\begin{gather}\label{Eq:AlmostConjecture}
	b_{n}\sim\htwo\big(\log_{2p}(n)\big)\cdot n^{\alvalue_{p}}
	\cdot 2^{n}.
\end{gather}
For $p=2$ this is proven in \autoref{P:MainTheorem}.

\item Recall from \autoref{T:MainTheorem} that 
$\frac{b_{n}}{n^{\alvalue_{p}}\cdot 2^{n}}$ is bounded 
by two constants, and it would be interesting to have good explicit values for these constants,
see \autoref{SS:ScalarsMainTheorem} for some possible values.
For a more precise asymptotic formula one would need to analyze the oscillation as in \autoref{E:Hfactor}.

\item One could try to compute other asymptotic formulas.
An example that comes to mind is to verify an analog of \autoref{T:MainTheorem} (with the same subexponential factor) 
for the three dimensional 
tilting $SL_{2}$-representation $\mathrm{Sym}^{2}V$ (in case $p\neq 2$). When $p=2$ this three dimensional $SL_{2}$-representation 
is not tilting and getting asymptotic formulas might be very hard.

\end{enumerate}


\noindent\textbf{Acknowledgments.} We thank
Michael Larsen for sharing a draft of \cite{La-char2-story},
which uses very interesting methods that are quite different from the ones in this work. The methods are exciting and important, so we all agreed that it makes sense to write two separate papers.
We also thank Henning Haahr Andersen, David He, Abel Lacabanne and Pedro Vaz for useful discussions and comments, and are very grateful to Kenichi Shimizu for pointing out a gap in the proof of \autoref{T:AppendixProjTensor} in the first version.

We express our appreciation to ChatGPT for their assistance with proofreading. In addition, D.T. extends heartfelt thanks to the tree constant for inspiration.

P.E.'s work was partially supported by the NSF grant DMS-2001318, D.T. was supported by the ARC Future Fellowship FT230100489.

\section{Fractal behavior of growth problems}

Before coming to the main parts of this paper, let us briefly 
indicate the main difficulty (and maybe the most exciting part)
of growth problems in the above sense: a certain type of fractal behavior 
of these sequences and of their generating functions.
For the sake of this paper, and partially justified by the discussion in this section, we use \emph{fractal behavior} to mean that the exponent $\avalue$ of the subexponential factor is
transcendental (note that e.g. $\alvalue_{p}$ from \autoref{Eq:Alpha} is transcendental), or at least irrational. For instance, $\alvalue_{p}$ could be the fractal dimension of some fractal.

\subsection{No fractal behavior}\label{SS:NoFractal}

For $\mathbf{k}$ of characteristic zero, let $\Gamma$ 
be a connected reductive algebraic group over $\mathbf{k}$.
\cite[Theorem 2.2]{Bi-asymptotic-lie} and \cite[Theorem 2.5]{CoEtOs-growth-mod-p} give 
the asymptotic for the numbers $b_{n}$ for an arbitrary finite dimensional $\Gamma$-representation.
The example of $\Gamma=SL_{2}(\mathbf{k})$ and its vector 
representation is given in \autoref{E:IntroCharZero} where $b_{n}\sim\sqrt{2/\pi}\cdot n^{-1/2}\cdot 2^{n}$.

In general, \cite[Theorem 2.2]{Bi-asymptotic-lie} and \cite[Theorem 2.5]{CoEtOs-growth-mod-p} prove (since this case is semisimple the same holds true for the length instead of the number of indecomposable summands):

\begin{Proposition}\label{P:ReductiveGroup}
In the above setting, the asymptotic takes the form
\begin{gather*}
b_{n}\sim
C\cdot n^{\avalue}\cdot(\dim_{\mathbf{k}}V)^{n}\text{ with }\avalue\in\tfrac{1}{2}\Z,C\in\R_{>0}.
\end{gather*}
Thus, the subexponential factor $n^{\avalue}$ always has some half integer exponent.\qed
\end{Proposition}

It is not a coincidence that $\avalue\in\tfrac{1}{2}\Z$ since: 
\begin{gather*}
\scalebox{0.89}{\text{``The nature of the generating function’s singularities determines 
the associated subexponential factor.''}}
\end{gather*}
The above strategy is well-known in symbolic dynamics, and under certain assumptions on the 
singularities of the generating 
function one always gets half integer powers for the subexponential growth term. 
For example, this works if the generating function is algebraic.

In fact,
there are well developed algorithms to compute the asymptotics if the generating function is 
sufficiently nice (e.g. meromorphic), 
see for example \cite[Section 7.7]{Mi-analytic-combinatorics}. 
The algorithm presented therein always produces 
an exponent $s\in\tfrac{1}{2}\Z$ and is enough to treat the case of a connected reductive algebraic group in
characteristic zero.

Here is another example where one always gets a half integer coefficient, namely zero,
regardless of the characteristic. As we will argue below, this 
means that these cases do not show fractal behavior with respect to 
the growth problems we consider.

Let $\Gamma$ be a finite group and 
let $V$ be a representation of $\Gamma$ over $\mathbf{k}$, with $\mathbf{k}$ of arbitrary characteristic.
We assume for simplicity that $V$ is a faithful $\Gamma$-representation; if $V$ is not faithful we replace 
$\Gamma$ by a quotient and continue as below.
Let $M\subset\Gamma$ be the
central cyclic subgroup of $\Gamma$ 
containing all scalar matrices. Let $m=|M|$. The category
$\mathrm{Rep}(\Gamma)$ of finite dimensional $\Gamma$-representations splits into a direct sum
\begin{gather*}
\mathrm{Rep}(\Gamma)=\bigoplus_{i\in \Z/m\Z}\mathrm{Rep}(\Gamma)_{i},
\end{gather*}
where $M$ acts on objects of $\mathrm{Rep}(\Gamma)_{i}$ in the same way as in $V^{\otimes i}$.

Let $\mathrm{S}^{i}(\Gamma)$ be the set of isomorphism 
classes of simple objects of $\mathrm{Rep}(\Gamma)_{i}$; 
for any $L\in\mathrm{S}^{i}(\Gamma)$ let 
$P(L)\in\mathrm{Rep}(\Gamma)_{i}$ be its
projective cover. For any $i\in\Z/m\Z$ we define
\begin{gather*}
\ell_{i}=\frac{m}{|\Gamma|}
\bigg(\sum_{L\in\mathrm{S}^{i}(\Gamma)}\dim_{\mathbf{k}}P(L)\bigg),
\quad
\nu_{i}=\frac{m}{|\Gamma|}
\bigg(\sum_{L\in\mathrm{S}^{i}(\Gamma)}\dim_{\mathbf{k}}L\bigg)
.
\end{gather*}
Let $A$ be the regular representation of $\Gamma$ and let $A=\oplus_{i\in \Z/m\Z}A_{i}$ be its
decomposition into summands $A_{i}\in\mathrm{Rep}(\Gamma)_{i}$.
A result of Bryant--Kov{\'a}cs \cite[Theorem 2]{BrKo-tensor-group} says that for 
sufficiently large $n$ with $n\equiv i\bmod{m}$
we have that $A_{i}$ is a direct summand of $V^{\otimes n}$. 
Similarly to $b_{n}$, let $l_{n}=l_{n}(V):=\ell(V^{\otimes n})$, where $\ell(\placeholder)$ denotes the \emph{length} of representations.
Assuming $\mathbf{k}=\overline{\mathbf{k}}$, 
it follows easily
that for $n\equiv i\bmod{m}$ we have
\begin{gather*}
l_{n}\sim\ell_{i}\cdot(\dim_{\mathbf{k}}V)^{n},
\quad
b_{n}\sim\nu_{i}\cdot(\dim_{\mathbf{k}}V)^{n}.
\end{gather*}
In particular, each of the sequences $l_{n}/(\dim_{\mathbf{k}}V)^{n}$ and $b_{n}/(\dim_{\mathbf{k}}V)^{n}$ has at most $m$ limit points. In particular, \cite[Section 7.7]{Mi-analytic-combinatorics} applies and the subexponential factor has half integer exponent.
A similar analysis works without the assumption 
$\mathbf{k}=\overline{\mathbf{k}}$.

This immediately proves that $\hone(n)$, in this case, is periodic 
with period $m$, the subexponential factor is trivial and the exponential factor is given by the dimension. Precisely:

\begin{Proposition}\label{P:FiniteGroups}
For a finite group $\Gamma$ and a $\Gamma$-representation $V$ we always 
have
\begin{gather*}
\begin{aligned}
l_{nm+r}&\sim
s(r)\cdot n^{0}\cdot(\dim_{\mathbf{k}}V)^{nm}
,
\\
b_{nm+r}&\sim
t(r)\cdot n^{0}\cdot(\dim_{\mathbf{k}}V)^{nm},
\end{aligned}
\quad\text{with }s(r),t(r)\in(0,1],t(r)\leq s(r),
\end{gather*}
for some $m\in\N$ and all $r\in\{0,\dots,m-1\}$.\qed
\end{Proposition}

The scalars $s(r),t(r)$ in \autoref{P:FiniteGroups} are easy to compute. For $t(r)$ see for example \cite[Proposition 2.1]{CoEtOs-growth-mod-p} and \cite[(2A.1)]{LaTuVa-growth-pfdim} for characteristic zero, and \cite[Section 5]{LaTuVa-growth-pfdim-two} for arbitrary characteristic.

\begin{Example}
For instance if $p=2$, and $\Gamma$ is the symmetric group $S_{3}$ or $S_{4}$ and $V$ is a faithful representation of $\Gamma$ we have
$l_{n}\sim\frac{2}{3}\cdot\dim_{\mathbf{k}}V^{n}$
and $b_{n}\sim\frac{1}{2}\cdot\dim_{\mathbf{k}}V^{n}$.
\end{Example}

\begin{Remark}\label{R:FiniteTensor}
Let $\mathbf{C}$ be a finite tensor category, see e.g. \cite[Chapter 6]{EtGeNiOs-tensor-categories} for details.
Let $X$ be a tensor generator of $\mathbf{C}$.
We have a decomposition of $\mathbf{C}$ over the universal grading group $U$ of $\mathbf{C}$. Write $X\cong\oplus_{i\in I}X_{i}$ for $X_{i}$ indecomposables with $X_{i}$ of degree $d_{i}\in U$. Let $U_{X}$ be the subgroup of $U$ generated by $d_{i}-d_{j}$. Then $U/U_{X}$ is a cyclic group $\Z/m\Z$, and statements similar to 
\autoref{P:FiniteGroups} hold. Details are omitted, but the key 
statement hereby is proven in \autoref{T:AppendixProjTensor} in the appendix.
\end{Remark}

To summarize, in the two above settings the subexponential
exponent $\avalue$ is a \emph{half-integer}, in particular not transcendental, and the function $\hone(n)$ is \emph{constant up to some period} $m\in\N$.

\subsection{Fractal behavior}

For $\Gamma=SL_{2}(\mathbf{k}^{2})$ our problem, where $\mathbf{k}$ is of prime characteristic, 
is difficult also because the generating function 
that we compute in \autoref{th2.4} does not have nice enough properties 
to run the classical strategies. For example, we will see that we 
have to face a dense set of singularities, see e.g. 
\autoref{P:SingularitiesPGeneral}.
And in fact, the exponent we get is not a half integer but rather the transcendental number $\alvalue_{p}$
from \autoref{Eq:Alpha}.

In a bit more details, see \autoref{ThmFunEq} for the precise statement,
we get a functional equation for the generating function $F$ of the numbers $b_{n}$
that takes the form
\begin{gather}\label{Eq:MahlerFirst}
F(w)=r_{1}(w)+r_{2}(w)\cdot F(w^{p}).
\end{gather}
Functions of this type are called \emph{2-Mahler functions of degree $p$}.

\begin{Remark}
The name originates in Mahler's approach to 
transcendence and algebraic independence results for
the values at algebraic points in the study of
power series satisfying functional equations of a
certain type. Mahler's original functional equation is of the form $F(w)=w+F(w^{2})$, which is 
an example of what we call a \emph{2-Mahler function of degree 2}. 

Let $r_{i}(w)$ denote rational functions.
The Mahler functions above have been 
generalized under the umbrella of \emph{$s$-Mahler functions of degree $p$} satisfying
\begin{gather}\label{Eq:MahlerFunction}
r_{0}(w)\cdot F(w)=r_{1}(w)+r_{2}(w)\cdot F(w^{p})+r_{3}(w)\cdot F(w^{p^{2}})+\dots+
r_{s}(w)\cdot F(w^{p^{s-1}}),
\end{gather}
but we will not need this generalization. 
We refer to \cite{Ni-mahler} for a nice discussion of Mahler
functions and a list of historical references. Mahler functions occur 
in combinatorics as generating functions of partitions
and related structures. 
\end{Remark}

A crucial fact is that such a Mahler function often grows with exponent $\avalue+1$ for a transcendental $\avalue$ when approaching its relevant singularity. We will use this 
in \autoref{S:Asym}.

Moreover, connections from Mahler functions that grow with 
transcendental $\avalue$ to fractals are well-studied.
To the best of our knowledge, general
theorems relating them are not known, and connections are 
only example-based. In the rest of this section, we give
several examples that are easier to deal with than our main result.

The common source for these examples is the well-known principle that projecting $p$-adic objects onto the real world leads to fractals. At the same time, modular representations of algebraic groups are known to exhibit $p$-adic patterns, stemming from Steinberg-type tensor product theorems. Therefore, one can expect that combinatorial invariants of modular representations of algebraic groups (such as dimensions, weight multiplicities, etc.) tend to exhibit fractal behavior. 
In fact, fractal behavior in the study of reductive groups in prime characteristic 
has been observed several times, and is part of folk knowledge.

Most relevant for this paper are fractal patterns within the study of tilting representations.
For example, for an algebraically closed field $\mathbf{k}$ of prime characteristic $p$,
let $\Gamma=SL_{2}(\mathbf{k})$.
The multiplicities of Weyl in tilting $\Gamma$-representations have fractal patterns of step size $p$ and so do their characters, 
by Donkin's tensor product formula \cite[Proposition 2.1]{Do-tilting-alg-groups} (this is a bit easier to see in the reformulation given in \cite[Section 3]{SuTuWeZh-mixed-tilting}).
The same is reflected in the Temperley--Lieb combinatorics, see e.g. 
\cite[Theorem 2.8]{BuLiSe-tl-char-p} 
or \cite[Figure 3]{Sp-modular-tl}, which is exploited in \cite{KhSiTu-monoidal-cryptography}.

For higher rank there are also several known instances of fractal behavior of tilting 
characters, see for example the picture of cells in \cite[Figure 1]{An-p-cells-affine-weyl} or 
billiards of tilting characters as in \cite[Section 4]{LuWi-billiard-conjecture}.

\subsubsection{Modular representations and the Cantor set}\label{SS:Cantor}

One of the simplest examples of a fractal is the 
(bounded) \emph{Cantor set} $\mathtt{C}$, which can be realized as the set of real numbers in $[0,1]$ which admit a ternary expansion without digit $1$ with the probably familiar picture
\begin{gather*}
\begin{tikzpicture}[anchorbase]
\node at (0,0) {\includegraphics[height=2cm]{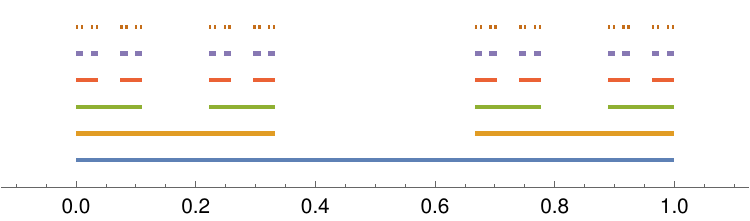}};
\end{tikzpicture}
.
\end{gather*}
This example is related to representation theory of $\Gamma=SL_{2}(\mathbf{k})$ for $\mathbf{k}$ being algebraically closed 
of characteristic $p$ (the standard choice is $\mathbf{k}=\bar{\F}_{p}$) as follows.

The finite dimensional simple $\Gamma$-representations 
$L_{n}$ are indexed by $n\in\N$, their highest weight.
Expressing the number $p$-adically, say $n=[n_{r},
\dots,n_{1},n_{0}]=n_{r}p^{r}+\dots+n_{1}p+n_{0}$ 
with $n_{i}\in\{0,\dots,p-1\}$, one gets
\begin{gather*}
L_{n}
\cong
L_{n_{r}}^{(r)}\otimes\dots\otimes L_{n_{1}}^{(1)}\otimes L_{n_{0}},
\end{gather*}
where the exponent ${}^{(i)}$ 
denotes Frobenius twists which 
acts on characters by $x\mapsto x^{p^{i}}$. This is a special case 
of Steinberg's tensor product theorem.
Since
the character of $L_{n}$ for $n=[n_{0}]$ 
is $\frac{x^{n}-x^{-n}}{x-x^{-1}}$, this gives a description of the weights of $L_{n}$ for $n\in\N$.

This also shows that we may consider, often infinite dimensional, 
representations $L_{n}$ where $n\in\Z_{p}$ 
is a $p$-adic integer, which is made explicit in Haboush's generalized Steinberg tensor product theorem for distribution algebras \cite[Theorem~4.9]{Haboush}. Concretely, given 
$n\in\Z_{p}$, 
we may define 
$L_{n}$ to be the infinite tensor product
\begin{gather*}
L_{n}:=\bigotimes_{j=0}^{\infty}L_{n_{j}}^{(j)},\quad\mbox{with}\quad n=[\dots,n_{r},\dots,n_{1},n_{0}]=\sum_{j=0}^{\infty}
n_{j}p^{j},
\end{gather*}
which, by definition, is the span of tensor products 
of vectors in $L_{n_{j}}^{(j)}$ such that 
almost all of them are (fixed) highest weight vectors. 

The group $\Gamma$ does not act in this space, but $L_{n}$ admits an action of the distribution algebra $\mathrm{Dist}=\mathrm{Dist}(\Gamma)$, for which this representation is generated by the highest weight vector $v_{n}$. This gives $L_{n}$ a $\N$-grading with even degrees, placing vectors of weight $\mu$ in degree $n-\mu$. The Hilbert series for this grading is thus
\begin{gather*}
h_{L_{n}}(w)=\prod_{j=0}^{\infty}
(1+w^{2\cdot p^{j}}+\dots+w^{2n_{j}\cdot p^{j}}).
\end{gather*}
This is a holomorphic function for $|w|<1$, and we will mostly consider it on the interval $(0,1)$, i.e. as $h_{L_{n}}\colon (0,1)\to\R$. 

Consider now the special case $p=3$, $n=-\frac{1}{2}$. Hence, 
$n_{j}=1$ for all $j$ and we get 
\begin{gather*}
h(w):=h_{L_{-1/2}}(w)=\prod_{j=0}^{\infty}(1+w^{2\cdot 3^{j}}).
\end{gather*}
This function satisfies the functional equation 
\begin{gather*}
h(w)
=(1+w^{2})\cdot h(w^{3}).
\end{gather*}
This is a 2-Mahler function of degree $3$, and the 
relevant singularity is at $w=1$.

The Taylor coefficients of
$h(w)$ form the \emph{Cantor set sequence} 
$(\mathrm{ca}_{n})_{n\in\N}$ defined by
\begin{gather*}
\mathrm{ca}_{n}
=\begin{cases}
1 & \text{if the ternary expansion of $n$ contains no $1$},
\\
0 & \text{otherwise}.
\end{cases}
\end{gather*}
This sequence is well-studied, see e.g. \cite[Section 1]{MR4216102} and \cite[A292686]{oeis}, and $h(w)$ is its generating function.

The corresponding fractal is constructed as follows. Let $D(L)$ 
denote the set of degrees of~$L$, so that $D(L_{-1/2})$ is the set of nonnegative integers with ternary 
representations where all digits of $m$ 
take values $0,2$ only.
For $N\in\N$, define $\Delta_{N}:=\frac{1}{3^{N}}D(L_{-1/2})$.
We then have $\Delta_{N}\subset\Delta_{N+1}$. 
Let $\Delta_{\infty}:=\cup_{N\in\N}\Delta_{N}$. 
Then the closure of $\Delta_{\infty}$ in $\R$ is the \emph{unbounded} Cantor 
set $\mathtt{C}_{\infty}$, which is invariant under multiplication by $3$. 
The set $\mathtt{C}$ is the intersection 
$\mathtt{C}_{\infty}\cap [0,1]$, and 
$\mathtt{C}_{\infty}=\cup_{N\in\N}3^{N}\mathtt{C}$ (nested union).
Note that the set $\mathtt{C}_{\infty}$ comes with a natural measure $\mu$, 
the Hausdorff measure of $\mathtt{C}_{\infty}$: the weak limit of the counting measures of $\Delta_{N}$ rescaled by dividing by $|\Delta_{N}\cap [0,1]|=2^{N}$.

Now, by the standard theory of Mahler functions, see e.g. 
\cite[Theorem 1]{MR3589306}, we get
\begin{gather*}
C_{1}(1-w)^{-\avalue}\big(1+o(1)\big)\leq h(w)\leq C_{2}(1-w)^{-\avalue}\big(1+o(1)\big),\ x\uparrow 1, 
\end{gather*}
for some $0<C_{1}<C_{2}$ with $C_{1},C_{2}\in\R$. Moreover, $\avalue$ is found by substituting $(1-w)^{-\avalue}$ into the Mahler equation, i.e. from the condition 
$(1-w)^{-\avalue}\sim (1+w^{2})(1-w^{3})^{-\avalue}$ as $w\uparrow 1$. This yields 
$1=2\cdot 3^{-\avalue}$, hence 
\begin{gather*}
\avalue=\log_{3}2\approx 0.631.
\end{gather*}
The number $\avalue$ is transcendental, and the Hausdorff (and Minkowski) dimension of the Cantor sets $\mathtt{C}$ and $\mathtt{C}_{\infty}$. 

Rewriting the above using $\asymp$, see \autoref{Eq:Asymp} for the notation,
we get
\begin{gather*}
h(w)=\sum_{n\in\N}\mathrm{ca}_{n}w^{n}
\asymp_{w\uparrow 1}(1-w)^{-\avalue}\big(1+o(1)\big).
\end{gather*}
Then Tauberian theory (as for example in \autoref{L:Tauber}) implies
\begin{gather*}
\sum_{k=0}^{n}\mathrm{ca}_{k}\asymp n^{\avalue}\big(1+o(1)\big).
\end{gather*}
This describes the asymptotics 
of Ces{\'a}ro sums of the Cantor set sequence, which 
counts the dimension of the subspace $L_{-1/2}[\leq n]$ spanned by vectors in $L_{-1/2}$ of degree $\leq n$.

One may further ask how the sequence $C_{n}=n^{-\avalue}\sum_{k=0}^{n}\mathrm{ca}_{k}$ 
behaves when $n\to\infty$ or, related, how the function $(1-x)^{\avalue}h(w)$ (equivalently, 
$\ln(w^{-1})^{\avalue}h(w)$) behaves when $w\uparrow 1$. This behavior can be analyzed as follows. Let $\theta(m)=1$ for $m<0$ and $\theta(m)=0$ for $m\geq 0$.
We have 
\begin{gather*}
\lim_{k\to\infty}\ln(w^{-3^{-k}})^{\avalue}h(w^{3^{-k}})=\ln(w^{-1})^{\avalue}\lim_{k\to\infty}2^{-k}\prod_{m=-k}^{\infty}
(1+w^{2\cdot 3^m})=h_{0}(w),
\end{gather*}
for the function 
\begin{gather*}
h_{0}(w):=\ln(w^{-1})^{\avalue}\prod_{m=-\infty}^{\infty}
\frac{1+w^{2\cdot 3^m}}{2^{\theta(m)}}.
\end{gather*}
This is a periodic function in the sense that $h_{0}(w)=h_{0}(w^{3})$ (note that the product is absolutely convergent). This implies that the function $\ln(w^{-1})^{\avalue}h(w)$ has no limit as $w\uparrow 1$ and asymptotically has oscillatory behavior: it approaches the periodic function $h_{0}(w)$. 

Plotting this illustrates the overall growth rate and the oscillation:
\begin{gather}\label{Eq:Devil}
\begin{tikzpicture}[anchorbase]
\node at (0,0) {\includegraphics[height=4cm]{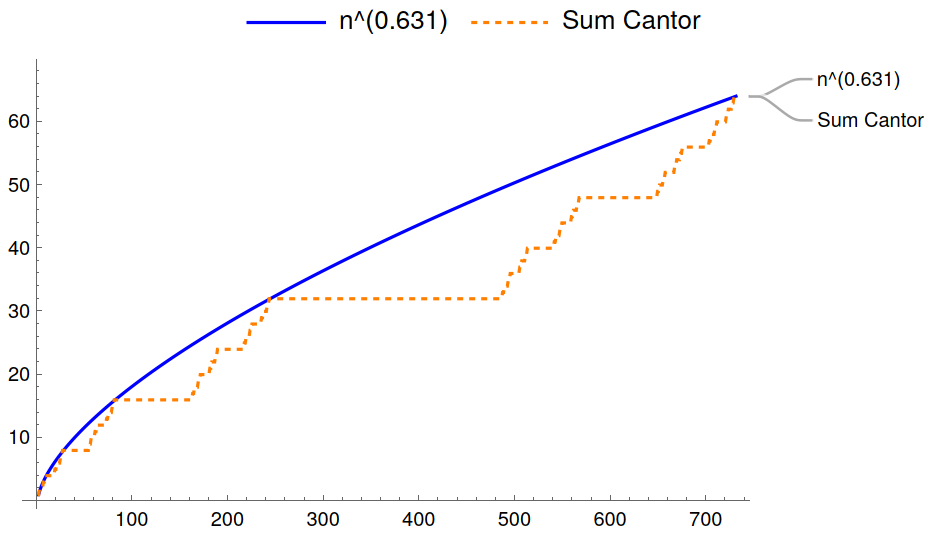}};
\node at (-1.5,1) {$p=3$};
\end{tikzpicture}
\hspace{-0.45cm},
\begin{tikzpicture}[anchorbase]
\node at (0,0) {\includegraphics[height=4cm]{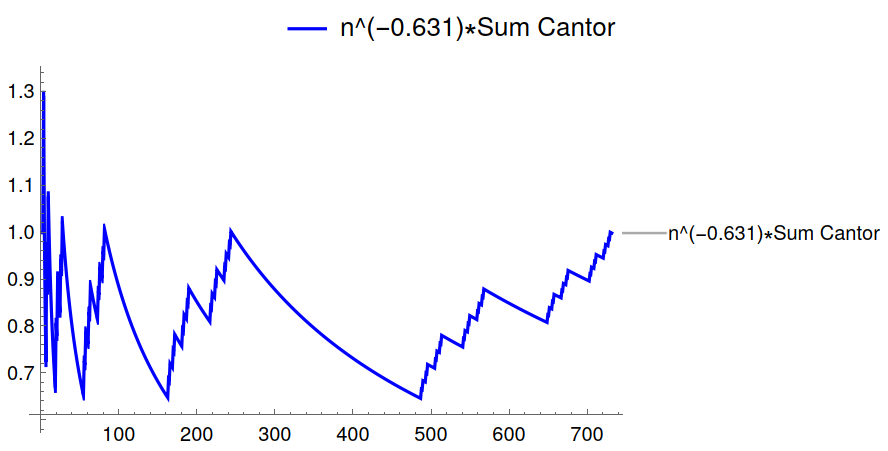}};
\node at (-1.5,1) {$p=3$};
\end{tikzpicture}
.
\end{gather}
To analyze this a bit further, it is easy to see that the analytic function $\mathtt{h}_{0}(v):=v^{-\avalue}h_{0}(e^{-v})$ in the region $\mathrm{Re}\,v>0$, which satisfies the equation $2\mathtt{h}_{0}(3v)=\mathtt{h}_{0}(v)$, is nothing but the Laplace transform 
of the Hausdorff measure $\mu$ of the set $\mathtt{C}_{\infty}$:
\begin{gather*}
\mathtt{h}_{0}(v)=\int_{0}^{\infty}e^{-tv}d\mu(t).
\end{gather*}
This implies that
\begin{gather*}
h_{0}(w)=\ln(w^{-1})^{\avalue}\int_{0}^{\infty}w^{t}d\mu(t).
\end{gather*}
Similarly, the sequence $C_{n}$ for large $n$ behaves as 
the \emph{devil's staircase function}:  
\begin{gather*}
\lim_{m\to\infty}C_{\lfloor 3^{m}w\rfloor}=w^{-\avalue}\mu([0,w]),
\end{gather*}
which is periodic under the map $w\mapsto 3w$. This is visualized in
\autoref{Eq:Devil}.
Hence, this sequence does not have a limit as $n\to\infty$ and exhibits oscillatory behavior, approaching the periodic function $\log_{3}(n)^{-\avalue}\mu\big([0,\log_3(n)]\big)$. 
Note that this function is continuous but not differentiable, nor
absolutely continuous, since it is the integral of a singular measure
(the Hausdorff measure of $\mathtt{C}_{\infty}$). It is, however, 
H{\"o}lder continuous with exponent $\avalue$. 
We can also find the range of oscillation of $C_{n}$. Indeed, it is easy to see that 
$\limsup_{n\to\infty}C_{n}=1$ (attained for $n=3^{k}$) while $\liminf_{n\to\infty}C_{n}=2^{-\avalue}\approx 0.646$ (approached for $n=2\cdot 3^{k}-1$). 

\begin{Remark} 
The function $h_{0}(e^{-3{^u}})$ is holomorphic in the strip 
$|\mathrm{Im}\,u|<\frac{\pi}{2\ln 3}$ and periodic under $u\mapsto u+1$, 
so we may consider its Fourier coefficients $A_{n}$. They are related to the Fourier coefficients 
$a_{n}$ of the measure $d\mu(e^{-3^{u}})$ by the formula 
\begin{gather*}
A_{n}=\Gamma(\avalue+\tfrac{2\pi in}{\ln 3})a_{n}.
\end{gather*}
Here and throughout $\Gamma(c)$ denotes the \emph{gamma function} evaluated at $c\in\C$ (not to be confused with the group $\Gamma$).
Since $\mu$ has unit volume on the period, we have $|a_{n}|\in O(1)$, $n\to\infty$. Hence:
\begin{gather*}
|A_{n}|\in O\big(|\Gamma(\avalue+\tfrac{2\pi in}{\ln 3})|\big)=O(n^{\avalue-\frac{1}{2}}e^{-\frac{\pi^{2} n}{\ln 3}}),\ n\to\infty.
\\
\begin{tikzpicture}[anchorbase]
\node at (0,0) {\includegraphics[height=4cm]{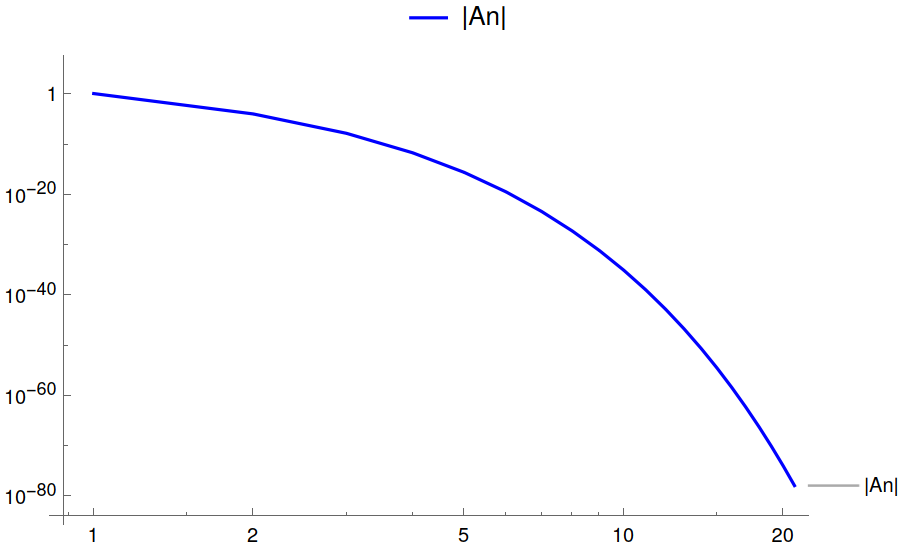}};
\node at (-1.5,0.7) {$p=3$};
\end{tikzpicture}
.
\end{gather*}
The numbers $|A_{n}|$ are tiny for large $n$, see the above loglogplot.
\end{Remark} 

\subsubsection{Generalization to general primes and highest weights} 

A similar analysis applies for any prime $p$ and rational highest weights $\lambda\in\Z_{p}\cap\Q$ (excluding positive integers, for which $L_{\lambda}$ is finite dimensional). Namely, in this case $\lambda$ has an infinite periodic expansion, which can be computed in the standard way: 
write $\lambda$ as $q-\frac{m}{n}$ where $0<\frac{m}{n}\leq 1$ and $q\in\Z$, and 
pick the minimal $N\in\Z_{>0}$ such that $n$ divides $p^{N}-1$. Then the repeating string of the expansion 
of $\lambda$ is the base $p$ expression of the number $r(\lambda):=\frac{m(p^N-1)}{n}$ (a string of length $N$, with zeros at the beginning if needed). If $r(\lambda)=[r_{N-1},\dots,r_{0}]$, where $r_{j}$ are its base $p$ digits, then 
\begin{gather*}
\avalue(\lambda)=\frac{1}{N}\sum_{j=0}^{N-1}\log_{p}(r_{j}+1).
\end{gather*}
The corresponding fractal $\mathtt{C}(\lambda)$ is the set of real numbers which have a base $p$ expansion with $j$th digit in $\{0,\dots,r_{j\text{ mod }N}\}$, which has Hausdorff (and Minkowski) dimension $\avalue(\lambda)$. Note that for $\lambda\in\Z_{<0}$, $\mathtt{C}(\lambda)=\R_{\geq 0}$ (so $\avalue(\lambda)=1$), while for $\lambda\notin\Z$ we get an actual fractal, i.e. $\avalue(\lambda)<1$ is transcendental. For example, for $p=3$ we have $\mathtt{C}(-1/2)=\mathtt{C}_{\infty}$ and $\avalue(-1/2)=\log_{3}2$ as explained above. 

The deeper analysis of the oscillating behavior of the character of $L_{\lambda}$
and the sequence of Ces{\'a}ro sums of its coefficients (i.e., the sequence $\dim_{\mathbf{k}}L_{\lambda}[\leq n]$ of dimensions of weight spaces $\leq n$) also extends mutatis mutandis to general $p$ and $\lambda$. 
In fact, the power behavior of the sequence $\dim_{\mathbf{k}}L_{\lambda}[\leq n]$ is observed for 
sufficiently generic $\lambda\in\Z_{p}$. Namely, if $\lambda=\sum_{j=0}^{\infty}r_{j}p^{j}$
and there exists a limit 
\begin{gather*}
\avalue(\lambda):=\lim_{N\to \infty}\frac{1}{N}\sum_{j=0}^{N-1}\log_{p}(r_{j}+1),
\end{gather*}
then we have 
\begin{gather*}
\lim_{n\to\infty}\frac{\log\dim_{\mathbf{k}}L_{\lambda}[\leq n]}{\ln n}=\avalue(\lambda).
\end{gather*}
This includes rational $\lambda$, and also if $\lambda$ is chosen randomly 
(all digits are independent and uniformly distributed) then
\begin{gather*}
\avalue(\lambda)=\frac{1}{p}\log_{p}(p!).
\end{gather*}
For example, for $p\neq 2$, let $\avalue=\log_{p}(p-1)$. Consider the \emph{base $p$ Cantor set sequence}:
\begin{gather*}
\mathrm{ca}_{n}^{p}
=\begin{cases}
1 & \text{if the base $p$ expansion of $n$ contains no $1$},
\\
0 & \text{otherwise}.
\end{cases}
\end{gather*}
Thus, we have:
\begin{gather*}
p=5\colon
\avalue\approx0.861,\quad
p=7\colon
\avalue\approx0.921,\quad
p=11\colon
\avalue\approx0.960.
\end{gather*}
All of these values are transcendental.

Then the analog of \autoref{Eq:Devil} for $p=5$ and $p=7$ is:
\begin{gather*}
\begin{tikzpicture}[anchorbase]
\node at (0,0) {\includegraphics[height=4cm]{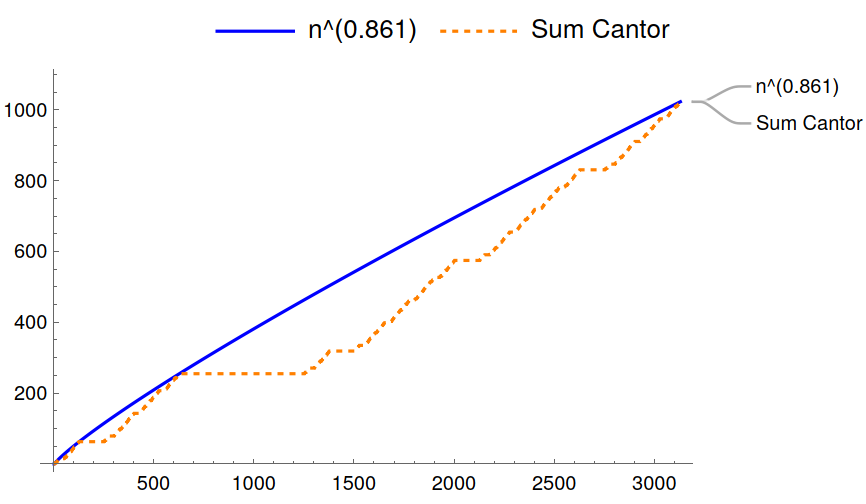}};
\node at (-1.5,1) {$p=5$};
\end{tikzpicture},
\begin{tikzpicture}[anchorbase]
\node at (0,0) {\includegraphics[height=4cm]{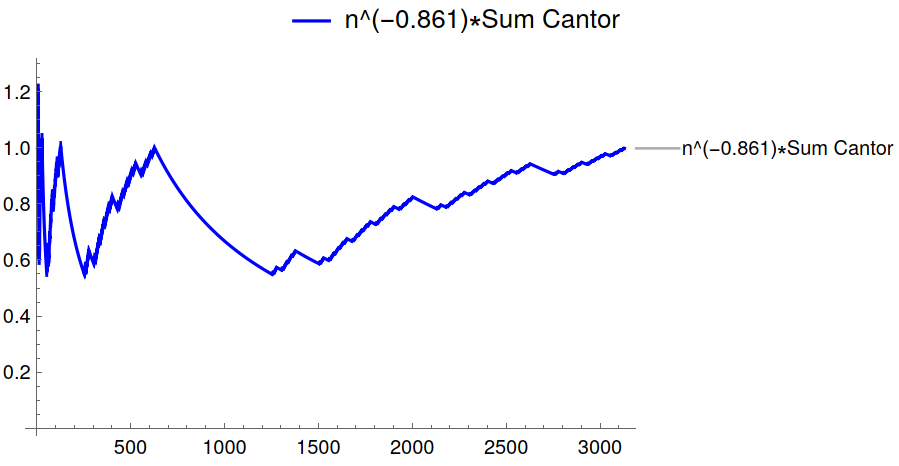}};
\node at (-1.5,1) {$p=5$};
\end{tikzpicture}
,
\\
\begin{tikzpicture}[anchorbase]
\node at (0,0) {\includegraphics[height=4cm]{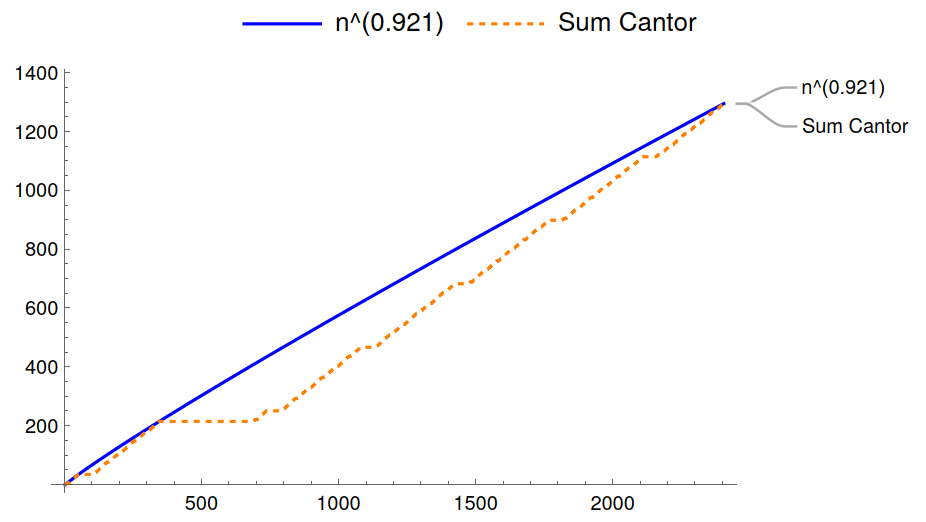}};
\node at (-1.5,1) {$p=7$};
\end{tikzpicture},
\begin{tikzpicture}[anchorbase]
\node at (0,0) {\includegraphics[height=4cm]{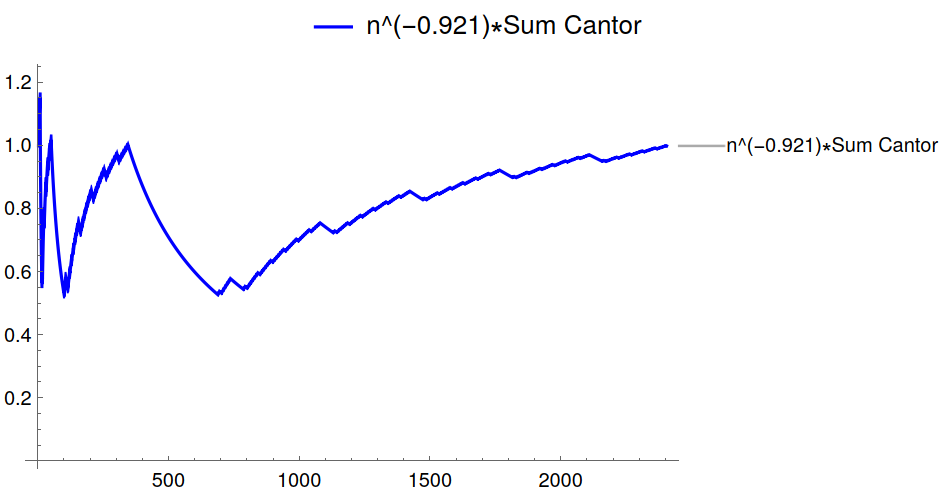}};
\node at (-1.5,1) {$p=7$};
\end{tikzpicture}
.
\end{gather*}

\subsubsection{Generalization to higher rank} 

The above analysis can also be extended to simply connected simple algebraic groups $G$ of arbitrary rank $r$. In this case, a weight $\lambda$ is an $r$-tuple $(\lambda_{1},\dots,\lambda_{r})$, $\lambda_{j}\in\Z_{p}$ (the coefficients of $\lambda$ 
with respect to fundamental weights), and a fractal attached to the simple highest weight module $L_{\lambda}$ over the distribution algebra $\mathrm{Dist}(\Gamma)$ can be built in a Euclidean space $\mathrm{Fun}(R_{+},\R)$, where $R_{+}$ is the set of positive roots of $\Gamma$, for which we choose an ordering. For example, suppose $\lambda=\frac{\mu}{1-p}$, where $\mu=(\mu_{1},\dots,\mu_{r})$ 
is an integral weight with $0\leq\mu_{i}\leq p-1$. Then the fractal 
$\mathtt{C}(\lambda,B)$ attached to $L_{\lambda}$ depends on a choice of a PBW basis $B$ 
of $L_{\mu}$. 

Namely, pick a collection $B\subset\mathrm{Fun}(R_{+},[0,p-1])$ such that the set of vectors $\prod_{\avalue\in R_{+}}e_{\avalue}^{b(\avalue)}v_{\mu}$, 
$b\in B$ (product in the chosen order) forms a basis of $L_{\mu}$ (where $v_{\mu}\in L_{\mu}$ is the highest weight vector). Then we define $\mathtt{C}(\lambda,B)$ as the set of functions 
$\phi\colon R_{+}\to\R_{\geq 0}$ which have a base $p$ expansion such that for all $j\in\Z$, the $j$th digit $\phi_{j}$ of $\phi$ belongs to $B$. The Hausdorff (and Minkowski) 
dimension of $\mathtt{C}(\lambda,B)$ equals
\begin{gather*}
\avalue(\lambda)=\log_{p}\dim_{\mathbf{k}}L_{\mu}
\end{gather*}
for any choice of $B$. (However, it must be mentioned that $\dim_{\mathbf{k}}L_{\mu}$ is notoriously hard to compute in general.)

Also, $\mathtt{C}(\lambda,B)$ comes equipped with a natural measure $\mu_{B}$, obtained by 
suitable rescaling of the counting measures as before (a multiple of the Hausdorff measure). Finally, we have a proper map $\pi\colon\mathtt{C}(\lambda,B)\to\R_{\geq 0}^{r}=\R_{\geq 0}R_{+}$ 
given by $\pi(\phi)=\sum_{\avalue\in R_{+}}\phi(\avalue)\avalue$, 
and the measure $\mu:=\pi_{\ast}\mu_{B}$ does not depend on $B$, nor on the ordering of positive roots (it expresses the large-scale asymptotics of the character of $L_{\lambda}$). This measure completely determines the large-scale behavior of $L_{\lambda}$. 

\subsubsection{Modular representations and Sierpinski gaskets}

Fix a prime $p$.
To give another example, recall \emph{Sierpinski's gasket or triangle}, which is Pascal's triangle modulo $p$. For $p=2$ we have 
(we only illustrate cutoffs):
\begin{gather*}
\begin{tikzpicture}[anchorbase]
\node at (0,0) {\includegraphics[height=5cm]{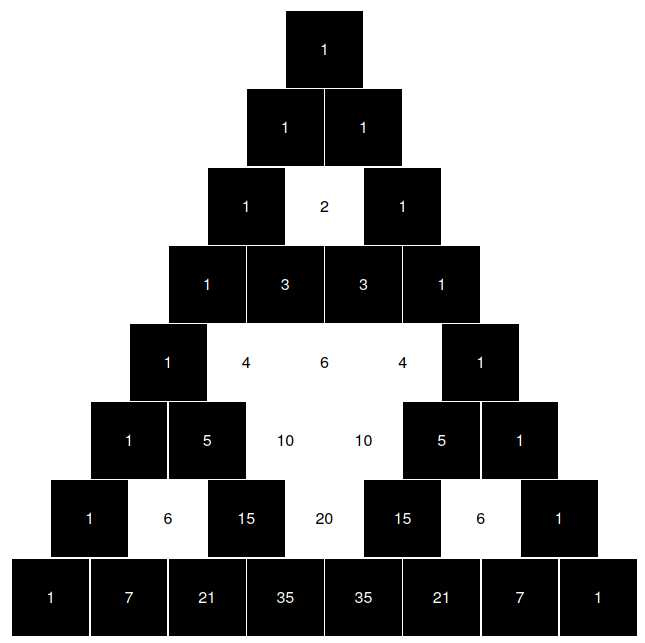}};
\node at (-1.5,1) {$p=2$};
\end{tikzpicture}
,\quad
\begin{tikzpicture}[anchorbase]
\node at (0,0) {\includegraphics[height=5cm]{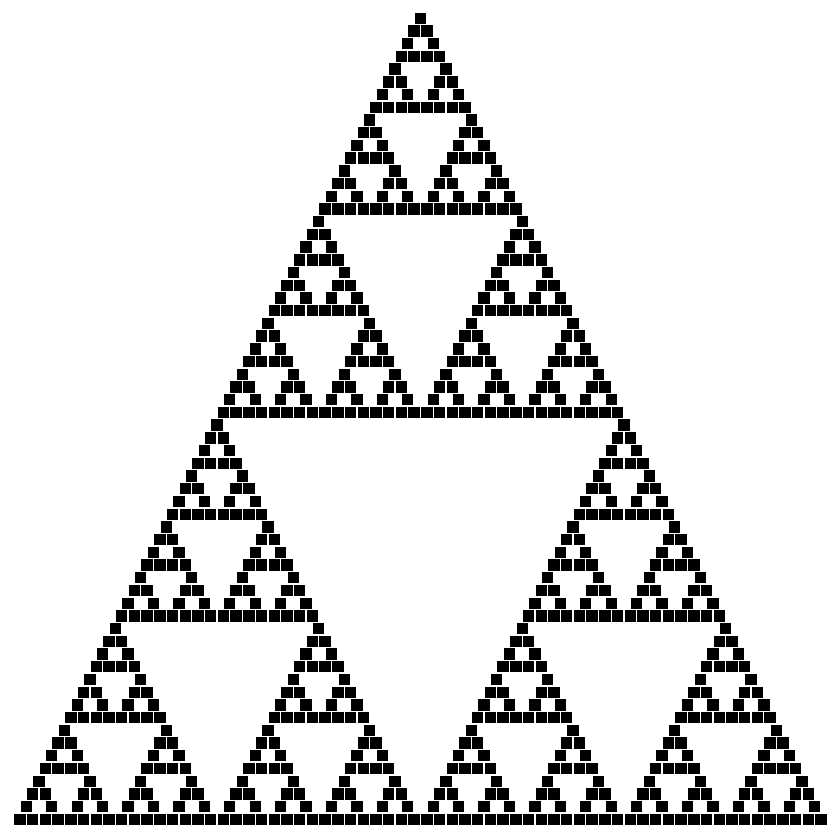}};
\node at (-1.5,1) {$p=2$};
\end{tikzpicture}
.
\end{gather*}
Sierpinski's triangle is combinatorially 
given by the limit of the above pictures, keeping the black boxes and 
disregarding the white boxes.

This example is related to representation theory of $\Gamma=SL_2(\mathbf{k})$ as follows.

For each $N$ let us plot the weights of $L_{n}$ for $n\leq p^{N}-1$ 
on the line $x+y=n$ on the coordinate plane: the weight $n-2k$ corresponds to the point $(n-k,k)$. 
Let us rescale the obtained set by the factor $p^{-N}$ and denote the resulting set by $\Delta_{N}(p)$. Thus $\Delta_{N}(p)$ is the set of pairs $(x,y)$, $x=\frac{X}{p^N}$, $y=\frac{Y}{p^N}$ where $X,Y$ are $\leq N$-digit numbers in base $p$ such that 
there are no carries in computing $X+Y$ (i.e., for all $i$, $X_{i}+Y_{i}\leq p-1$). 
This shows that $\Delta_{N-1}(p)\subset\Delta_{N}(p)$. 

Let $\Delta_{\infty}(p)=\cup_{N\geq 1}\Delta_{n}(p)$ and $\Delta(p)$ be the closure of $\Delta_{\infty}(p)$. Hence, $\Delta(p)$ is the compact set of pairs $(x,y)\in\R^{2}$ such that $x=0.x_{1}x_{2}\dots$, $y=0.y_{1}y_{2}\dots$ in base $p$, and $x_{i}+y_{i}\leq p-1$. This set is the Sierpinski triangle from above.

\begin{Remark}
The set $\Delta_\infty(p)$ also has a direct representation-theoretic interpretation. Namely it corresponds to the set of weights in the simple representations of the (infinite type) affine group scheme $(SL_{2})_{\mathrm{perf}}$, obtained by `perfecting' $SL_{2}$, see \cite[\S 7.1]{CW}. It would be interesting to have a similar direct interpretation for $\mathtt{C}_{\infty}$ from \autoref{SS:Cantor}. This is less obvious since distribution algebras of perfected group schemes are trivial, see \cite[Lemma~3.1.2]{CW}.
\end{Remark}

It is easy to see that $|\Delta_{N}|=\binom{p}{2}^{N}$. From this it is not hard to deduce the 
well-known fact that the Hausdorff dimension of $\Delta(p)$ 
is the number given by
\begin{gather*}
\avalue=
\log_{p}\binom{p}{2}=1+\log_{p}\frac{p+1}{2}=1+\frac{\ln\frac{p+1}{2}}{\ln p},\quad\text{e.g. }\avalue\approx 1.631\text{ for $p=3$}.
\end{gather*}
This number is transcendental, and can also be seen at the level of the generating function given by
\begin{gather*}
f(z):=\sum_{n=0}^{\infty}(\dim_{\mathbf{k}}L_{n})z^{n}.
\end{gather*}
Steinberg's tensor product theorem
yields that $f$ is a 2-Mahler equation of degree $p$: 
\begin{align*}
f(z)&=\prod_{m=0}^{\infty}\frac{1-(p+1)z^{p^{m+1}}+pz^{(p+1)p^{m}}}{(1-z^{p^{m}})^{2}}=(1+2z+\dots+pz^{p-1})\cdot f(z^{p})
\\
&=\frac{1-(p+1)z^{p}+pz^{p+1}}{(1-z)^2}\cdot f(z^{p}).
\end{align*}
We have $(1+2z+\dots+pz^{p-1})|_{z=1}=\binom{p}{2}$,
so we have for some $C>1$
\begin{gather*}
f\asymp(1-z)^{-\avalue} 
\colon\quad
C^{-1}\cdot(1-z)^{-\avalue}\leq f(z)\leq C\cdot(1-z)^{-\avalue}.
\end{gather*}
This implies that for some $\widetilde{C}>1$, we have, for large enough $n$:
\begin{gather*}
\sum_{j=0}^{n}
\dim_{\mathbf{k}}L_{j}\asymp
n^{\avalue}
\colon\quad
\widetilde{C}^{-1}\cdot n^{\avalue}\leq\sum_{j=0}^{n}\dim_{\mathbf{k}}L_{j}\leq\widetilde{C}\cdot n^{\avalue}.
\end{gather*}
This is illustrated in \autoref{Eq:Gasket} below.
A more detailed asymptotics of $f(z)$ can be obtained as follows. 
Note that 
\begin{gather*}
(1-z)^{\avalue}\cdot f(z)=
\prod_{m=0}^{\infty}
\frac{1+2z^{p^{m}}+\dots+pz^{(p-1)p^{m}}}{(1+z^{p^{m}}+\dots+z^{(p-1)p^{m}})^{\avalue}}. 
\end{gather*}
Thus, we see that 
\begin{gather*}
\lim_{m\to\infty}(1-z^{p^{-m}})^{\avalue}\cdot f(z^{p^{-m}})=g_{p}(z), 
\end{gather*}
where (note that the product below is absolutely convergent)
\begin{gather*}
g_{p}(z):=
\prod_{m=-\infty}^{\infty}
\frac{1+2z^{p^{m}}+\dots+pz^{(p-1)p^{m}}}{(1+z^{p^{m}}+\dots+z^{(p-1)p^{m}})^{\avalue}}
\end{gather*}
is a periodic function in the sense that $g_{p}(z)=g_{p}(z^{p})$. 
This implies that the sequence $L(n)=n^{-\avalue}\sum_{j=0}^{n}\dim_{\mathbf{k}}L_{j}$ also exhibits oscillatory behavior:
\begin{gather}\label{Eq:Gasket}
\begin{tikzpicture}[anchorbase]
\node at (0,0) {\includegraphics[height=4cm]{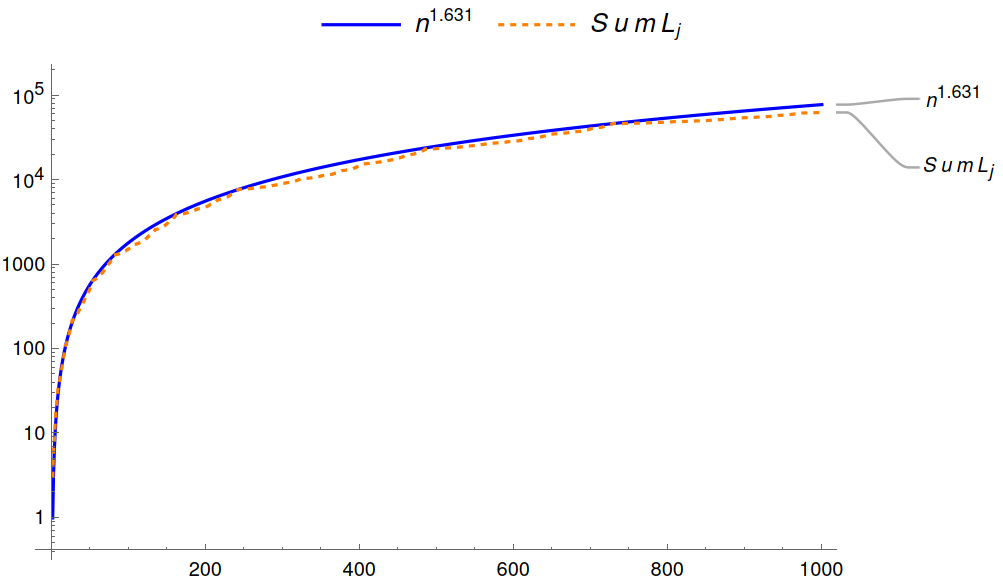}};
\node at (-1.5,1.2) {$p=3$};
\end{tikzpicture}
,\quad
\begin{tikzpicture}[anchorbase]
\node at (0,0) {\includegraphics[height=4cm]{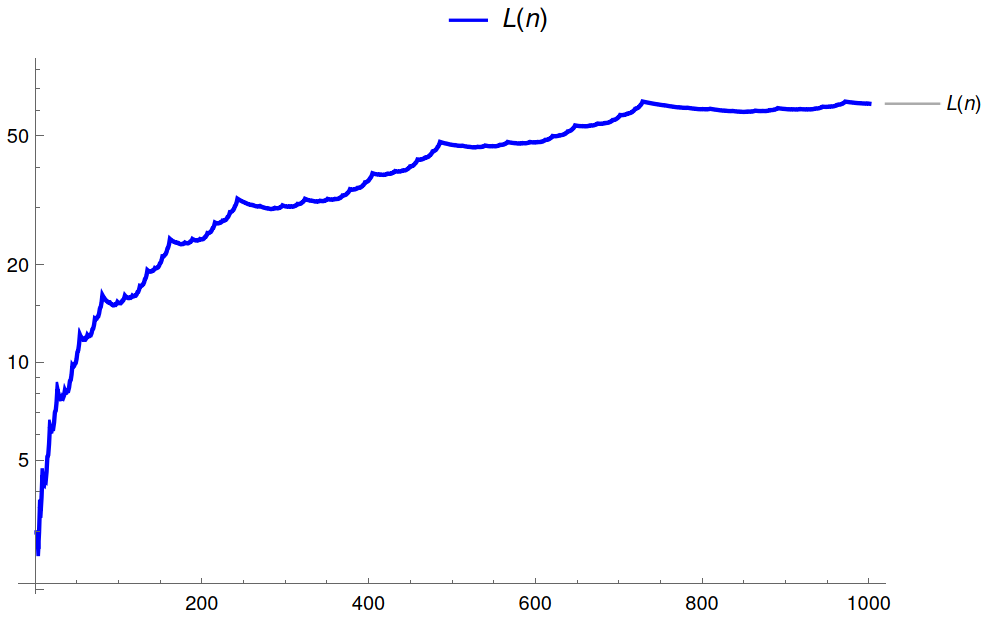}};
\node at (-1.5,1.2) {$p=3$};
\end{tikzpicture}
.
\end{gather}
Let $\mu_{p}$ be the direct image of the (suitably normalized) Hausdorff measure on $\Delta_{\infty}(p)$ under the 
map $(x,y)\mapsto x+y$. Then the analytic function $\mathtt{h}_{p0}(v):=v^{-\avalue}h_{p0}(e^{-v})$ on $(0,\infty)$ is the Laplace transform 
of $\mu_{p}$: 
\begin{gather*}
\mathtt{h}_{p0}(v)=\int_{0}^{\infty}e^{-tv}d\mu_{p}(t),
\\
h_{p0}(x)=\ln(x^{-1})^{\avalue}\int_{0}^{\infty}x^{t}d\mu_{p}(t).
\end{gather*}
Similarly, the sequence $C_{n}(p):=n^{-\avalue}\sum_{j=0}^{n}\dim_{\mathbf{k}}L_{j}$ for large $n$ behaves as 
follows: 
\begin{gather*}
\lim_{m\to\infty}C_{[p^{m}w]}=w^{-\avalue}\mu_{p}([0,w]),
\end{gather*}
Thus, this sequence approaches the periodic function $\log_{p}(n)^{-\avalue}\mu_{p}([0,\log_{p}(n)])$ at infinity. As before, this function is continuous but not differentiable, nor
absolutely continuous.  

\subsubsection{Lengths of tensor powers}\label{SS:Length}

Recall that $\ell(\placeholder)$ denotes the length of a representation.
Finally let us discuss the fractal behavior of the integer sequence $l_{n}=\ell(V^{\otimes n})$ 
where $V\cong\mathbf{k}^{2}$ is the vector representation of $\Gamma$. We will not explicitly consider any fractals in this example, however the origin of what we, based on the previous examples, call `fractal behavior' is still found in the same principles, such as Steinberg's tensor product theorem.

\textit{$p=2$.} For simplicity consider the case $p=2$ first. In this case one can show that the generating function 
\begin{gather*}
L(z):=\sum_{n\geq 0}l_{n}z^{-n-1}
\end{gather*}
is holomorphic for $|z|>2$ and satisfies the functional equation 
\begin{gather*}
L(z)=(1+z)\cdot L(z^{2}-2).
\end{gather*}
In particular, this shows that $l_{2k-1}=l_{2k}$ for all $k\in\Z_{>0}$. 
This is not a Mahler equation, but it turns into one under a simple change of variable. 
Namely, setting $z=w+w^{-1}$ and $wh(w)=L(z)$, we have 
\begin{gather*}
h(w)=(1+w+w^{2})\cdot h(w^{2}),
\end{gather*}
i.e., $h$ is a 2-Mahler function of degree $2$. 

One gets
\begin{gather*}
h(w)=\prod_{n\geq 0}(1+w^{2^{n}}+w^{2\cdot 2^{n}})\text{ for }|w|<1. 
\end{gather*}
Using this, we can compute as previously the asymptotics of $h$ as 
$w\uparrow 1$. Set
\begin{gather*}
\avalue=\log_{2}3\approx 1.585,
\end{gather*}
which is transcendental.
We then have, again by the theory of Mahler functions,
\begin{gather*}
C_{1}\cdot (1-w)^{-\avalue}\big(1+o(1)\big)\leq 
h(w)\leq C_{2}\cdot (1-w)^{-\avalue}\big(1+o(1)\big)\text{ for }w\uparrow 1, 
\end{gather*}
for some $0<C_{1}<C_{2}$ with $C_{1},C_{2}\in\R$. 
Hence, we get
\begin{gather*}
C_{1}\cdot (z-2)^{-\avalue/2}\big(1+o(1)\big)\leq L(z)\leq C_{2}\cdot (z-2)^{-\avalue/2}\big(1+o(1)\big)\text{ for }z\downarrow 2.
\end{gather*}
So from the Tauberian theory, cf. \autoref{SS:Tauber}, it follows that 
\begin{gather*}
\frac{2^{1-\frac{\avalue}{2}}C_{1}}{\Gamma(\frac{\avalue}{2}+1)}
\cdot n^{\frac{\avalue}{2}}
\big(1+o(1)\big)
\leq\sum_{j=0}^{n}\frac{l_{j}}{2^{j}}\leq
\frac{2^{1-\frac{\avalue}{2}}C_{2}}{\Gamma(\frac{\avalue}{2}+1)}
\cdot n^{\frac{\avalue}{2}}
\big(1+o(1)\big)\text{ for }n\to\infty. 
\end{gather*}
Or, equivalently, 
\begin{gather*}
\frac{2C_{1}}{3\Gamma(\frac{\avalue}{2}+1)}\cdot n^{\frac{\avalue}{2}}\big(1+o(1)\big)\leq\sum_{k=1}^{n}\frac{l_{2k}}{2^{2k}}\leq
\frac{2C_{2}}{3\Gamma(\frac{\avalue}{2}+1)}
\cdot n^{\frac{\avalue}{2}}\big(1+o(1)\big)\text{ for }n\to\infty.
\end{gather*}
As previously, the behavior of $h(w)$ as $w\uparrow 1$ can be analyzed as follows. 
We have 
\begin{gather*}
\lim_{k\to\infty}\ln(w^{-2^{-k}})^{\avalue}h(w^{2^{-k}})=\ln(w^{-1})^{\avalue}\lim_{k\to\infty}3^{-k}\prod_{m=-k}^{\infty}
(1+w^{2^{m}}+w^{2\cdot 2^{m}})=\mathtt{h}_{0}(w),
\end{gather*}
where (for $\theta(m)$ as in \autoref{SS:Cantor} above)
\begin{gather*}
\mathtt{h}_{0}(w):=\ln(w^{-1})^{\avalue}\prod_{m=-\infty}^{\infty} 
\frac{1+w^{2^{m}}+w^{2\cdot 2^{m}}}{3^{\theta(m)}}.
\end{gather*}
This is a periodic function in the sense that $\mathtt{h}_{0}(w)=\mathtt{h}_{0}(w^{2})$, and the function $\ln(w^{-1})^\avalue h(w)$ asymptotically approaches the periodic function $\mathtt{h}_{0}(w)$ 
as $w\uparrow 1$. Writing $w=e^{-v}$, we obtain that the function $v^{\avalue}h(e^{-v})$ 
approaches $\mathtt{h}_{0}(e^{-v})$ as $v\downarrow 0$, i.e., 
\begin{gather*}
\lim_{k\to\infty}
(2^{-k}v)^{\avalue}
h(e^{-2^{-k}v})=\mathtt{h}_{0}(e^{-v}).
\end{gather*}
Or, equivalently, 
\begin{gather*}
\lim_{k\to\infty}(2^{-k}v)^{\avalue}L(2\cosh 2^{-k}v)=\mathtt{h}_{0}(e^{-v}).
\end{gather*}
But 
\begin{gather*}
2\cosh 2^{-k}v=2e^{2^{-2k-1}v^{2}}+O(2^{-4k})\text{ for }k\to\infty,
\end{gather*}
so we get 
\begin{gather*}
\lim_{k\to\infty}(2^{-k}v)^{\avalue} 
L(2e^{2^{-2k-1}v^{2}})=\mathtt{h}_{0}(e^{-v}).
\end{gather*}
Thus, setting $x:=e^{v^{2}/2}$, we obtain 
\begin{gather*}
\lim_{k\to\infty}(\ln x^{4^{-k}})^{\avalue/2} 
L(2x^{4^{-k}})=\widetilde{h_{0}}(x):=2^{-\avalue}\mathtt{h}_{0}(e^{-\sqrt{2\ln x}}).
\end{gather*}
Hence, the function 
\begin{gather*}
\widetilde{h}_{0}(e^{\frac{1}{2}4^{u}})=2^{-\avalue}\mathtt{h}_{0}(e^{-2^{u}})
\end{gather*}
is periodic with period $1$ and analytic for $|\mathrm{Im}\,u|\leq\frac{\pi}{\ln 4}$, i.e., the strip of holomorphy is twice as wide as in previous examples. This happens because the role of the prime $p$ in this example is played by the number $4$ (rather than $2$), i.e., as $z\downarrow 2$, the function $L(z)$ behaves as a Mahler function of degree four, rather than a Mahler function of degree two. As we will see, this will lead to much greater regularity of the coefficient sequence $l_{2n}$. 

We would now like to understand the asymptotics of the sequence $l_{2n}$ in more detail. 
If $l_{2k}$ was known to behave sufficiently regularly, we could read off its asymptotics from the asymptotics of the Ces{\'a}ro sums $\sum_{k=1}^{2n}\frac{l_{2k}}{2^{2k}}$ by Abel resummation: 
\begin{gather*}
\frac{2C_{1}}{3\Gamma(\frac{\avalue}{2})}\cdot n^{-1+\frac{\avalue}{2}}\cdot 2^{2n}\big(1+o(1)\big)\leq l_{2n}\leq\frac{2C_{2}}{3\Gamma(\frac{\avalue}{2})}\cdot n^{-1+\frac{\avalue}{2}}\cdot 2^{2n}\big(1+o(1)\big)\text{ for }n\to\infty. 
\end{gather*}
The required regularity is guaranteed by the following lemma, showing that the sequence 
$\frac{l_{2n}}{2^{2n}}$ is decreasing. 

\begin{Lemma}\label{L:Length}
We have $4l_{n}\geq l_{n+2}$. 
\end{Lemma} 

\begin{proof} 
Let 
$Q(z):=-(1-4z^{-2})L(z)=\sum_{n\geq 0}(4l_{n-2}-l_{n})z^{-n-1}$, where we agree that 
$l_{-1}=l_{-2}:=0$. Then the functional equation for $L$ implies that
\begin{gather*}
Q(z)=z^{-4}(1+z)(z^{2}-2)^{2}Q(z^{2}-2).
\end{gather*}
So the coefficients $c_{n}=4l_{n}-l_{n+2}$ of $Q(z)$ satisfy the recursion 
\begin{gather*}
\sum_{n\geq 0}c_{n}z^{-n-3}=z^{-3}+z^{-4}+z^{-5}(1+z^{-1})\sum_{n\geq 0}\frac{c_{n}z^{-2n}}{(1-2z^{-2})^{n+1}}.
\end{gather*}
This shows that $4l_{n}-l_{n+2}\geq 0$ for all $n\geq 0$, which implies the statement.  
\end{proof} 

From this it follows (with some work) that 
\begin{gather*}
\lim_{k\to\infty}2^{-2\cdot[4^{k}y]}[4^{k}y]^{1-\avalue/2}l_{2\cdot[4^{k}y]}=\phi(\log_{4}y),
\end{gather*}
where $\phi(y)$ is a periodic function in the sense that $\phi(y)=\phi(4y)$ (period doubling with respect to the previous examples). Thus, $l_{2n}$ 
behaves roughly like $(2n)^{0.21}\cdot 2^{2n}$ with $1-\frac{1}{2}\log_{2}3\approx0.21$
(as expected, growing faster than in characteristic zero, where 
it is $(2n)^{-1/2}\cdot 2^{2n}$, see \autoref{E:IntroCharZero}). Here is the plot 
(for $l_{n}$; note that $l_{2k+1}=l_{2k}$):
\begin{gather*}
\begin{tikzpicture}[anchorbase]
\node at (0,0) {\includegraphics[height=4.5cm]{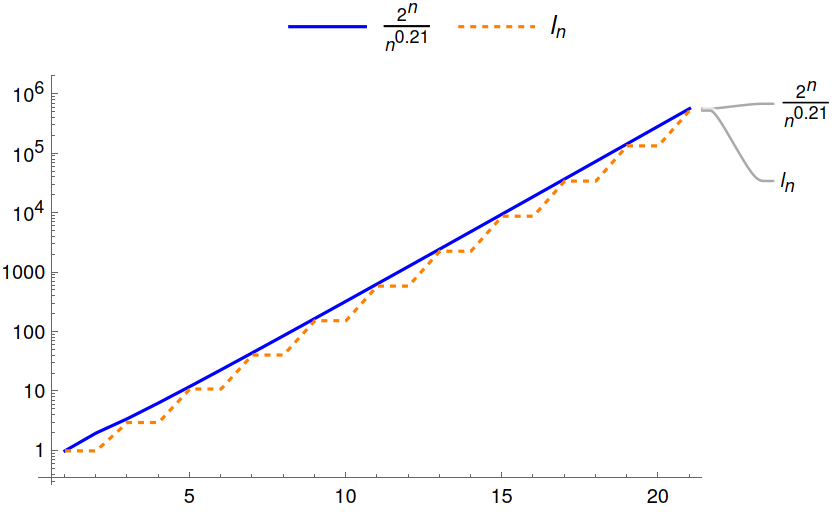}};
\node at (-1.5,1.2) {$p=2$};
\end{tikzpicture}
,\quad
\begin{tikzpicture}[anchorbase]
\node at (0,0) {\includegraphics[height=4.5cm]{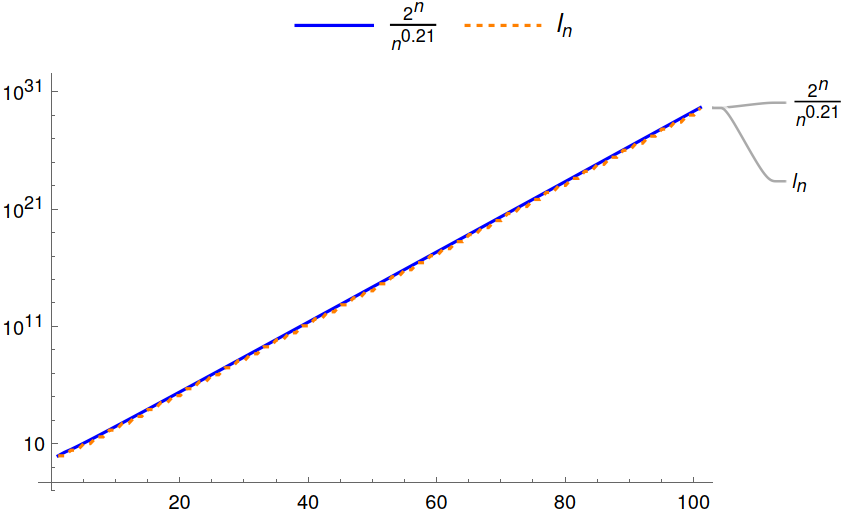}};
\node at (-1.5,1.2) {$p=2$};
\end{tikzpicture}
.
\end{gather*}
To see the oscillation, we zoom into 
$l_{2n}/((2n)^{0.21}\cdot 2^{2n})$:
\begin{gather*}
\begin{tikzpicture}[anchorbase]
\node at (0,0) {\includegraphics[height=4.5cm]{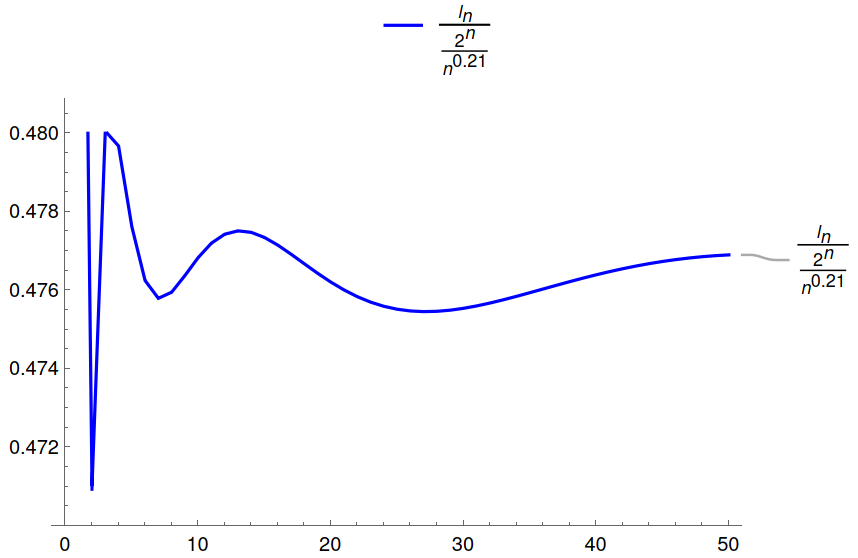}};
\node at (0.5,1.2) {$p=2$, even $n$ zoom};
\end{tikzpicture}
.
\end{gather*}
As before, the Fourier coefficients $a_{n}$ of the 1-periodic function $\phi(4^{u})$ are related to the Fourier coefficients $A_{n}$ of the function 
$2^{-\avalue}\mathtt{h}_{0}(e^{-2^{u}})$ by the formula 
\begin{gather*}
A_{n}=\Gamma(\tfrac{\avalue}{2}+1+\tfrac{2\pi in}{\ln 4})a_{n}.
\end{gather*}
So since $\mathtt{h}_{0}(e^{-2^{u}})$ is analytic in the strip  $|\mathrm{Im}\,u|\leq
\frac{\pi}{\ln 4}$,
we have $|A_{n}|=O(e^{-(\frac{2\pi^2}{\ln 4}-\varepsilon)|n|})$ 
for all $\varepsilon>0$, hence 
\begin{gather*}
|a_{n}|=O(e^{-(\frac{\pi^{2}}{\ln 4}-\varepsilon)|n|})\text{ for }|n|\to\infty.  
\end{gather*}
Thus, the function $\phi(4^{u})$ is analytic in the strip $|\mathrm{Im}\,u|\leq
\frac{\pi}{2\ln 4}$. 

This is a significant difference from the previous examples, where the analogous function was not even absolutely continuous (nor differentiable). This happens due to presence of the change of variable 
$z=w+w^{-1}$ which has a quadratic branch point at $w=1$ and thus causes doubling of periods and widths of strips of holomorphy. 

\textit{General $p$.} The analysis is essentially the same as for $p=2$, 
and we will omit a discussion. Let us simply point out that 
the generating function satisfies
\begin{gather*}
L(w+w^{-1})=w^{1-p}(1+w+w^{2}+\dots+w^{2p-2})\cdot L(w^{p}+w^{-p})
\end{gather*}
which gives
\begin{gather*}
h(w)=(1+w+w^{2}+\dots+w^{2p-2})\cdot h(w^{p}).
\end{gather*}
Hence, we have again a 2-Mahler function of degree $p$. From this we get
\begin{gather*}
\avalue=\log_{p}(2p-1)
\end{gather*}
as the exponent of the subexponential factor.

\subsubsection{Conclusion and goal}

The primary objective of this paper is to study a similar but more complicated problem 
of precisely estimating the \emph{number of indecomposable summands} of $V^{\otimes n}$.
Unlike the study of the length, which grows faster than the number of indecomposable summands due to the non semisimple nature of the category of representations of $SL_{2}(\mathbf{k})$, this problem introduces subtleties. For example, since the 2-Mahler equation of degree $p$ it leads to is inhomogeneous, the generating function $h(w)$ is not a single product but rather a sum of products. Still, the asymptotic behavior of the sequences and functions we consider ends up being very similar to the above examples. There will be an especial similarity with the example in this subsection (e.g. length of $V^{\otimes n}$). Namely, we also observe doubling of periods and strip widths, resulting in a very high degree of regularity of the sequence of interest. 

We will however abstain from exploring specific fractals within this context.

\section{Generating function}\label{SecGen}

We fix a prime $p$. Our first goal is to 
explicitly describe the generating function for 
the sequence $b_{n}$ as in \autoref{Eq:MainSequence}. All objects will depend on $p$, but to keep notation light we usually omit $p$ from notation.

\subsection{The function \texorpdfstring{$F$}{F}}\label{SubSecFunF}

Let $w$ and $z$ be formal variables. By 
expansion, we have a ring inclusion
\begin{gather*}
\Q[[z^{-1}]]\hookrightarrow\Q[[w]],
\text{ corresponding to $z\mapsto w+w^{-1}$}.
\end{gather*}
Moreover, the above restricts to inclusions
\begin{gather*}
\Z[[z^{-1}]]\hookrightarrow\Z[[w]]\quad\text{and}\quad 
\Q(z)=\Q(z^{-1})\hookrightarrow\Q(w).
\end{gather*}
The image of the latter inclusion consists precisely of those rational functions in $w$ invariant under $w\leftrightarrow w^{-1}$. We will use this several times 
tacitly below.

We are now ready to study the generating function for the sequence $b_{n}$.
It is a bit more convenient to shift the generating function for 
$b_{n}$ as follows. Let
\begin{gather*}
H(z):=H_{p}(z):=\sum_{n\geq 0}b_{n}z^{-n-1}\in\Z[[z^{-1}]].
\end{gather*}
We will also regard $H(z)$ as a holomorphic function with domain $|z|>2$ which vanishes at infinity.

It will be convenient for formal manipulations to focus on $H(w+w^{-1})$.
We therefore set
\begin{gather*}
F(w):=F_{p}(w):=\sum_{n\geq 0}
b_{n}(w+w^{-1})^{-n-1}\in\Z[[w]].
\end{gather*}
Of course, by construction, we can also interpret $F(w)$ as a holomorphic function (the singularity of this function that will be important later is now at $w=1$) on
\begin{gather}\label{Eq:Discs}
\C\supset\Omega:=\{w\in\C\mid |w+w^{-1}|>2\}.
\quad\text{Plot of $\Omega$}\colon
\begin{tikzpicture}[anchorbase]
\node at (0,0) {\includegraphics[height=4cm]{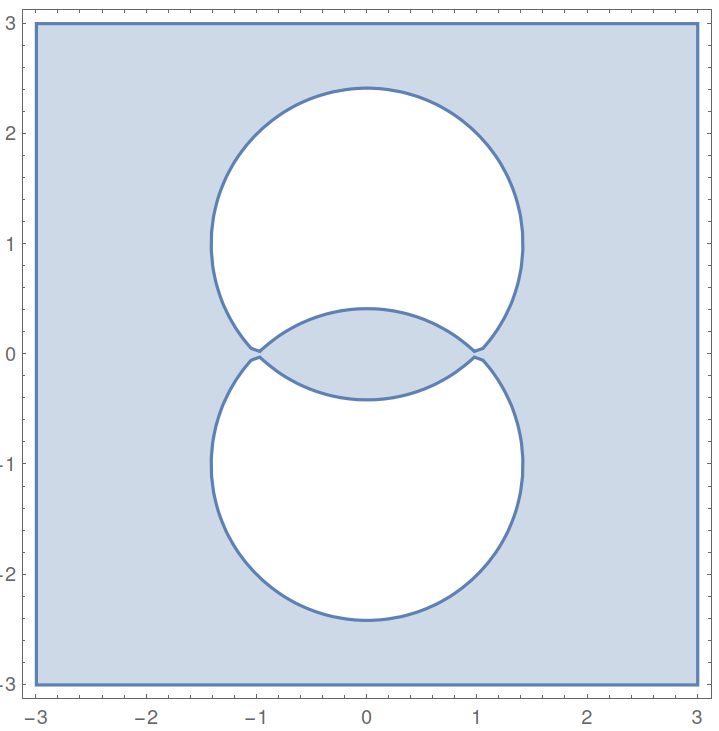}};
\end{tikzpicture}
.
\end{gather}
Indeed, direct manipulation shows that $\Omega$ is the 
set of complex numbers $x+iy$ with
\begin{gather*}
\big(x^{2}+(y-1)^{2}-2\big)
\big(x^{2}+(y+1)^{2}-2\big)>0,
\end{gather*}
yielding the intersection and joint exterior of two disks, as displayed in \autoref{Eq:Discs}.

Our main result of this section is:

\begin{Theorem}\label{th2.4}
We have the following explicit formula for $F(w)$:
\begin{gather*}
F(w)=\sum_{k=0}^{\infty}\frac{w^{p^{k}}(1-w^{p^{k}(p-1)})}{(1-w^{p^{k}})(1+w^{p^{k+1}})}\prod_{j=0}^{k-1} \frac{1-w^{(p+1)p^{j}}}{(1-w^{p^{j}})(1+w^{p^{j+1}})}.
\end{gather*}
\end{Theorem} 

\begin{Remark}
As pointed out to us by Henning Haahr Andersen, a finer version (with a quite different formulation) of \autoref{th2.4} has appeared in \cite{Erd}.
\end{Remark}

The remainder of this section is devoted to the proof of \autoref{th2.4}.

\begin{Example}\label{E:GrowthCompare}
If $p=2$, then the expression in \autoref{th2.4} simplifies to
\begin{gather*}
F(w)=\sum_{k=0}^{\infty}\frac{1}{w^{2^{k}}+w^{-2^{k}}}\prod_{j=0}^{k-1}\left(1+\frac{1}{w^{2^{j}}+w^{-2^{j}}}\right),
\end{gather*}
as one can directly verify.
\end{Example}

\begin{Example}
We can also consider the generating function for a field of characteristic zero for the numbers $b_{n}$, see 
\autoref{E:IntroCharZero}. Then we find
\begin{align*}
H_{\infty}(z)=&\sum_{n\geq 0}\binom{n}{\lfloor n/2\rfloor}z^{-n-1}
= z^{-1}+z^{-2}+3z^{-3}+6z^{-5}+\dots
\\
=&\frac{1}{2}\left(\sqrt{\frac{z+2}{z-2}}-1\right).
\end{align*}
This is well-known and can be derived from \cite[A001405]{oeis}. 
\end{Example}

\subsection{Recursion relations}\label{SS:Recursion}

Recall the category of tilting representations for $SL_{2}(\mathbf{k})$, see e.g. \cite{Do-tilting-alg-groups} and \cite{Ri-good-filtrations} (additional details can be found in \cite{AnStTu-cellular-tilting}), or, using the identification with the 
Temperley--Lieb calculus, \cite{An-simple-tl} or \cite{TuWe-quiver-tilting}.
All the terminology and facts about tilting representations that we use below can be found or derived from these references. We denote the category by  $\Tilt\big(SL_2(\bk)\big)$.

We identify the weight lattice of $SL_{2}(\mathbf{k})$ with $\Z$ and 
the dominant weights with $\N$.
Let $K=K_{p}$ be the split Grothendieck ring of $\Tilt\big(SL_2(\bk)\big)$. Then $K$ has a $\Z$-basis $\{[\mathtt{T}_{i}]\mid i\in\N\}$, where $\mathtt{T}_{i}$ is the indecomposable tilting representation of highest weight $i$.
We then have ring isomorphisms
\begin{gather*}
\Z[u,u^{-1}]^{s}\xleftarrow{\,\mathrm{ch}\,} 
K\xrightarrow{[V]\mapsto x}\Z[x],
\end{gather*}
where $\Z[u,u^{-1}]^{s}$ is the subring of Laurent polynomials symmetric under $u\leftrightarrow u^{-1}$, and $ch$ associates the formal character to a representation. In particular, the composed isomorphism identifies $x$ and $u+u^{-1}$. We will freely use both isomorphisms and switch notation accordingly.

Let $\nu=\nu_{p}\colon\Z[u,u^{-1}]^{s}\simeq\Z[x]\simeq K\to\Z$ be the group homomorphism
\begin{gather*}
\nu\colon K\to\Z,\;[\mathtt{T}_{i}]\mapsto 1,\text{ for all $i\in\N$}.
\end{gather*}
Hence, we have $\nu(x^{n})=b_{n}$.

By Donkin's tensor product formula, for $0\leq k< p$
the operation $\mathtt{T}\mapsto\mathtt{T}^{(1)}\otimes \mathtt{T}_{p-1+k}$ maps indecomposable tilting representations to indecomposable tilting
representations; concretely:
\begin{gather*}
\mathtt{T}_{i}^{(1)}\otimes\mathtt{T}_{p-1+k}\;\simeq\; \mathtt{T}_{pi+p-1+k}.
\end{gather*}
Here $\mathtt{T}^{(1)}$ is the Frobenius twist of $\mathtt{T}$. In particular, we have
\begin{gather*}
ch(\mathtt{T}^{(1)})=\mP\big(ch(\mathtt{T})\big),
\end{gather*}
where $\mP=\mP_{p}$ is the endomorphism of $\Z[u,u^{-1}]$ determined by $u\mapsto u^{p}$.
Consequently, we have
\begin{gather}\label{eqRec}
\nu=\nu\big(\mP(\placeholder)\,ch(\mathtt{T}_{p-1+k})\big),
\end{gather}
for all $0\leq k<p$.

We will interpret \eqref{eqRec} as recursion relations 
to compute $b_{n}=\nu(x^{n})$. To this end, we need to know 
the following few characters: it
is well-known and easy to compute that
\begin{gather*}
ch(\mathtt{T}_{p-1})=\frac{u^{p}-u^{-p}}{u-u^{-1}}
\end{gather*}
and, for $0<k<p$,
\begin{gather}\label{Eqch}
ch(\mathtt{T}_{p-1+k})=(u^{k}+u^{-k})\frac{u^{p}-u^{-p}}{u-u^{-1}}.
\end{gather}

We can reformulate this in terms of the Chebyshev polynomials of the first kind, $T_{n}$ (not to be confused with the notation for tilting modules) and of the second kind, $U_{n}$. Indeed, they can be defined via
\begin{gather}\label{Eq:TPoly}
T_{n}\Big(\frac{w+w^{-1}}{2}\Big)=\frac{w^n+w^{-n}}{2},\quad
\text{and}\quad U_{n}\Big(\frac{w+w^{-1}}{2}\Big)=\frac{w^{n+1}-w^{-n-1}}{w-w^{-1}}.
\end{gather}
We then find:
\begin{align*}
b_{n}=\nu(x^{n})=&\nu\big(2^{n}T_{p}(x/2)^{n}U_{p-1}(x/2)\big)
\\
=&\nu\big(2^{n+1}T_{p}(x/2)^{n}T_{k}(x/2)U_{p-1}(x/2)\big),\quad\text{for $0<k<p$.}
\end{align*}

\begin{Example}\label{E:CharTwoThing}
For $p=2$, the two equations become
\begin{gather*}
\nu\big(x(x^{2}-2)^{n}\big)=\nu(x^{n})\quad\text{and}\quad\nu\big(x^{2}(x^{2}-2)^{n}\big)=\nu(x^{n}).
\end{gather*}
The first relation allows one to compute $b_{2n+1}$ in terms of lower cases and the second relation allows one to compute $b_{2n+2}$. Concretely, one obtains
\begin{gather*}
b_{2n-1}=b_{2n}=\sum_{k=0}^{n-1}\binom{n-1}{k}2^{n-1-k}b_{k}.
\end{gather*}
This gives a very efficient way of computing the numbers $b_{2n-1}=b_{2n}$.
\end{Example}

Let $\mA=\mA_{p}$ be the group endomorphism of $\Z[[w]]$ given by
\begin{gather}\label{DefA}
\mA(r)(w)=\frac{1}{p}\sum_{j=0}^{p-1}r(\zeta^{j}w),\quad\text{for }\; r\in\Z[[w]],
\end{gather}
where $\zeta=e^{2\pi i/p}$ is a primitive complex $p$th root of unity.
In particular, we have
\begin{gather*}
\mA(w^{m})=
\begin{cases} 
w^{m} &\text{if $p$ divides $m$},
\\
0 &\text{otherwise}.
\end{cases}
\end{gather*}

\begin{Lemma}\label{LempEq}
The function $F(w)$ satisfies
\begin{align*}
F(w^{p})&=\mA\big(F(w)\big),\quad\text{and}
\\
F(w^{p})+1&=\mA\left((w^{k}+w^{-k})F(w)\right),\text{ for }1\leq k\leq p-1.
\end{align*}
\end{Lemma}

\begin{proof}
For convenience, we extend the definition of $\nu\colon\Z[u,u^{-1}]^{s}=K\to\Z$ to
\begin{gather*}
\nu\colon K[[w]]\to \Z[[w]],
\end{gather*}
$\Z[[w]]$-linearly. Then we can write
\begin{gather*}
F(w)=\nu\left(\sum_{n\geq 0}\frac{(u+u^{-1})^{n}}{(w+w^{-1})^{n+1}}\right)=
\nu\left(\frac{1}{w+w^{-1}-u-u^{-1}}\right).
\end{gather*}
Using \eqref{eqRec} and \eqref{Eqch}, we can then conclude that
\begin{gather}\label{eqFwp}
F(w^{p})=\nu\left(\frac{u^{p}-u^{-p}}{u-u^{-1}}\frac{R_{k}}{w^{p}-w^{-p}-u^{p}-u^{-p}}\right),\text{ for }0\leq k<p,
\end{gather}
where $R_{0}=1$ and $R_{k}=u^{k}+u^{-k}$ for $i>0$.
We now complete the proof of the first displayed equation, corresponding to $k=0$, the other cases being similar.

Writing
\begin{gather*}
\frac{w}{w-u}=-\sum_{i>0}\frac{w^{i}}{u^{i}}
\end{gather*}
we find that
\begin{gather*}
\mA\left(\frac{w}{w-u}\right)=\frac{w^{p}}{w^{p}-u^{p}}.
\end{gather*}
Consequently, we find
\begin{align*}
\mA\left(\frac{1}{w+w^{-1}-u-u^{-1}}\right)=&\mA\left(\frac{w}{u-u^{-1}}\left(\frac{1}{w-u}-\frac{1}{w-u^{-1}}\right)\right)
\\
=&\frac{w^{p}}{u-u^{-1}}\left(\frac{1}{w^{p}-u^{p}}-\frac{1}{w^{p}-u^{-p}}\right)
\\
=&\frac{u^{p}-u^{-p}}{u-u^{-1}}\frac{1}{w-w^{-p}-u^{p}-u^{-p}}.
\end{align*}
Comparing this expression with~\eqref{eqFwp} indeed proves the first equation of the lemma.
\end{proof}

\subsection{Functional equation}

The following proposition shows that $F(w)$ is a 2-Mahler function of degree $p$.

\begin{Proposition}\label{ThmFunEq}
We have
\begin{gather*}
F(w)=\rp[p](w+w^{-1})+\left(1+\rp[p](w+w^{-1})\right)\cdot F(w^{p}),
\end{gather*}
where $\rp[p]$ is explicitly given by
\begin{gather*}
\rp[p](w+w^{-1})=\frac{w(1-w^{p-1})}{(1-w)(1+w^{p})}.
\end{gather*}
For $p>2$, $\rp[p](z)$ is also determined as the rational function in $z^{-1}$ such that $z\rp[p](z)$ is the Pad{\'e} approximant of order $[\frac{p-1}{2}/\frac{p-1}{2}]$ of 
$\frac{H_{\infty}(z)z^{p+1}}{H_{\infty}(z)z^{2}+z^{p}}$.
\end{Proposition}



\begin{Example}\label{ExThmFunEq}
Consider the case $p=2$.
The two equations in \autoref{LempEq} are
\begin{gather*}
F(w^{2})=\frac{1}{2}\big(F(w)+F(-w)\big)\quad\text{and}\quad 
F(w^{2})+1=\frac{w+w^{-1}}{2}\big(F(w)-F(-w)\big).
\end{gather*}
Taking a linear combination of the two equations neutralising the terms in $F(-w)$ yields
\begin{gather*}
\underbrace{(w+w^{-1})}_{=\rp[2](w+w^{-1})^{-1}}F(w)=
1+\underbrace{(1+w+w^{-1})}_{=\rp[2](w+w^{-1})^{-1}+1}F(w^{2}).
\end{gather*}
This is equivalent to the functional equation in \autoref{ThmFunEq} specialized at $p=2$. The proof of the case $p>2$ given below is a refinement of this argument.
\end{Example}

\begin{proof}[Proof of \autoref{ThmFunEq}]
We can obtain the functional equation by taking an appropriate linear combination of the $p$ equalities in \autoref{LempEq}. Indeed, we use the interpretation of $\mA$ in \eqref{DefA}, and we will sum with appropriate $w$-dependent coefficients to make all terms containing $F(\zeta^{k}w)$ for $k>0$ cancel. We can start by adding up the equalities in \autoref{LempEq} for $1\leq k\leq p-1$ with coefficients $\gamma_{k}(w)\in\Q((w))$ so that the resulting coefficients of $F(\zeta^{k}w)$ yield a value $-C(w)\in\Q((w))$ independent of $1\leq k\leq p-1$.

So we must have
\begin{gather}\label{Eq:C}
\sum_{k=1}^{p-1}(\zeta^{kj}w^{k}+\zeta^{-kj}w^{-k})\gamma_{k}(w)=-C(w)
\end{gather}
for all $j\in [1,p-1]$. Consider
\begin{gather*}
L_{w}(t):=\sum_{k=1}^{p-1}(t^{k}+t^{-k})\gamma_{k}(w)+C(w)\in\Q((w))[t,t^{-1}].
\end{gather*}
By \eqref{Eq:C} and its complex conjugate, $L_{w}(\zeta^{j}w)=L_{w}(\zeta^{j}w^{-1})=0$, for all $j\in[1,p-1]$. It follows that, up to scaling, 
\begin{gather*}
L_{w}(t)=t^{1-p}\prod_{j=1}^{p-1}(t-\zeta^{j}w)(t-\zeta^{-j}w)=t^{1-p}\frac{(t^{p}-w^{p})(t^{p}-w^{-p})}{(t-w)(t-w^{-1})}=
\\
t^{1-p}\left(\sum_{j=0}^{p-1}t^{j}w^{-j}\right)\left(\sum_{i=0}^{p-1}t^{i}w^{i}\right)=
\sum_{i,j=0}^{p-1}t^{i+j-p+1}w^{i-j}.
\end{gather*}
It follows that, setting $[j]_{w}=\frac{w^{j}-w^{-j}}{w-w^{-1}}\in\Z((w))$, we can choose:
\begin{gather*}
\gamma_{k}(w)=\sum_{i=0}^{p-1-k}w^{2i+k+1-p}
=[p-k]_{w},
\quad C(w)=[p]_{w}.
\end{gather*}
Summing over the $p-1$ equations as intended we thus get
\begin{gather*}
p\sum_{j=1}^{p-1}[j]_{w}\big(F(w^{p})+1\big)
=\left(\sum_{j=1}^{p-1}[j]_{w}(w^{p-j}+w^{j-p})\right)F(w)-[p]_{w}\big(F(\zeta w)+\cdots+F(\zeta^{p-1}w)\big).
\end{gather*}
So multiplying the first equation, $F(w^{p})=\mA\big(F(w)\big)$, by $p[p]_{w}$, and adding it to the above equation, we obtain
\begin{gather*}
p[p]_{w}F(w^{p})+p\sum_{j=1}^{p-1}[j]_{w}\big(F(w^{p})+1\big)
=\left([p]_{w}+\sum_{j=1}^{p-1}[j]_{w}(w^{p-j}+w^{j-p})\right)F(w).
\end{gather*}
This can be written as
\begin{gather*}
\sum_{j=1}^{p}(w^{j}-w^{-j}) F(w^{p})+\sum_{j=1}^{p-1}(w^{j}-w^{-j})
=(w^{p}-w^{-p})F(w),
\end{gather*}
subsequently as
\begin{gather*}
\left(\frac{1-w^{p+1}}{1-w}-\frac{1-w^{-p-1}}{1-w^{-1}}\right)F(w^{p})+
\left(\frac{1-w^{p}}{1-w}-\frac{1-w^{-p}}{1-w^{-1}}\right)=(w^{p}-w^{-p})F(w),
\end{gather*}
and finally as the expression in the theorem.

The values $b_{n}$ for $n<p$ are equal to the corresponding numbers in characteristic zero. Consequently, we have
\begin{gather*}
\frac{1}{H_{p}(z)}=\frac{1}{H_{\infty}(z)}+O(z^{1-p}).
\end{gather*}
If $T_{n}$ is the Chebyshev polynomial of the first kind, then we set
\begin{gather*}
P_{p}(x):=2T_{p}(x/2).
\end{gather*}
In other words:
\begin{gather*}
P_{p}(w+w^{-1})=w^{p}+w^{-p}.
\end{gather*}
The functional equation then implies
\begin{gather*}
H_{p}(z)=\rp[p](z)+\big(1+\rp[p](z)\big)\cdot H_{p}\big(P_{p}(z)\big).
\end{gather*}
Since $P_{p}(z)$ is of degree $p$ (as a polynomial in $z$), we have (as functions in $z^{-1}$)
\begin{gather*}
H_{p}\big(P_{p}(z)\big)=z^{-p}+O(z^{-p-2}).
\end{gather*}
Rewriting the functional equation, we thus find
\begin{gather*}
\rp[p](z)^{-1}=\frac{1+H_{p}\big(P_{p}(z)\big)}{H_{p}(z)-H_{p}\big(P_{p}(z)\big)}=\frac{1}{H_{\infty}(z)}+z^{2-p}+O(z^{1-p}).
\end{gather*}

It thus follows that $z^{-1}\rp[p](z)^{-1}$ is the Pad{\'e} approximant of order $[a/b]$ of
\begin{gather*}
\frac{z^{-1}}{H_{\infty}(z)}+z^{1-p}
\end{gather*}
if $a+b<p$ and $z^{-1}\rp[p](z)$ is a rational function in $z^{-1}$ of degree $(a,b)$. 

It remains to verify the final statement.
For $p>2$ consider therefore the expansion
\begin{align*}
\rp[p](w+w^{-1})^{-1}=&\frac{1}{1+w+w^2+\cdots+w^{p-2}}\cdot\frac{1+w^p}{w}
\\
=&\frac{(1+w)w^{\frac{p-3}{2}}}{1+w+w^2+\cdots+w^{p-2}}\cdot\frac{1+w^p}{w^{\frac{p-1}{2}}(1+w)}.
\end{align*}
The left factor on the last line is the inverse of a polynomial in $w+w^{-1}$ of degree $(p-3)/2$, while the right factor is a polynomial in $w+w^{-1}$ of degree $(p-1)/2$. In particular, $z^{-1}\rp[p](z)^{-1}$ is a rational function in $z^{-1}$ of degree $(\frac{p-1}{2},\frac{p-1}{2})$, concluding the proof.
\end{proof}

\subsection{Conclusion}

\begin{proof}[Proof of \autoref{th2.4}]
We can now employ the functional equation from \autoref{ThmFunEq} to derive the closed expression for $F(w)$ in \autoref{th2.4}. Set 
\begin{gather*}
G_{p}(w):=\prod_{j=0}^{\infty}
\frac{1-w^{(p+1)p^{j}}}
{(1-w^{p^{j}})(1+w^{p^{j+1}})}\;\in\Q((w))^{\times}.
\end{gather*}
So 
\begin{gather}\label{eqGG}
G_{p}(w)=\frac{1-w^{p+1}}{(1-w)(1+w^{p})}G_{p}(w^{p})=\left(1+\frac{1}{\rp[p](w+w^{-1})}\right)G_{p}(w^{p}).
\end{gather}
Thus setting 
$F(w)=G_{p}(w)F_{\ast}(w)$, \autoref{ThmFunEq} can be rewritten as
\begin{gather*}
F_{\ast}(w)=F_{\ast}(w^{p})+\frac{w(1-w^{p-1})}{(1-w)(1+w^{p})G_{p}(w)}.
\end{gather*}
So we get
\begin{gather*}
F_{\ast}(w)=\sum_{k=0}^{\infty}
\frac{w^{p^{k}}(1-w^{p^{k}(p-1)})}{(1-w^{p^{k}})(1+w^{p^{k+1}})
G_{p}(w^{p^{k}})},
\end{gather*}
or equivalently
\begin{gather*}
F(w)=\sum_{k=0}^{\infty}
\frac{w^{p^{k}}(1-w^{p^{k}(p-1)})}{(1-w^{p^{k}})(1+w^{p^{k+1}})}
\frac{G_{p}(w)}{G_{p}(w^{p^{k}})}.
\end{gather*}
The conclusion of \autoref{th2.4} then follows from substituting (iterations of) \autoref{eqGG}.
\end{proof}

\section{Asymptotics of the generating function}\label{S:Asym}

Recall that $p$ is a fixed prime.
In order to understand the sequence $b_{n}$ from 
\autoref{Eq:MainSequence} better, we focus on the generating function $H(z)$ as $z$ approaches the radius of convergence $2$, see \autoref{SecGen}, as usual in the theory of asymptotics of generating functions.

\subsection{Asymptotics}\label{SS:AAsymptotics}

Concretely, we will focus on the behavior for $w\downarrow 1$ of $F(w)=H(w+w^{-1})$, viewed as a smooth function
\begin{gather*}
F\colon(0,1)\to\R.
\end{gather*}
This is well-defined, as $(0,1)$ lies in 
$\Omega\subset\mC$ from \autoref{Eq:Discs}. This can be seen directly 
from the formula $(x^{2}+(y-1)^{2}-2)
(x^{2}+(y+1)^{2}-2)>0$, or the plot \autoref{Eq:Discs}.
We will now prove one of our main results:

\begin{Theorem}\label{T:AsymGenFunction}
We have
\begin{gather*}
F(w)=F_{0}(w)\cdot(1-w)^{-\log_{p}(\frac{p+1}{2})}\cdot\big(1+o(1)\big).
\end{gather*}
Here $F_{0}(w)\colon(0,1)\to\R_{>0}$ is real analytic, $F_{0}(w^{p})=F_{0}(w)$ and bounded away from $0$ and $\infty$.
\end{Theorem}


\begin{proof}
Recall from \autoref{ThmFunEq} that
we have
\begin{gather*}
F(w)=\frac{w(1-w^{p-1})}{(1-w)(1+w^{p})}+
\left(1+\frac{w(1-w^{p-1})}{(1-w)(1+w^{p})}\right)\cdot F(w^{p}).
\end{gather*}
A calculation gives
\begin{gather*}
(1-w)(1+w^{p})\cdot F(w)=
w(1-w^{p-1})+(1-w^{p+1})\cdot F(w^{p}).
\end{gather*}
This is a $2$-Mahler function of degree $p$, cf. 
\autoref{Eq:MahlerFirst}, meaning it is of the form
\begin{gather*}
F(w)=r_{1}(w)+r_{2}(w)\cdot F(w^{p}).
\end{gather*}
Let $\lambda$ be a variable.
After clearing denominators so that the 
$r_{i}(w)$ are polynomials, 
the so-called characteristic polynomial 
(as recalled e.g. in \cite{MR3589306}) of a Mahler functional equation as in \autoref{Eq:MahlerFunction} is
\begin{gather*}
r_{0}(w)\cdot\lambda^{s-1}-
r_{1}(w)\cdot\lambda^{s-2}-
r_{2}(w)\cdot\lambda^{s-3}-
\cdots-
r_{s}(w)\cdot\lambda^{0}\in\C[w,\lambda].
\end{gather*}
This in our example is
\begin{gather*}
(1-w)(1+w^{p})\cdot\lambda
-(1-w^{p+1})
\end{gather*}
which has a root at $\lambda(w)=\frac{w^{p+1}-1}{(w-1)(1+w^{p})}$.
At the relevant singularity we get
\begin{gather*}
\lim_{w\uparrow 1}\lambda(w)=\frac{p+1}{2},
\end{gather*}
which is the eigenvalue of the Mahler function $F(w)$.
Now the classical theory of Mahler functions, see e.g.
\cite[Theorem 1]{MR3589306}, implies that the 
log with base the degree of the Mahler equation of this 
eigenvalue is minus the exponent of $(1-w)$ in the asymptotic expansion. The remaining parts of the theorem follow 
also from \cite[Theorem 1]{MR3589306}.
\end{proof}

In the reminder of this section we make \autoref{T:AsymGenFunction} more explicit.

\subsection{The oscillating factor}\label{SS:ALimit}

We define the following function on $(0,1)$, rescaling $F$:
\begin{gather*}
\mathbf{F}(w):=\big(\ln(w^{-1})\big)^{\log_{p}(\frac{p+1}{2})}F(w)\colon(0,1)\to\R_{>0}.
\end{gather*}

\begin{Lemma}
As a function on $(0,1)$ we have:
\begin{gather*}
\mathbf{F}(w^{p^{-r}})=\left(\tfrac{p+1}{2}\right)^{-r}\big(\ln(w^{-1})\big)^{\log_{p}(\frac{p+1}{2})}\sum_{k=-r}^{\infty}\frac{w^{p^{k}}(1-w^{p^{k}(p-1)})}{(1-w^{p^{k}})(1+w^{p^{k+1}})}\prod_{j=-r}^{k-1} \frac{1-w^{(p+1)p^{j}}}{(1-w^{p^{j}})(1+w^{p^{j+1}})}.
\end{gather*}
\end{Lemma}

\begin{proof}
Directly from \autoref{th2.4}.
\end{proof}

We will now show that there exists a function $\mathbf{F}_{0}(w)\colon(0,1)\to\R$ with $\lim_{r\to\infty}\mathbf{F}(w^{p^{-r}})=\mathbf{F}_{0}(w)$, defined pointwise, 
satisfying $\mathbf{F}_{0}(w^{p})=\mathbf{F}_{0}(w)$. 
To this end, consider the product in $\Z[[w]]$ given by
\begin{gather*}
\Pi(w):=\prod_{j=1}^{\infty}\frac{2}{p+1}\frac{1-w^{(p+1)p^{-j}}}{(1-w^{p^{-j}})(1+w^{p^{-j+1}})}=
\prod_{j=1}^{\infty}\frac{2}{p+1}\frac{1+w^{p^{-j}}+w^{2p^{-j}}+\dots+w^{pp^{-j}}}{1+w^{p^{-j+1}}}.
\end{gather*}
Define the power series
\begin{gather*}
\mathbf{F}_{0}(w):=\big(\ln(w^{-1})\big)^{\log_{p}(\frac{p+1}{2})}\sum_{k=-\infty}^{\infty}\frac{w^{p^{k}}(1-w^{p^{k}(p-1)})}{(1-w^{p^{k}})(1+w^{p^{k+1}})}\Pi(w^{p^{k}})\left(\frac{p+1}{2}\right)^{k}.
\end{gather*}

\begin{Proposition}\label{P:FzeroThing}
For $0<w<1$, we have 
\begin{gather*}
\lim_{r\to\infty}\mathbf{F}(w^{p^{-r}})=\mathbf{F}_{0}(w).
\end{gather*}
Moreover, $\mathbf{F}_{0}\colon(0,1)\to\R_{>0}$ is a real analytic and oscillatory term, $\mathbf{F}_{0}(w^{p})=\mathbf{F}_{0}(w)$, and is bounded away from $0$ and $\infty$.
\end{Proposition}

\begin{Remark}
The function $\mathbf{F}_{0}$ in \autoref{P:FzeroThing} is 
a rescaling of $F_{0}$ from \autoref{T:AsymGenFunction}, so describes the oscillation of $F_{0}$.
\end{Remark}

\begin{proof}[Proof of \autoref{P:FzeroThing}]
It is easy to see that $\mathbf{F}_{0}(w)$ is well-defined for $w\in(0,1)$:
Firstly, the factors in the expression $\Pi(w^{p^{k}})$ converge to $1$ 
rapidly, and this implies that $\Pi(w^{p^{k}})$ itself converges to some number 
in $(0,1)$, and we can assume that $\Pi(w^{p^{k}})$ is equal to $1$.
Second, the left term in the sum goes to $(p-1)/2$ and takes values in $\big(1,(p-1)/2\big)$ for negative $k$, so the negative part of the sum converges.
Finally, for positive $k$ the summands go to $0$ rapidly and the sum also converges.

As in the previous paragraph, $\Pi(w^{p^{k}})$ converges to $1$ for $w\in(0,1)$.
We then obtain 
\begin{gather*}
\mathbf{F}(w^{p^{-r}})\Pi(w^{p^{-r}})=\big(\ln(w^{-1})\big)^{\log_{p}(\frac{p+1}{2})}\sum_{k=-r}^{\infty}\frac{w^{p^{k}}(1-w^{p^{k}(p-1)})}{(1-w^{p^{k}})(1+w^{p^{k+1}})}\Pi(w^{p^{k}})\left(\frac{p+1}{2}\right)^{k}.
\end{gather*}
The claim follows.
\end{proof}

\begin{Proposition}\label{P:SingularitiesPGeneral}
The series $\mathbf{F}_{0}(p^{-v})$ 
converges absolutely and uniformly on compact sets in the region $\mathrm{Re}\,v>0$
but has a dense set of singularities on the imaginary axis. 
\end{Proposition}

\begin{proof}
This follows as for $p=2$ proven in \autoref{SS:ConvergencepIsTwo} below.
Details are omitted.
\end{proof}

\subsection{Example for \texorpdfstring{$p=2$}{p=2}}\label{SS:ConvergencepIsTwo}

For $p=2$, the expression for $\mathbf{F}_{0}$ simplifies to
\begin{gather*}
\mathbf{F}_{0}(w)=\big(\ln(w^{-1})\big)^{\log_{2}(\frac{3}{2})}\sum_{k=-\infty}^{\infty}\left(\frac{3}{2}\right)^{k}\frac{1}{w^{2^{k}}+w^{-2^{k}}}\prod_{j=1}^{\infty}\frac{2}{3}
\left(1+\frac{1}{w^{2^{k-j}}+w^{-2^{k-j}}}\right).
\end{gather*}
In particular, we have 
\begin{gather*}
\mathbf{F}_{0}(2^{-v})=\ln 2\sum_{k=-\infty}^{\infty}\frac{(2^{k}v)^{\log_{2}(\frac{3}{2})}}{2^{2^{k}v}+2^{-2^{k}v}}\Phi(2^{k}v),
\text{ where }
\Phi(v):=\prod_{j=1}^{\infty}
\left(1-\frac{(2^{2^{-j-1}v}-2^{-2^{-j-1}v})^{2}}{3(2^{2^{-j}v}+2^{-2^{-j}v})}\right).
\end{gather*}

\begin{Lemma} 
The product $\Phi(v)$ converges absolutely and uniformly on compact sets (not containing zeros and poles) to a meromorphic function of $v\in\C$.
\end{Lemma} 

\begin{proof}
This holds since $2^{2^{-j}v}\to 1$ exponentially fast as $j\to\infty$.
\end{proof}

The poles of factors in $\Phi$ are solutions of the equation $e^{2^{-j}v\ln 2}=\pm i$, 
i.e.,
\begin{gather*}
v=\frac{2^{j-1}(2n+1)\pi i}{\ln 2}.
\end{gather*}
On the other hand, zeros occur when $2^{2^{-j+1}v}$ nontrivial cube roots of $1$,  
so 
\begin{gather*}
v=\frac{2^{j}(3n\pm 1)\pi i}{3\ln 2},
\end{gather*}
and they have multiplicities (the number of factors $2$ in $3n\pm 1$). 

\begin{Proposition}\label{P:Singularities}
The series $\mathbf{F}_{0}(2^{-v})$ 
converges absolutely and uniformly on compact sets in the region $\mathrm{Re}\,v>0$
but has a dense set of singularities on the imaginary axis. 
\end{Proposition}

\begin{proof}
This holds by the above discussion.
\end{proof}

\section{Monotonicity}\label{S:Monoton}

In this section we address the monotonicity of our main sequence.

\subsection{Neighboring values of \texorpdfstring{$b_{n}$}{bn}}

We will now prove:

\begin{Theorem}\label{T:Mono}
We have $b_{n+2}\leq 4b_{n}$, or equivalently 
$b_{n+2}/2^{n+2}\leq b_{n}/2^{n}$,  for all $n\in\N$.
\end{Theorem} 

\begin{proof}
The proof will occupy this section, and is split into a few lemmas.
For the first lemma up next, let $L_{n}(w)\in\Z[w,w^{-1}]$ for $n\in\N$ be a sequence of Laurent polynomials such that, for some $\nu\in\Z$, we have
\begin{gather}\label{eqq1}
L_{n}(w^{-1})=w^{\nu}L_{n}(w)
\quad\text{and}\quad
L_{n+1}(w)+L_{n-1}(w)=(w+w^{-1})L_{n}(w),n\in\Z_{\geq 1}.
\end{gather}
Note that $L_{n}/L_{n+1}\in\Q(z)\subset\Q[[z^{-1}]]$, where we use $\Q(z)\subset \Q(w)$ as in \autoref{SubSecFunF}.
Now we will prove that actually $L_{n}/L_{n+1}\in\N[[z^{-1}]]$.

\begin{Lemma}\label{L:AuxLem}
If $L_{0}/L_{1}$ is a positive power series in $z^{-1}$, 
then the same is true for $L_{n}/L_{n+1}$ for all $n\in\N$. 
\end{Lemma} 

\begin{proof} 
We have 
\begin{gather*}
\frac{L_{n}}{L_{n+1}}=\frac{L_{n}}{zL_{n}-L_{n-1}}=\frac{z^{-1}}{1-z^{-1}\frac{L_{n-1}}{L_{n}}}.
\end{gather*}
So the statement follows by induction on $n$. 
\end{proof} 

For any $n\in\frac{1}{2}\Z_{>0}$ and an integer $1\leq r\leq n$, let us define
\begin{gather*}
K_{r,n}(z):=\frac{w^{n-r}+w^{-n+r}}{w^{n}+w^{-n}}
\in\Z[[z^{-1}]].
\end{gather*}
Similarly, define 
\begin{gather*}
M_{r,n}(z):=\frac{w^{n-r}+w^{n-r-2}+...+w^{-n+r}}{w^{n}+w^{n-2}+...+w^{-n}}=\frac{w^{n-r+1}-w^{-n+r-1}}{w^{n+1}-w^{-n-1}}\in\Z[[z^{-1}]].
\end{gather*}

\begin{Lemma}\label{L:AuxLemTwo}
\leavevmode
\begin{enumerate}

\item The function 
$K_{r,n}(z)$ has positive Taylor coefficients in $z^{-1}$.

\item The function 
$M_{r,n}(z)$ has positive Taylor coefficients in $z^{-1}$.

\end{enumerate}
\end{Lemma}

\begin{proof}
\textit{(a).} If $n$ is an integer, let $L_{n}(w)=w^{n}+w^{-n}$, and 
if $n$ is an honest half integer, let $L_{n-\frac{1}{2}}(w)=w^{\frac{1}{2}}(w^{n}+w^{-n})$. Then $L_{n}$ satisfy \eqref{eqq1} with $\nu=0$ or $\nu=1$. Moreover, $L_{0}/L_{1}=2z^{-1}$ 
in the integer case and $L_{0}/L_{1}=\frac{1}{z-1}
=\frac{z^{-1}}{1-z^{-1}}$ in the non-integer case. 

Thus, the result follows from 
\autoref{L:AuxLem}, as we can write
\begin{gather*}
K_{r,n}(z)=\frac{L_{n-r}(z)}{L_n(z)}
=\prod_{j=1}^{r}\frac{L_{n-j}(z)}{L_{n-j+1}(z)}.
\end{gather*}

\textit{(b).} Let $L_{n}(w)=w^{n}+w^{n-2}+\dots+w^{-n}$. 
Then $L_{n}$ satisfy \eqref{eqq1}, and $L_{0}/L_{1}=z^{-1}$. Thus, as in part (a), the result follows from 
\autoref{L:AuxLem}.
\end{proof} 

We let
\begin{gather*}
f(z)=\sum_{n\geq 0}\frac{b_{n}}{2^{n}}z^{-n-1}=
2H(2z).
\end{gather*}
Moreover, we consider the renormalization
\begin{gather*}
\psi(z):=-(1-z^{-2})f(z)=\sum_{n\geq 0}c_{n}z^{-n-1}.
\end{gather*}
Note that $b_{0}=b_{1}=1$,
so that $c_{0}=-1$ and $c_{1}=-1/2$ are not positive. But we prove:

\begin{Lemma}\label{L:AuxLemFour}
The Taylor coefficients $c_n$ in $z^{-1}$ of the function $\psi(z)$ satisfy $c_n>0$ for $n>1$.
\end{Lemma}

\begin{Example}\label{E:AuxLemFour}
For $p=2$ the function $\psi(z)$ Taylor expands in $z^{-1}$ as
\begin{gather*}
\psi(z)=-z-\frac{1}{2}z^{2}+\frac{3}{4}z^{3}
+\frac{1}{8}z^{4}+\frac{1}{16}z^{5}
+\frac{3}{32}z^{6}+\frac{3}{64}z^{7}
+\frac{7}{128}z^{8}+\frac{7}{256}z^{9}+\dots
\end{gather*}
which can be obtained from the data listed in \autoref{SS:MainData}.
\end{Example}

\begin{proof}[Proof of \autoref{L:AuxLemFour}]
Recall the 2-Mahler coefficients $\rp[p](z)$ from \autoref{ThmFunEq}, and also the functional equation for $H(z)$ given in the proof of \autoref{ThmFunEq}, namely:
\begin{gather*}
H(z)=\frac{1}{\rp[p](z)}+\left(1+\frac{1}{\rp[p](z)}\right)\cdot H\big(2T_{p}(z/2)\big).
\end{gather*}
Here $T_{m}(x)$ is the 
$m$th Chebyshev polynomial of the first kind with $T_{-1}=0$, as recalled in \autoref{Eq:TPoly}. Set also $s_{p}(z)=\rp[p](2z)$.

Generalzing the calculation in 
\autoref{SS:GrowthpisTwo} below, we get the following expression (note that $\psi(z)=-2(1-z^{-2})H(2z)$ which gives a rescaling of the above functional equation):
\begin{gather*}
\psi(z)=
-2\frac{1-z^{-2}}{s_{p}(z)}
+\left(1+\frac{1}{s_{p}(z)}\right)
\frac{1-z^{-2}}{1-T_{p}(z)^{-2}}
\cdot\psi\big(T_{p}(z)\big)
.
\end{gather*}
We let $M_{n}(z)=2\big(T_{n}(z)+T_{n-2}(z)+\dots+T_{0\text{ or 1}}(z)\big)$, with the last terms being either $T_{0}(z)$ 
or $T_{1}(z)$, depending on the parity. A calculation gives
\begin{gather*}
\frac{1-z^{-2}}{1-T_{p}(z)^{-2}}
=
\frac{z^{-2}T_{p}(z)^{2}}{M_{p-1}(z)^{2}}.
\end{gather*}
To see this observe that this is equivalent to
$\big(T_{p}(z)^{2}-1\big)/(z^{2}-1)=M_{p-1}(z)^{2}$. This in turn follows by setting $2z=w+w^{-1}$ and we get
\begin{gather*}
\frac{T_{p}(z)^{2}-1}{z^{2}-1}
=\frac{(w^{p}+w^{-p})^{2}-4}{(w+w^{-1})^{2}-4}
=\left(\frac{w^{p}-w^{-p}}{w-w^{-1}}\right)^{2}
\\
=\left(w^{p-1}+w^{p-3}+\dots+w^{-p+3}+w^{-p+1}\right)^{2}=M_{p-1}(z)^{2}.
\end{gather*}
We thus get
\begin{gather*}
\psi(z)=
-2\frac{1-z^{-2}}{s_{p}(z)}
+\left(1+\frac{1}{s_{p}(z)}\right)
\frac{z^{-2}T_{p}(z)^{2}}{M_{p-1}(z)^{2}}
\cdot\psi\big(T_{p}(z)\big)
.
\end{gather*}
We get:
\begin{align*}
-z^{-1}-\frac{1}{2}z^{-2}+\psi_{\geq 2}(z)
=&-2\frac{1-z^{-2}}{s_{p}(z)}
-\left(1+\frac{1}{s_{p}(z)}\right)
\frac{z^{-2}\big(2T_{p}(z)+1\big)}{4M_{p-1}(z)^{2}}
\\
&+\left(1+\frac{1}{s_{p}(z)}\right)
\frac{z^{-2}T_{p}(z)^{2}}{M_{p-1}(z)^{2}}\cdot
\psi_{\geq 2}\big(T_{p}(z)\big)
.
\end{align*}
We rewrite this as
\begin{gather*}
\psi_{\geq 2}(z)=R(z)
+S(z)\frac{z^{-2}T_{p}(z)^{2}}{M_{p-1}(z)^{2}}\cdot\psi_{\geq 2}\big(T_{p}(z)\big),
\end{gather*}
where
\begin{gather*}
R(z)=z^{-1}+\frac{1}{2}z^{-2}-2\frac{1-z^{-2}}{s_{p}(z)}-\left(1+\frac{1}{s_{p}(z)}\right)\frac{z^{-2}\big(2T_{p}(z)+1\big)}{4M_{p-1}(z)^{2}},
\\
S(z)=1+\frac{1}{s_{p}(z)}.
\end{gather*}
Recall that the polynomial $T_{m}(z)$ has real roots. Furthermore, $T_{m}(z)$ or 
$T_{m}(z)/z$ depends on $z^{2}$. Both observations together imply that we have that $1/T_{m}(z)$ 
is a positive series in $z^{-1}$. The same applies to $M_{m}(z)$. 
So it suffices to check that the series 
$S(z)$ and $R(z)$ are positive. The next two lemmas imply these facts.

\begin{Lemma}\label{L:AuxLemThree} 
The function $\rp[p](z)^{-1}=\frac{w(1-w^{p-1})}{(1-w)(1+w^{p})}$ 
has positive Taylor coefficients in $z^{-1}$. The same holds for 
$s_{p}(z)^{-1}$.
\end{Lemma}  

\begin{proof} 
For $p=2$ we have 
\begin{gather*}
\frac{1}{\rp[p](z)}=\frac{w(1-w)}{(1-w)(1+w^{2})}=\frac{w}{1+w^{2}}
=\frac{1}{z},
\end{gather*}
which is manifestly positive.

Hence, we can assume that $p$ is odd.
We have 
\begin{gather*}
\frac{1}{\rp[p](z)}=\frac{w+w^{2}+\dots+w^{p-1}}{1+w^{p}}=K_{1,\frac{p}{2}}+K_{2,\frac{p}{2}}+\dots+K_{\frac{p-1}{2},\frac{p}{2}}.
\end{gather*}
So the result follows from \autoref{L:AuxLemTwo}.(a).

That $s_{p}(z)^{-1}$ is positive is then immediate from
$s_{p}(z)=\rp[p](2z)$.
\end{proof} 

\begin{Lemma}\label{L:AuxLemFive}
The function $z^{2}R(z)$ has positive Taylor coefficients in $z^{-1}$.
\end{Lemma} 

\begin{proof}
A calculation shows that
\begin{gather*}
z^2R(z)=
\underbrace{2\frac{w(1-w^{2p-2})}{1-w^{2p}}}_{Y_{1}}
+
\underbrace{\frac{w^2(1-w^{2p-4})}{1-w^{2p}}}_{Y_{2}}
+
\underbrace{\frac{w^{p-1}(1-w^2)}{1-w^{2p}}}_{Y_{3}}
\\
+
\underbrace{\frac{w^{p}(1-w^{2p-2})(1-w^2)}{(1-w^{2p})^2}}_{Y_{4}}
+
\underbrace{\frac{w(1-w^{p-1})}{(1-w)(1+w^p)}\frac{w^{2p-2}(1-w^2)^2}{(1-w^{2p})^2}}_{Y_{5}}
.
\end{gather*}
The $Y_{j}$, 
as indicated above, are positive for $j=1,2,3,4,5$. Indeed, $Y_{1}$, $Y_{2}$, $Y_{3}$, $Y_{4}$ 
are positive by \autoref{L:AuxLemTwo} and $Y_{5}$ is positive 
by \autoref{L:AuxLemTwo} and \autoref{L:AuxLemThree}.
\end{proof}

\autoref{L:AuxLemFive} implies that $R(z)$ is positive, while $S(z)$ is positive 
by \autoref{L:AuxLemThree}. This finishes the proof of \autoref{L:AuxLemFour}.
\end{proof}

Taking all together yields \autoref{T:Mono}.
\end{proof}

\subsection{Example for \texorpdfstring{$p=2$}{p=2}}\label{SS:GrowthpisTwo}

For $p=2$ the calculation in \autoref{L:AuxLemFive} is rather straightforward.
In this case \autoref{ThmFunEq} implies the functional equation
\begin{align*}
\psi(z)&=-(1-z^{-2})f(z)
=-2(1-z^{-2})H(2z)
\\
&=-2(1-z^{-2})\big(
\tfrac{1}{2z}+(1+\tfrac{1}{2z})
\cdot H(2T_{2}(z))\big)
\\
&=-2(1-z^{-2})\big(\tfrac{1}{2z}+\tfrac{1}{2}(1+\tfrac{1}{2z})
\tfrac{(2z^{2}-1)^2}{4z^{4}}
\cdot f(2z^{2}-1)\big)
\\
&=
-z^{-1}+z^{-3}+(1+\tfrac{1}{2z})
\tfrac{(2z^{2}-1)^2}{4z^{4}}\cdot\psi(2z^{2}-1)
\\
&=-2\frac{1-z^{-2}}{s_{2}(z)}
+\left(1+\frac{1}{s_{2}(z)}\right)\frac{z^{-2}T_{2}(z)^{2}}{M_{1}(z)^{2}}
\cdot\psi\big(T_{2}(z)\big)
.
\end{align*}


This gives a recursion of the coefficients $c_{n}$, namely:
\begin{gather*}
\sum_{n\geq 2}c_{n}z^{-n-1}
=\frac{3z^{-3}}{4}+\frac{z^{-4}}{8}+\frac{z^{-5}}{16}+
\frac{1}{4}z^{-4}\left(1+\frac{1}{2}z^{-1}\right)
\sum_{m\geq 2}c_{m}\frac{(2z^2)^{-m+1}}{(1-\frac{1}{2z^2})^{m-1}}.
\end{gather*}
This recursion has positive coefficients, hence, \autoref{L:AuxLemFive} follows. 

\section{The main theorem}\label{S:MainTheorem}

We now prove \autoref{T:MainTheorem}
after an auxiliary lemma for which we use \autoref{Eq:Asymp}.

\subsection{A Tauberian lemma}\label{SS:Tauber}

The following type of result, often called \emph{Tauberian theory}, is standard and just reformulated 
to suit our needs, see for instance 
\cite[\S 7.53]{MR3155290} or \cite[Theorem 2.10.2]{bigotebook} for related results.

\begin{Lemma}\label{L:Tauber}
\leavevmode
\begin{enumerate}

\item Consider a sequence $(a_{n})_{n\in\Z_{>0}}$ with $a_{n}\in\R_{\geq 0}$, for which the series $f(z)=\sum_{n=1}^{\infty}a_{n}z^{n}$ has radius of convergence $1$. If, for some $t\in(0,1)$, we have,
\begin{gather*}
f(z)\;\asymp_{[t,1)}\;(1-z)^{-\beta},
\end{gather*}
for some $\beta\in\R_{>0}$, then:
\begin{gather*}
\sum_{k=1}^{n}a_{k}\;\asymp_{\Z_{>0}}\;n^{\beta}.
\end{gather*}

\item Assume that, additionally to (a), for some $r\in\Z_{>0}$ and $B\in\R_{>0}$ we have
\begin{gather*}
a_{n+r}\leq a_{n}\leq B\cdot a_{n+r-1},
\end{gather*}
for all $n\in\Z_{>0}$. Then:
\begin{gather*}
a_{n}\;\asymp_{\Z_{>0}}\;n^{\beta-1}.
\end{gather*}
\end{enumerate}
\end{Lemma}

\begin{proof}
\textit{Claim (a)} We set $s_{n}:=\sum_{k=1}^{n}a_{k}$.

We first prove the upper bound on $s_{n}$. For every $z\in[t,1)$ we have
\begin{gather*}
s_{n}z^{n}=\sum_{k=1}^{n}a_{k}z^{n} 
\leq\sum_{k=1}^{n}a_{k}z^{k}
\leq f(z)\leq\frac{A}{(1-z)^{\beta}},
\end{gather*}
for some $A\in\R_{>0}$. The function $1/\big(z^{n}(1-z)^{\beta}\big)$ attains its minimum at $z=n/(\beta+n)$ (which is larger than $t$ for $n$ sufficiently large), allowing us to reformulate the above inequality as
\begin{gather*}
s_{n}\leq A\frac{(\beta+n)^{\beta+n}}{\beta^{\beta}n^{n}}.
\end{gather*}
The latter then implies that for some $C>0$, we have
$s_{n}\leq C\cdot n^{\beta}$. 

Deriving the inequality in the other direction is more subtle. However, as explained in \cite{353943}, this (in fact, both inequalities) is a special case of 
the de Haan--Stadtm{\"u}ller theorem, see \cite[Theorem~2.10.2]{bigotebook}.

\textit{Claim (b)} By the conclusion of part (a), we know that for some $C>1$
\begin{gather*}
C^{-1}\cdot n^{\beta}\leq\sum_{k=1}^{n}a_{k}
\leq C\cdot n^{\beta},\quad\text{for all }n>0.
\end{gather*}
We will use this freely.

We again start with the upper bound.
For all $n\geq r$ we find from monotonicity that
\begin{gather*}
\tfrac{n}{r}a_{n}\leq\sum_{j=0}^{\lfloor (n-1)/r\rfloor}a_{n-rj}\leq C\cdot n^{\beta}.
\end{gather*}
From this we can derive 
$a_{n}\leq C^{\prime}\cdot n^{\beta-1}$ for 
all $n\in\Z_{>0}$ for $C^{\prime}=rC>1$.

For the lower bound we consider the case $r=1$ first (in which case we can take $B=1$ and the second inequality in (b) is trivial). For any $N\geq n$ we have
\begin{gather*}
(N-n)a_{n}\geq\sum_{i=n+1}^{N}a_{i}\geq C^{-1}\cdot N^{\beta}-C\cdot n^{\beta}.
\end{gather*}
Since $\beta>0$, we can choose $m\in\Z_{>0}$ for which $m^{\beta}>C^{2}$. For $N=mn$, we can rewrite the above inequality as
\begin{gather*}
a_{n}\geq\frac{C^{-1}m^{\beta}-C}{m-1}\cdot n^{\beta-1},\quad\text{for all }n>0.
\end{gather*}
By construction the factor in front of $n^{\beta-1}$ is positive and we are done.

In case $r>1$ we can similarly show that
\begin{gather}\label{sumeq}
a_{n}+a_{n-1}+\cdots+a_{n-r+1}\geq A\cdot n^{\beta-1}
\end{gather}
for some $A>0$. On the other hand, we have
\begin{gather*}
r\cdot a_{n}\geq a_{n}+B^{-1}\cdot a_{n+1-r}+\dots+B^{1-r}\cdot a_{n-(r-1)^2}\geq a_{n}+\sum_{i=1}^{r}B^{-i}\cdot a_{n+i-r}.
\end{gather*}
In particular, we can take $A^{\prime}>0$ so that
\begin{gather*}
a_{n}\geq A^{\prime}\cdot(a_{n}+a_{n-1}+\dots+a_{n+1-r}),
\end{gather*}
which together with \eqref{sumeq} concludes the proof.
\end{proof}

\begin{Remark}
The conclusion in \autoref{L:Tauber}.(b) does not follow without the additional assumption $a_{n}\geq B\cdot a_{n+r-1}$ in case $r>1$. Indeed, it suffices to take the sequence
\begin{gather*}
a_{n}=
\begin{cases}
1&\text{if $n$ is odd},
\\
0&\text{if $n$ is even}.
\end{cases}
\end{gather*}
Then $\sum_{k=0}^{n}a_{k}\asymp_{\N}n$ and $a_{n+2}\leq a_{n}$, but 
$a_{n}\not\asymp_{\N}1$.

We in fact already had an example of this type: in \autoref{SS:Cantor} the Cantor set sequence 
$(\mathrm{ca}_{n})_{n\in\N}$ satisfies
\begin{gather*}
\sum_{n\in\N}\mathrm{ca}_{n}w^{n}
\asymp_{w\uparrow 1}(1-w)^{-\avalue}\big(1+o(1)\big)
\quad\text{and}\quad
\sum_{k=0}^{n}\mathrm{ca}_{k}\asymp n^{\avalue}\big(1+o(1)\big),
\end{gather*}
for $\avalue=\log_{3}2$. However, $\mathrm{ca}_{n}\not\asymp_{\N}n^{\avalue-1}$.
\end{Remark}

\subsection{Conclusion}

We are ready to prove our main result:

\begin{proof}[Proof of \autoref{T:MainTheorem}]
Set $a_{n}:=\frac{b_{n-1}}{2^{n-1}}$ for $n\in\Z_{>0}$. Then \autoref{T:Mono} implies that $a_{n+2}\leq a_{n}$. We also know that $b_{n}\leq b_{n+1}$, since each indecomposable summand in $V^{\otimes n}$ is responsible for at least one 
in $V^{\otimes n+1}$. Hence, $a_{n}\leq 2\cdot a_{n+1}$, so the sequence $a_{n}$ satisfies all conditions in parts (a) and (b) of \autoref{L:Tauber} with $r=2$. 

For the purpose of this proof we will write $\asymp$ to mean $\asymp_{[t,1)}$ for an arbitrary $t\in(0,1)$. Recall that
\begin{gather*}
F(w)=\sum_{n\geq 0}
b_{n}(w+w^{-1})^{-n-1}=H(w+w^{-1}).
\end{gather*}
\autoref{T:AsymGenFunction} implies that
\begin{gather*}
F(w)\asymp (1-w)^{-2\beta},\quad
\beta=\frac{1}{2}\log_{p}\left(\frac{p+1}{2}\right)=\alvalue_{p}+1.
\end{gather*}
Recall further that
\begin{gather*}
H(z)=\sum_{n\geq 0}b_{n}z^{-n-1}\quad\text{and thus}\quad 2H(2z^{-1})=\sum_{n>0}
a_{n}z^{n}.
\end{gather*}
Solving $2z^{-1}=w+w^{-1}$ for $w$ in the region $w<1$ then yields
\begin{gather*}
2H(2z^{-1})\asymp\bigg(\frac{\sqrt{1-z^{2}}+z-1}{z}\bigg)^{-2\beta}
=(1-z)^{-\beta}
\bigg(\frac{z^{2}}{2(1-\sqrt{1-z^{2}})}\bigg)^{\beta}.
\end{gather*}
Since the last factor is bounded on $(0,1]$ away from $0$, we find
\begin{gather*}
\sum_{n>0}a_{n}z^{n}
=2H(\tfrac{2}{z})\asymp(1-z)^{-\beta}.
\end{gather*}
\autoref{L:Tauber}.(b) then implies
\begin{gather*}
a_{n}\asymp_{\Z_{>0}}n^{\beta-1}=n^{\alvalue_{p}},
\end{gather*}
and we are done.
\end{proof}

\section{Additional results in characteristic two}\label{S:PisTwo}

Throughout this section let $p=2$. This restriction is mostly for convenience: with some work the statements and proofs below generalize to all primes.

\subsection{Neighboring values of \texorpdfstring{$b_{n}$}{bn} in characteristic two}\label{SS:CharTwo}

We will now strengthen \autoref{T:Mono}, where we can focus on difference two for the indexes since $b_{2n-1}=b_{2n}$ for $n\in\Z_{>0}$. The auxiliary sequence that we use is $d_{n}=1-\frac{b_{n+2}}{4b_{n}}$, for 
$n\in\N$. 

\begin{Proposition}\label{P:Mono}
We have $d_{n}\geq 0$ and $d_{n}\in O(n^{-\alvalue_{2}-1})$.
\end{Proposition}

Since $-\alvalue_{2}-1\approx -0.293$, \autoref{P:Mono} shows that $\lim_{n\to\infty}d_{n}=0$.

\begin{proof}[Proof of \autoref{P:Mono}]
We start with an analysis of $c_{n}:=b_{n+4}-8b_{n+2}+16b_{n}$.

\begin{Lemma}\label{L:NeighboringLemma}
We have $0\leq c_{n}$ for all $n\in\N$.
\end{Lemma}

\begin{proof}
By \autoref{ExThmFunEq}, we have 
$H(z)=\sum_{n\geq 0}\frac{b_{n}}{2^{n}}z^{-n-1}=z^{-1}+(1+z^{-1})H(z^{2}-2)$.
So, if 
$\eta(z)=(1-4z^{-2})^2H(z)=
z^{-1}+z^{-2}-7z^{-3}-5z^{-4}+\xi(z)$,
then we get
\begin{gather*}
z^{-1}+z^{-2}-7z^{-3}-5z^{-4}+\xi(z)=z^{-1}(1-4z^{-2})^{2}
\\
+z^{-8}(1+z^{-1})\big((z^{2}-2)^{3}+(z^{2}-2)^{2}-7(z^{2}-2)-5+(z^{2}-2)^{4}\xi(z^{2}-2)\big).
\end{gather*}
Hence, we get for $\xi(z)=\sum_{n=0}^{\infty}c_{n}z^{-n-5}$ that: 
\begin{gather*}
\sum_{n=0}^{\infty}c_{n}z^{-n-5}=11z^{-5}+z^{-6}+z^{-7}+5z^{-8}+5z^{-9}+(1+z^{-1})\sum_{n=0}^{\infty}c_{n}z^{-2(n+5)}(1-2z^{-2})^{-n-1}.
\end{gather*}
This gives a positive recursion for $c_{n}$, which implies the statement.
\end{proof}

The numbers $d_{n}$ are nonnegative by \autoref{T:Mono}. Moreover, \autoref{L:NeighboringLemma} implies that $2\leq 1-d_{n+2}+\frac{1}{1-d_{n}}$
which gives
\begin{gather*}
1-d_{n}\leq\frac{1}{1+d_{n+2}}.
\end{gather*}
This yields $1-d_{n-4}\leq\frac{1}{1+d_{n-2}}$ and $1-\frac{1}{1+d_{n}}\leq d_{n-2}$. Putting these together gives $1-d_{n-4}\leq\frac{1}{1+1-\frac{1}{1+d_{n}}}=\frac{1+d_{n}}{1+2d_{n}}$, and iterating this procedure then gives
\begin{gather*}
1-d_{n-2k}\leq\frac{1+(k-1)d_{n}}{1+kd_{n}}.
\end{gather*}
We also have, say for the even values,
\begin{gather*}
(1-d_{0})\cdot\ldots\cdot(1-d_{n-2})=\tfrac{1}{4}\cdot a_{n}
\geq\tfrac{1}{4}C_{1}\cdot n^{\alvalue_{2}}
\end{gather*}
for some $C_{1}\in\R_{>0}$,
where the final inequality is \autoref{T:MainTheorem}.
Thus,
\begin{gather*}
C_{1}^{\prime}\cdot n^{\alvalue_{2}}\leq\frac{1}{1+nd_{n}},
\end{gather*}
for some $C_{1}^{\prime}\in\R_{>0}$,
which, by rewriting, proves the statement.
\end{proof}

\subsection{The main theorem revisited}\label{SS:CharTwoTwo}

The following generalizes \autoref{T:MainTheorem}:

\begin{Proposition}\label{P:MainTheorem}
\autoref{Eq:AlmostConjecture} is true.
\end{Proposition}

\begin{proof}
Let $a_{n}^{\prime}=\frac{b_{2n}}{2^{2n}}$. Note that this ignores the odd values since $a_{n}$ runs only over the even values of $b_{n}$, but since $b_{2n-1}=b_{2n}$ for all $n\in\Z_{>0}$ this does not play a role below and just simplifies the notation.

We observe that $4^{\alvalue_{2}}=4^{-\frac{1}{2}\log_{2}\frac{8}{3}}=\frac{3}{8}$, so that \autoref{T:MainTheorem} implies
that we can sandwich both, $\tfrac{3}{8}\cdot a_{n}^{\prime}$ and $a_{4n}^{\prime}$, at the same time:
\begin{gather*}
C_{1}\cdot\tfrac{3}{8}(2n)^{\alvalue_{2}}\leq
\left\{
\begin{gathered}
\tfrac{3}{8}\cdot a_{n}^{\prime}
\\
a_{4n}^{\prime}
\end{gathered}
\right.
\leq C_{2}\cdot\tfrac{3}{8}(2n)^{\alvalue_{2}}.
\end{gather*}
However, this does not imply the result yet, we need to know a bit more about the sequence $a_{n}^{\prime}$. Firstly, it follows from \autoref{P:Mono} that there is a constant $C\in\R_{>0}$ such that
\begin{gather}\label{bou}
(1-\tfrac{C}{n})a_{n}^{\prime}<a_{n+1}^{\prime}<a_{n}^{\prime}.
\end{gather} 
Also $a_{0}=1$ and, by \autoref{E:CharTwoThing}, for $n>0$ we have
\begin{gather*}
a_{n}^{\prime}=
2^{-n-1}\left(\sum_{k=0}^{\lfloor(n-1)/2\rfloor}\bigg(\binom{n-1}{2k}
+2\binom{n-1}{2k-1}\bigg)a_{k}^{\prime}+\delta_{0,n\bmod 2}\cdot 2a_{n/2}^{\prime}\right).
\end{gather*}
Thus, when running over even values, we get
\begin{gather}\label{bouu}
a_{2n}^{\prime}=
2^{-2n-1}\sum_{k=0}^{n}\bigg(\binom{2n-1}{2k}
+2\binom{2n-1}{2k-1}\bigg)a_{k}^{\prime}.
\end{gather}
Note also that, for $n>0$, we have
\begin{gather*}
2^{-2n-1}\sum_{k=0}^{n}\bigg(\binom{2n-1}{2k}
+2\binom{2n-1}{2k-1}\bigg)=\frac{3}{8}.
\end{gather*}
So using this, combined with
\eqref{bou} and \eqref{bouu}, we get
\begin{gather*}
a_{4n}^{\prime}\sim\tfrac{3}{8}\cdot a_{n}^{\prime}.
\end{gather*}
Hence, for all $x\in\R_{>0}$ there exists a limit
\begin{gather*}
\lim_{m\to\infty}a_{\lfloor 4^{m}x\rfloor}^{\prime}(4^{m}x)^{-\alvalue_{2}}=
\mathbf{f}_{0}(x),
\end{gather*}
and $\mathbf{f}_{0}$ is continuous on $\R_{>0}$ with $\mathbf{f}_{0}(4x)=\mathbf{f}_{0}(x)$, and bounded away from $0$ and $\infty$.
\end{proof}

\section{Fourier coefficients}\label{S:Fourier}

The functions that we have met above are oscillating, and in this section we analyze 
their Fourier coefficients.

\subsection{An Abelian lemma}\label{SS:Abelian}

Tauberian theory as in \autoref{SS:Tauber}, roughly speaking, says that 
given a certain behavior of a generating function, we get a certain behavior
of the associated sequence. There is an ``inverse'' to Tauberian theory often called \emph{Abelian theory}. 

Below we will need the following well-known result from Abelian theory:

\begin{Lemma}\label{L:Abelian}
Assume that we have two functions
$F(q)=\sum_{n=1}^{\infty}a_{n}q^{n}\in\R_{\geq}[[q]]$ and $G(q)=\sum_{n=1}^{\infty}b_{n}q^{n}\in\R_{\geq}[[q]]$
which converge for $q\in[0,1)$ and diverge for $q=1$. Then
\begin{gather*}
\big(a_{n}\sim_{n\to\infty}b_{n}\big)\Rightarrow
\big(F(q)\sim_{q\uparrow 1}G(q)\big).
\end{gather*}
\end{Lemma}

\begin{proof}
This is for example explained at the beginning of \cite[Section 7.5]{MR3155290}.
\end{proof}

\subsection{Some generalities on asymptotics of Fourier coefficients}\label{S:FourierGeneral}

Let $\per\in\Z_{\geq 2}$ and $\alnum>-1$.
Suppose we have a sequence $(a_{n})_{n\in\Z_{>0}}$ 
such that we have asymptotically
\begin{gather*}
a_{n}\sim\htwo(\log_{\per}n)\cdot n^{\alnum}
\end{gather*}
for some continuous $1$-periodic function $\htwo\colon\R\to\R_{>0}$ with $\htwo(x+1)=\htwo(x)$. The associated 
generating function is the series
\begin{gather*}
f(q)=\sum_{n=1}^{\infty}a_{n}q^{n}\in\R[[q]].
\end{gather*}

\begin{Lemma}
The series $f(q)$ absolutely converges for $0\leq|q|<1$ with singularity at $q=1$.
\end{Lemma}

\begin{proof}
This follows since $a_{n}\sim\htwo(\log_{\per}n)\cdot n^{\alnum}$.
\end{proof}

As before, let $\Gamma$ denote the gamma function. Recall 
that the \emph{Fourier coefficient formula} of a $\period$-periodic function $g$
(so that the integral expression up next makes sense) are given by $c_{n}=\frac{1}{\period}\int_{-\period/2}^{\period/2}g(x)e^{-2\pi ixn/\period}dx$, for $n\in\Z$. That is, $g(x)=\sum_{n=-\infty}^{\infty}
c_{n}e^{2\pi ixn/\period}$ for $x\in[-\period/2,\period/2]$.

We get the following asymptotic of $f$:

\begin{Proposition}\label{P:GeneralFourier}
Retain the assumptions above, and denote 
the Fourier coefficients of $h$ by $h_{n}$.
We have
\begin{gather*}
f(q)\sim_{q\uparrow 1}\big(\ln(q^{-1})\big)^{-\alnum-1}L\big(\log_{\per}(\ln(q^{-1}))\big)
\end{gather*}
where $L$ is the $1$-periodic function given by
\begin{gather*}
L(y)=
\sum_{n\in\Z}\Gamma(\alnum+1-\tfrac{2\pi in}{\log\per})e^{2\pi iny}h_{-n}.
\end{gather*}
\end{Proposition}

\begin{proof}
For the argument below we note that $q=1$ is the relevant singularity of
$f$, and we need to analyze the growth rate of $f$ for $q\uparrow 1$. 

Recalling that $a_{n}\sim\htwo(\log_{\per}n)\cdot n^{\alnum}$, we use \autoref{L:Abelian} and get
\begin{gather*}
f(q)\sim_{q\uparrow 1}
\sum_{n=1}^{\infty}\htwo(\log_{\per}n)n^{\alnum}q^{n}.
\end{gather*}
We rewrite this as
\begin{gather*}
f(q)\sim_{q\uparrow 1}\sum_{r\in\N}
\sum_{1\leq\per^{-r}n<\per}\htwo(\log_{\per}n)n^{\alnum}q^{n}.
\end{gather*}
So setting $x=\per^{-r}n$ and using the periodicity of $\htwo$, namely $\htwo(x+1)=\htwo(x)$, we get 
\begin{gather*}
f(q)\sim_{q\uparrow 1}\sum_{r\in\N}\per^{r\cdot \alnum}
\sum_{1\leq x<\per}\htwo(\log_{\per}x)x^{\alnum}q^{\per^{r}x},
\end{gather*}
where the inner sum is over (not necessarily reduced) 
fractions with denominator $\per^{r}$. Note that
the main contribution to the asymptotics $q\uparrow 1$ comes from large $r$, as $\alnum>-1$. 
Hence, we may replace the inner sum by an integral times $\per^{r}$ (reciprocal length of step), and the range $\N$ of summation for the outer sum by $\Z$. After doing this we get:
\begin{gather*}
f(q)\sim_{q\uparrow 1}\sum_{r\in\Z_{\ge 0}}\per^{r(\alnum+1)}\int_{1}^{\per}
\htwo(\log_{\per}x)x^{\alnum}q^{\per^{r}x}dx.
\end{gather*}
Making the substitution $u=\log_\nu(x)$ allows us to rewrite the above as
$$\ln(\nu)\sum_{r\in\Z_{\ge 0}}\int_0^1\nu^{(u+r)(\alnum+1)}q^{\nu^{r+u}}\tilde{h}(u)du\;=\; \ln(\nu)\int^\infty_0 \nu^{u(\alnum+1)}q^{\nu^u}\tilde{h}(u)du.$$
Substituting back to $x$, the above yields
\begin{gather*}
f(q)\sim_{q\uparrow 1}\int_0^\infty x^{\alnum}q^x \tilde{h}(\log_{\nu}(x))dx.
\end{gather*}
Replacing the dummy variable $x$ by $x/\ln(q^{-1})$ yields
\begin{gather*}
f(q)\sim_{q\uparrow 1}(\ln(q^{-1}))^{-\alnum-1}\int_0^\infty x^{\alnum}e^{-x}\tilde{h}\left(\log_{\nu}(x)-\log_\nu(\ln(q^{-1}))\right)dx.
\end{gather*}
We can thus define the 1-periodic function
$$L(y)\;:=\;\int_0^\infty x^{\alnum}e^{-x}\tilde{h}\left(\log_{\nu}(x)-y\right)dx.$$
That this periodic function corresponds to the one in the proposition follows easily from 
the definition $\Gamma(z)=\int_{0}^{\infty}x^{z-1}e^{-x}dx$ of the Gamma function.
\end{proof}

\begin{Theorem}\label{T:GeneralFourier}
With the notation as in \autoref{P:GeneralFourier}, we have
\begin{gather*}
h_{n}=\frac{L_{n}}{\Gamma(\alnum+1+\frac{2\pi in}{\log\per})}
\end{gather*}
for the Fourier coefficients $L_{n}$ of $L(y)$.
\end{Theorem} 

\begin{proof}
By \autoref{P:GeneralFourier} and Fourier inversion.
\end{proof}

\subsection{Analysis of the Fourier coefficients of the generating function for \texorpdfstring{$b_{n}/2^{n}$}{bn/2n}}\label{SS:CharTwoFourier}

Assume that \autoref{Eq:AlmostConjecture} holds. (For $p=2$ this assumption is true by \autoref{S:PisTwo}.)
We now apply the results in
\autoref{S:FourierGeneral} to our example of 
$b_{n}/2^{n}$ in characteristic $p$. 
That is, we look at $a_{n}=b_{n}/2^{n}$ for the $a_{n}$ in \autoref{S:FourierGeneral}.

Consider the generating function of $a_{n}$:
\begin{gather*}
f(q)=\sum_{n=0}^{\infty}a_{n}q^{n}\in\R_{>0}[[q]].
\end{gather*}
The relevant singularity is at $q=1$.
Below we will also use the generating function in $w$ defined mutatis mutandis as in \autoref{SubSecFunF}.

Moreover, associated to $f$ we can define 
$\mathbf{f}_{0}$ similarly as $\mathbf{F}_{0}$
from \autoref{SS:ALimit} but defined now in $q$, 
and we get the same type of statements 
for $\mathbf{f}_{0}$ as for $\mathbf{F}_{0}$.

Recall $\alvalue_{p}=\frac{1}{2}\log_{p}\frac{2p^{2}}{p-1}$ from \autoref{Eq:Alpha}. 
Setting $\per=2p$ (which we expect to be correct based on numerical data) and $\alnum=\alvalue_{p}$, we get:

\begin{Proposition}\label{P:FGrowth}
We have
\begin{gather*}
f\sim_{q\uparrow 1}\big(\ln(q^{-1})\big)^{-\alvalue_{p}-1}\mathbf{f}_{0}\big(\log_{2p}\ln(q^{-1})\big)
\end{gather*}
and the Fourier coefficients of $\mathbf{f}_{0}$ are given by 
\autoref{E:GrowthCompareNumerical}.
\end{Proposition}

\begin{proof}
Recall from \autoref{SS:ALimit} that 
$\lim_{r\to\infty}(\ln w^{-p^{-r}})^{\log_{p}(\frac{p+1}{2})}F(w^{p^{-r}})=\mathbf{F}_{0}(w)$. 
We can write this as
$\lim_{q\uparrow 1}\big(\ln(q^{-1})\big)^{-2(\alvalue_{p}+1)}f(q)=\mathbf{f}_{0}\big(\log_{2p}\ln(q^{-1})\big)$.
Then \autoref{P:GeneralFourier} as well as \autoref{T:GeneralFourier} apply.
\end{proof}

\section{Numerical data}\label{S:Data}

All of the below can be found on the GitHub page associated to this project 
\cite{CoEtOsTu-growth-code}, 
where the reader can find code that can be run online. For convenience, we list a 
few numerical examples here as well.

\begin{Remark}
All plots below are \emph{logplots} (also called semi-log plots). This means that the plots have one axis (the y-axis for us) on a logarithmic scale, the other on a linear scale.
\end{Remark}

\subsection{The main sequence}\label{SS:MainData}

Below we will always use $p\in\{2,3,5,7,\infty\}$ 
with $p=\infty$ meaning that the characteristic is large compared to the considered range, which leads to the same data as in characteristic zero.

We first give tables for $b_{n}$ for $n\in\{0,\dots,15\}$:

\begin{center}
\begin{tabular}{c|cccccccccccccccc}
& $b_{0}$ & $b_{1}$ & $b_{2}$ & $b_{3}$ & $b_{4}$ & $b_{5}$ & $b_{6}$ & $b_{7}$ & $b_{8}$ & $b_{9}$ & $b_{10}$ & $b_{11}$ & $b_{12}$ & $b_{13}$ & $b_{14}$ & $b_{15}$ \\
\hline
\rowcolor{orchid!25}
$p=2$ & 1 & 1 & 1 & 3 & 3 & 9 & 9 & 29 & 29 & 99 & 99 & 351 & 351 & 1273 & 1273 & 4679 \\
\rowcolor{lava!10}
$p=3$ & 1 & 1 & 2 & 2 & 5 & 6 & 15 & 21 & 50 & 77 & 176 & 286 & 637 & 1066 & 2340 & 3978 \\
\rowcolor{orchid!25}
$p=5$ & 1 & 1 & 2 & 3 & 6 & 9 & 19 & 28 & 62 & 91 & 208 & 308 & 716 & 1079 & 2522 & 3886 \\
\rowcolor{lava!10}
$p=7$ & 1 & 1 & 2 & 3 & 6 & 10 & 20 & 34 & 69 & 117 & 242 & 407 & 858 & 1431 & 3069 & 5085 \\
\rowcolor{orchid!25}
$p=\infty$ & 1 & 1 & 2 & 3 & 6 & 10 & 20 & 35 & 70 & 126 & 252 & 462 & 924 & 1716 & 3432 & 6435 \\
\end{tabular}
.
\end{center}

For completeness, we give here the $b_{n}$ for $n\in\{0,\dots,30\}$ (in order top to bottom, excluding $p=\infty$ where the sequence is \cite[A001405]{oeis}) in a copy-able fashion:

\resizebox{0.95\textwidth}{!}{$\{1,1,1,3,3,9,9,29,29,99,99,351,351,1273,1273,4679,4679,17341,17341,64637,64637,242019,242019,909789,909789,3432751,3432751,12998311,12998311,49387289,49387289\}$.}

\resizebox{0.95\textwidth}{!}{$\{1,1,2,2,5,6,15,21,50,77,176,286,637,1066,2340,3978,8670,14859,32301,55575,120822,208221,453399,781794,1706301,2942460,6438551,11103665,24357506,42015664,92376280\}$.}

\resizebox{0.95\textwidth}{!}{$\{1,1,2,3,6,9,19,28,62,91,208,308,716,1079,2522,3886,9061,14297,33098,53448,122551,202181,458757,771443,1732406,2962284,6587959,11428743,25193027,44250404,96775581\}$.}

\resizebox{0.95\textwidth}{!}{$\{1,1,2,3,6,10,20,34,69,117,242,407,858,1431,3069,5085,11066,18258,40205,66215,147136,242420,542202,895390,2011165,3334125,7505955,12507121,28174255,47229893,106315770\}$.}

\noindent
Here is one picture to compare their growth, including the case $p=\infty$:
\begin{gather*}
\begin{tikzpicture}[anchorbase]
\node at (0,0) {\includegraphics[height=5.0cm]{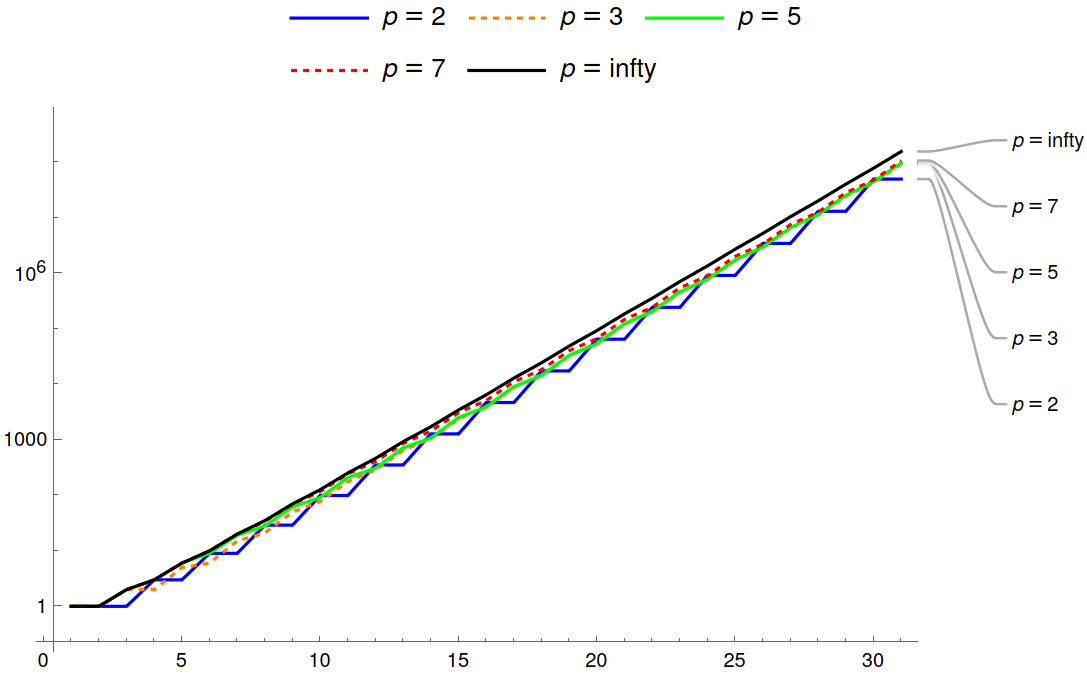}};
\end{tikzpicture}
.
\end{gather*}
Note that the sequences get a bit more regular if one runs over the, say, even values. This is in particular true in characteristic $p=2$ where $b_{2n-1}=b_{2n}$. Below we will strategically sometimes only illustrates the even values.

Next, the logplots of $b_{n}/2^{n}$ compared with $n^{\alvalue_{p}}$, and of $b_{n}$ compared with $n^{\alvalue_{p}}\cdot 2^{n}$ are:
\begin{gather*}
\begin{tikzpicture}[anchorbase]
\node at (0,0) {\includegraphics[height=4.6cm]{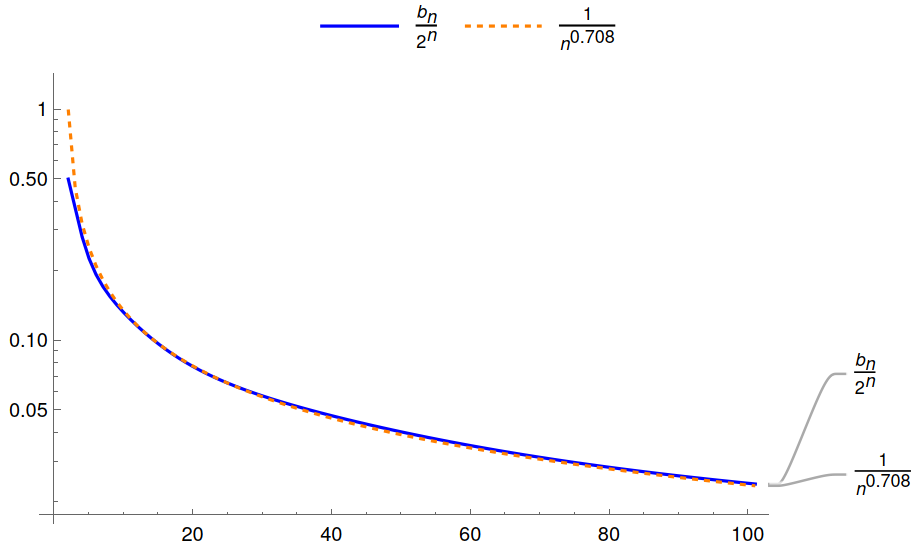}};
\node at (0,0) {$p=2$,even};
\end{tikzpicture}
,\quad
\begin{tikzpicture}[anchorbase]
\node at (0,0) {\includegraphics[height=4.6cm]{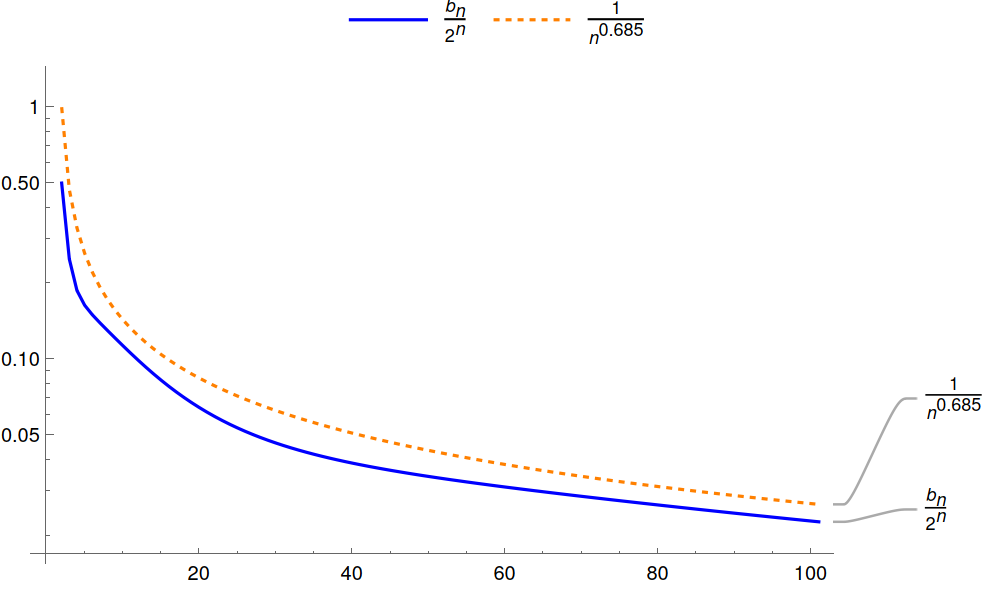}};
\node at (0,0) {$p=3$,even};
\end{tikzpicture}
,
\\
\begin{tikzpicture}[anchorbase]
\node at (0,0) {\includegraphics[height=4.6cm]{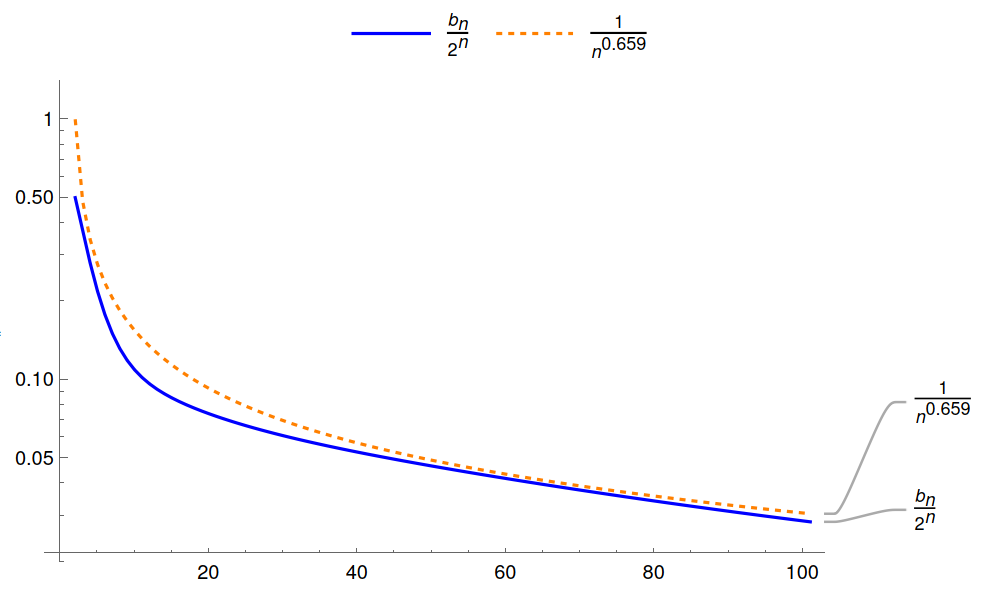}};
\node at (0,0) {$p=5$,even};
\end{tikzpicture}
,\quad
\begin{tikzpicture}[anchorbase]
\node at (0,0) {\includegraphics[height=4.6cm]{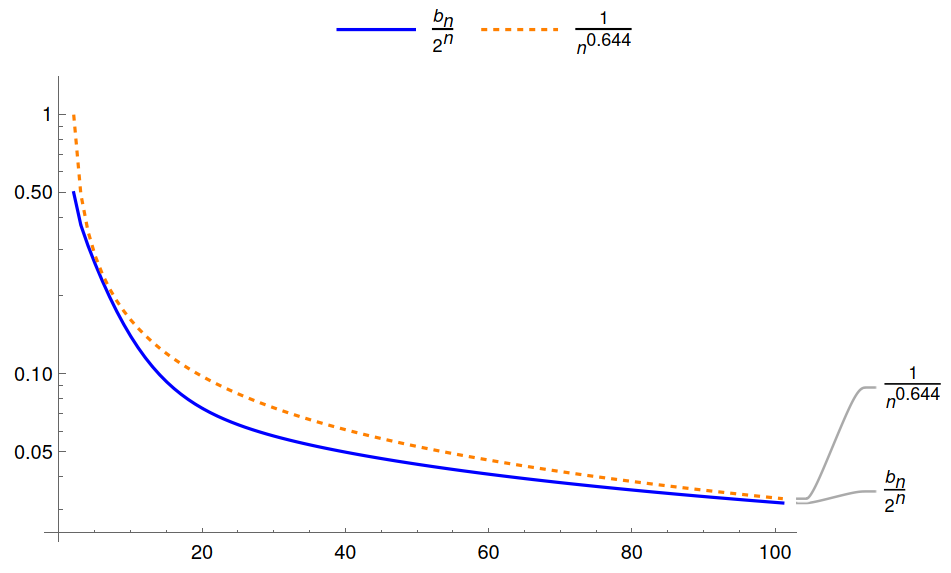}};
\node at (0,0) {$p=7$,even};
\end{tikzpicture}
,
\\
\begin{tikzpicture}[anchorbase]
\node at (0,0) {\includegraphics[height=4.6cm]{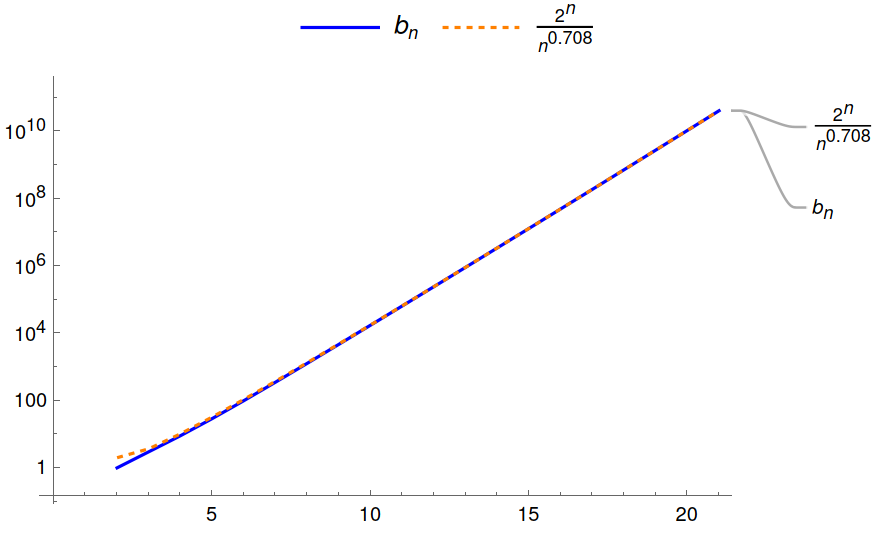}};
\node at (-1.5,0) {$p=2$,even};
\end{tikzpicture}
,\quad
\begin{tikzpicture}[anchorbase]
\node at (0,0) {\includegraphics[height=4.6cm]{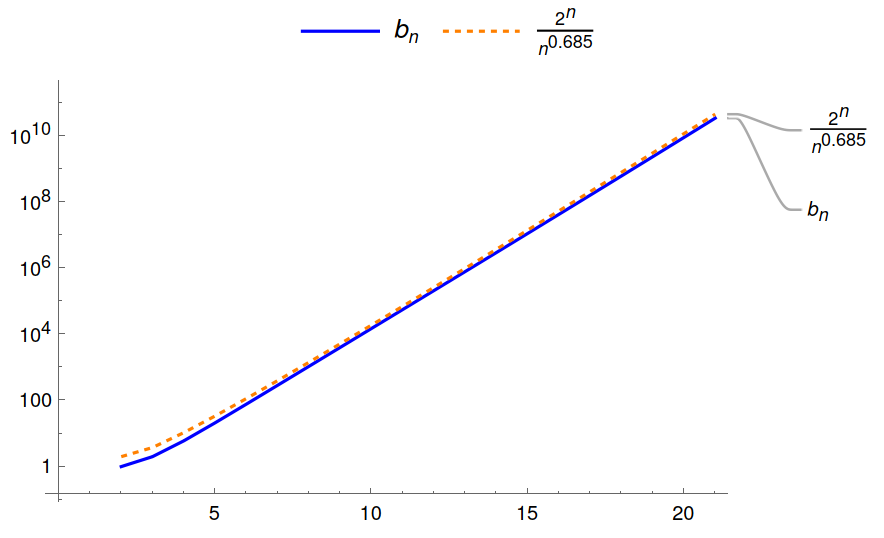}};
\node at (-1.5,0) {$p=3$,even};
\end{tikzpicture}
,
\\
\begin{tikzpicture}[anchorbase]
\node at (0,0) {\includegraphics[height=4.6cm]{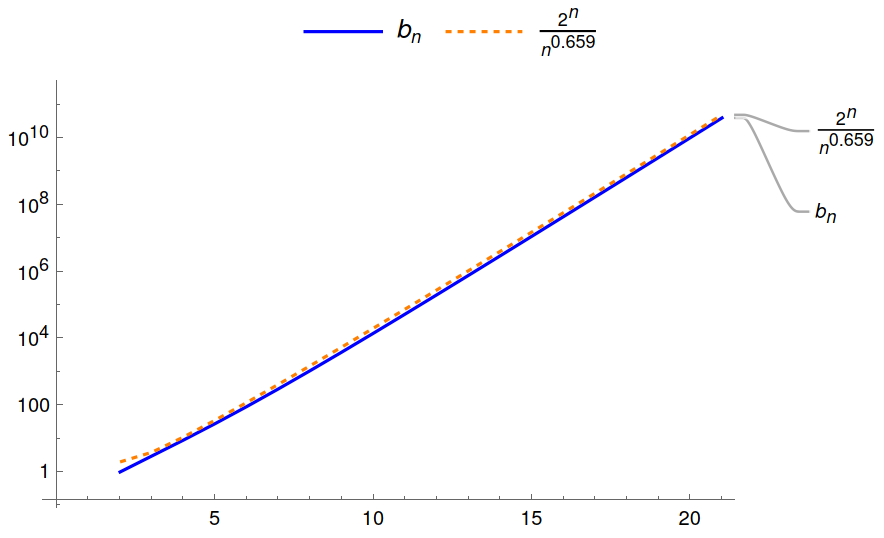}};
\node at (-1.5,0) {$p=5$,even};
\end{tikzpicture}
,\quad
\begin{tikzpicture}[anchorbase]
\node at (0,0) {\includegraphics[height=4.6cm]{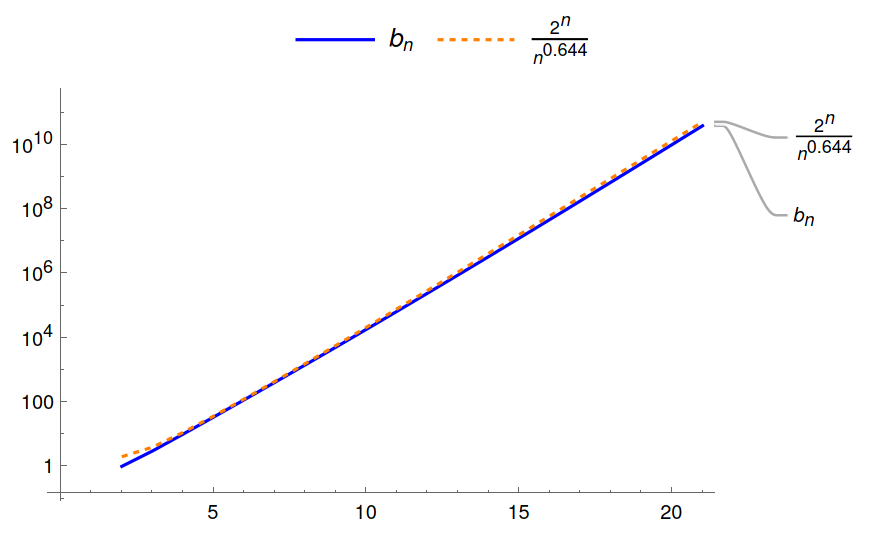}};
\node at (-1.5,0) {$p=7$,even};
\end{tikzpicture}
.
\end{gather*}

\subsection{The length sequence}

Let us now also give here the $l_{n}$ for $n\in\{0,\dots,30\}$ (again in order top for $p=2$ to bottom for $p=7$, excluding $p=\infty$ where the sequence is \cite[A001405]{oeis}) so that it is copy-able:

\resizebox{0.95\textwidth}{!}{$\{1,1,3,3,11,11, 41, 41, 155, 155, 593, 593, 2289, 2289, 8891,8891, 34683, 34683, 135697, 135697, 532041, 532041, 2089363, 2089363,8215553, 8215553, 32339011, 32339011, 127417011, 127417011, 502458289\}$.}

\resizebox{0.95\textwidth}{!}{$\{1, 1, 2, 4, 7, 14, 26, 50, 97, 184, 364, 692, 1378, 2641, 5264,10181, 20267, 39523, 78524, 154187, 305728, 603614, 1194758, 2368906,4682134, 9313411, 18387902, 36663241, 72331456, 144466892, 28488466\}$.}

\resizebox{0.95\textwidth}{!}{$\{1, 1, 2, 3, 6, 11, 21, 42, 78, 161, 297, 617, 1144, 2366, 4432, 9088,17223, 34986, 67049, 135013, 261326, 522271, 1019427, 2024828,3979781, 7866186, 15547861, 30614847, 60783158, 119345091, 237790431\}$.}

\resizebox{0.95\textwidth}{!}{$\{1, 1, 2, 3, 6, 10, 20, 36, 71, 135, 262, 517, 990, 2001, 3796, 7786,14690, 30379, 57188, 118712, 223515, 464341, 875955, 1817598,3439375, 7119305, 13522875, 27902564, 53222511, 109424657, 209629719\}$.}

\noindent
As before for $b_{n}$, here is one picture to compare their growth, including the case $p=\infty$:
\begin{gather*}
\begin{tikzpicture}[anchorbase]
\node at (0,0) {\includegraphics[height=5.0cm]{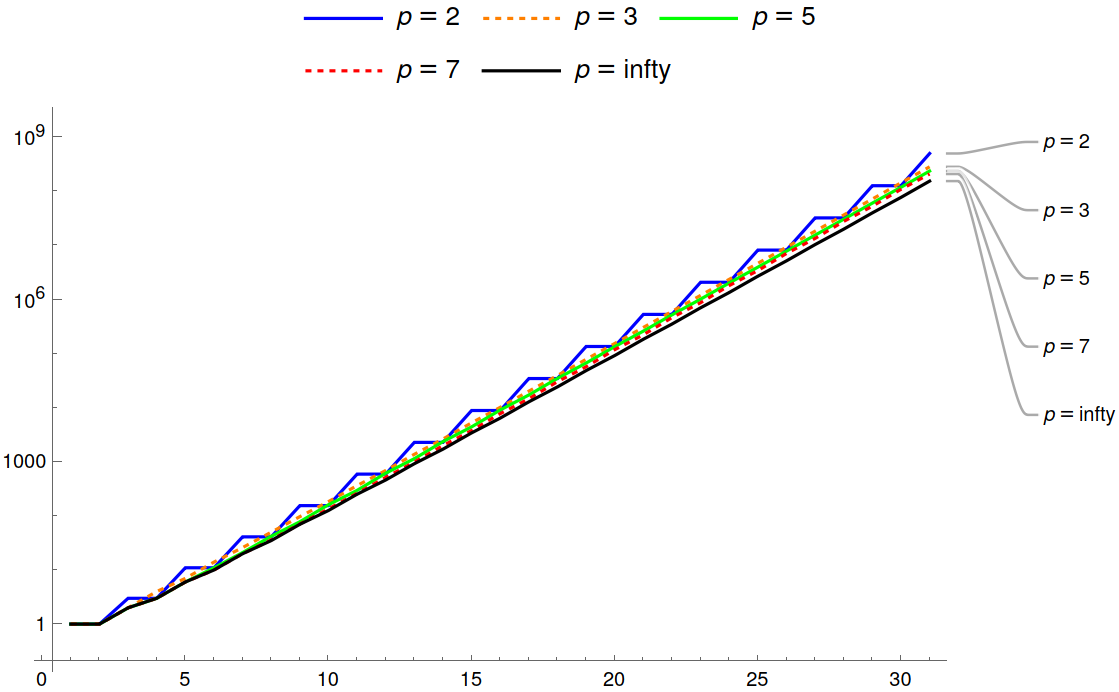}};
\end{tikzpicture}
.
\end{gather*}
Note the following comparison of the $l_{n}$ plot to the $b_{n}$ plot from above: for $p<p^{\prime}$, the $l_{n}$ grows a faster for $p$ and slower for 
$p^{\prime}$, and vice versa for the $b_{n}$.

\subsection{The generating function}\label{E:GrowthCompareNumerical}

Continuing with \autoref{E:GrowthCompare}, recall from \autoref{Eq:Alpha} that $2(\alvalue_{2}+1)\approx0.585$.
We compare $F(w)$, with the sum cut-off at $k=10$, with $(w-1)^{-0.585}$:
\begin{gather*}
\begin{tikzpicture}[anchorbase]
\node at (0,0) {\includegraphics[height=4cm]{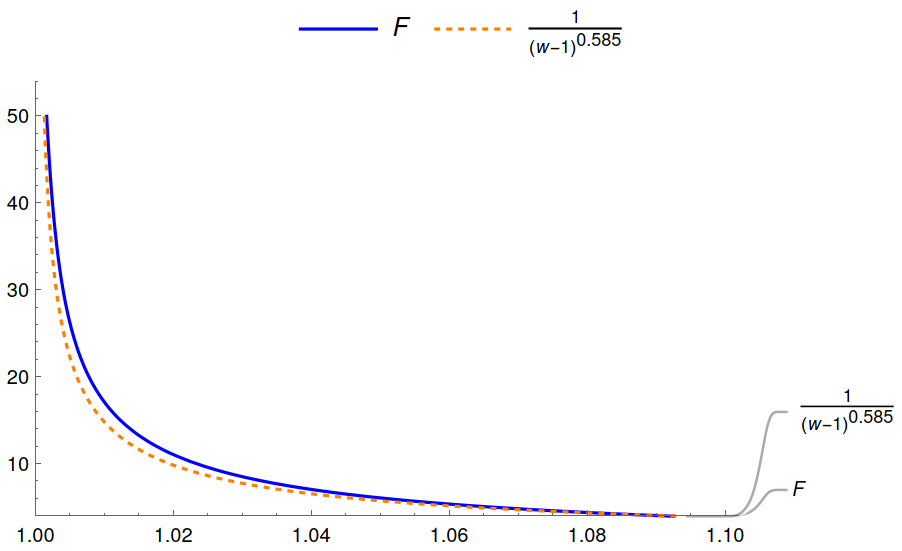}};
\end{tikzpicture}
\end{gather*}
The growth rates towards the singularity $w=1$ (left side of the picture) are almost the same and
this is crucial in \autoref{S:MainTheorem}.

Recall from \autoref{SS:ALimit} that
\begin{gather*}
	\mathbf{F}_{0}(w)=\big(\ln(w^{-1})\big)^{\log_{p}(\frac{p+1}{2})}\sum_{k=-\infty}^{\infty}\frac{w^{p^{k}}(1-w^{p^{k}(p-1)})}{(1-w^{p^{k}})(1+w^{p^{k+1}})}\Pi(w^{p^{k}})\left(\frac{p+1}{2}\right)^{k}.
\end{gather*}
We now give formulas to compute the Fourier coefficients $L_{n}$ of $\mathbf{f}_{0}$, and thus of $\mathbf{F}_{0}$, (fairly) efficiently.
It turns out that $L_{0}$ is moderate but $L_{n}$ are tiny for $n\neq 0$, 
which causes the behavior of $\mathbf{F}_{0}(2^{-p^{x}})$ as in the next example, plotted on $[0.1,0.9]$:
\begin{gather*}
\begin{tikzpicture}[anchorbase]
\node at (0,0) {\includegraphics[height=4.8cm]{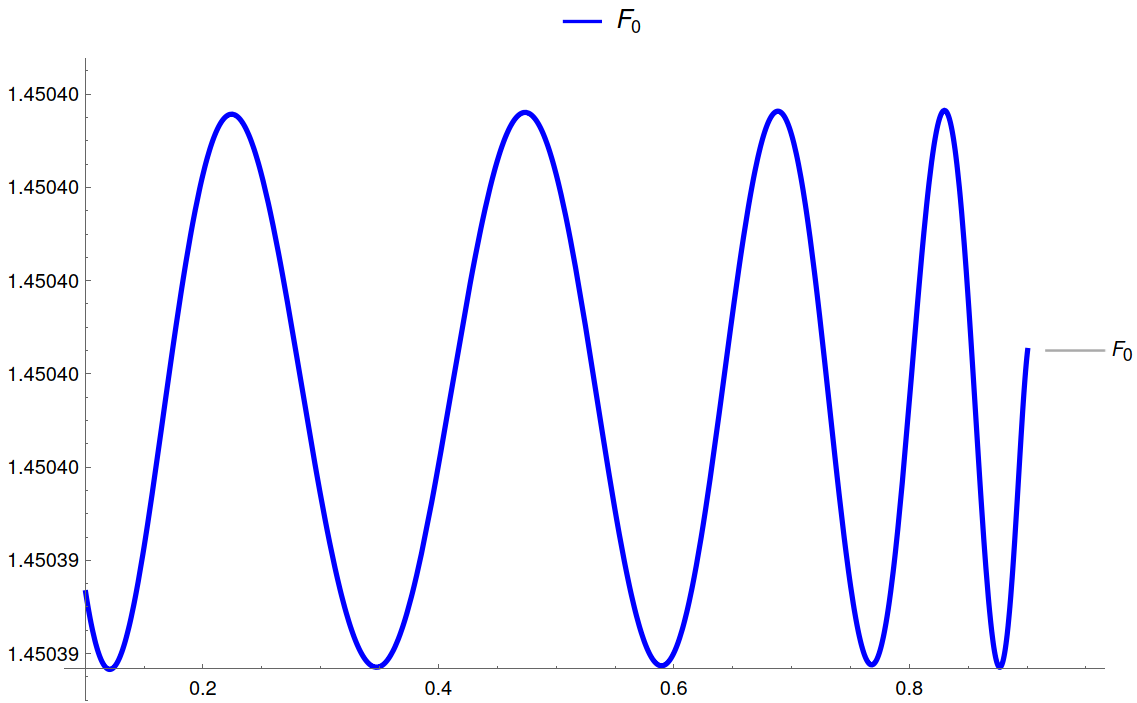}};
\node at (0.6,0.75) {$p=2$};
\end{tikzpicture}
.
\end{gather*}
This looks like the sine curve because $\sum_{n>1}|L_{n}|$ is much smaller than $|L_{1}|$ as we will see momentarily.
To wrap-up, for $p=2$ one can show that
\begin{gather*}
	\scalebox{0.95}{$\displaystyle L_{n}=\frac{1}{(\ln 2)^{2-\log_{2}(3/2)}}
		\int_{0}^{1}\frac{\ln(w^{-1})^{\log_{2}(\frac{3}{2})+in\log_{2}(e^{2\pi})}}{w+w^{-1}}\prod_{j=1}^{\infty} 
		\left(1-\frac{(w^{2^{-j-1}}-w^{-2^{-j-1}})^2}{3(w^{2^{-j}}+w^{-2^{-j}})}\right)\frac{dw}{w\log_{2}w^{-1}}.$}
\end{gather*}
Let $\zeta(x,s)=\sum_{k=0}^{\infty}(k+x)^{-s}$ denote the 
\emph{Hurwitz zeta function}.
Consider the following reparametrization:
\begin{gather*}
	\xi(\beta,u)=4^{-\beta}\zeta(\beta,\tfrac{u}{4}).
\end{gather*}
Recalling that $\Gamma(z)=\int_{0}^{\infty}x^{z-1}e^{-x}dx=\int_{0}^{1}(\ln\frac{1}{x})^{z-1}dx$, we get 
\begin{gather*}
	L_{n}=\frac{2\Gamma\big(\log_{2}(\tfrac{3}{2})+1+\tfrac{2\pi i n}{\ln 2}\big)}{(\ln 2)^{1+
			\frac{2\pi in}{\ln 2}}}\lim_{N\to\infty}\frac{1}{3^{N}}
	\sum_{-1\leq k_{1},\dots,k_{N}\leq 1}\xi\big(\log_{2}(\tfrac{3}{2})+1+\tfrac{2\pi i n}{\ln 2},2+{\textstyle\sum_{j=1}^{N}}k_{j}2^{-j}\big).
\end{gather*} 
The Hurwitz zeta function is quite easy to compute, so this gives 
a fairly efficient way to compute the Fourier coefficients $L_{n}$.

A bit of work shows that, for general $p$, we have:
\begin{gather*}
	\scalebox{0.95}{$L_{n}=\frac{2\Gamma(\log_{p}(\frac{p+1}{2})+1+
			\frac{2\pi in}{\ln p})}{(\ln p)}\lim_{N\to\infty}\frac{1}{(p+1)^{N}}
		\sum_{-1\leq k_{1},\dots,k_{N}\leq 1}\xi(\log_{p}(\tfrac{p+1}{2})+1+\tfrac{2\pi i n}{\ln p},2+{\textstyle\sum_{j=1}^{N}}k_{j}p^{-j}).$}
\end{gather*}
Moreover, using the classical formulas for the Hurwitz zeta function one can show
that $|L_{n}|$ is nonzero. Furthermore, let
\begin{gather*}
	S(p,n,N):=
	\bigg|\sum_{-1\leq k_{1},\dots,k_{N}\leq 1}\xi(\log_{p}(\tfrac{p+1}{2})+1+\tfrac{2\pi i n}{\ln p},2+{\textstyle\sum_{j=1}^{N}}k_{j}p^{-j})\bigg|.
\end{gather*}
One can also show that $\max\{S(p,n,N)\}$,
for fixed $N$, is obtained at $p=2$ and $n=0$. Indeed, the function $S(p,n,N)$ behaves as follows when varying $n$ and $p$, respectively, while keeping $N$ fixed: 
\begin{gather*}
\begin{tikzpicture}[anchorbase]
\node at (0,0) {\includegraphics[height=4.8cm]{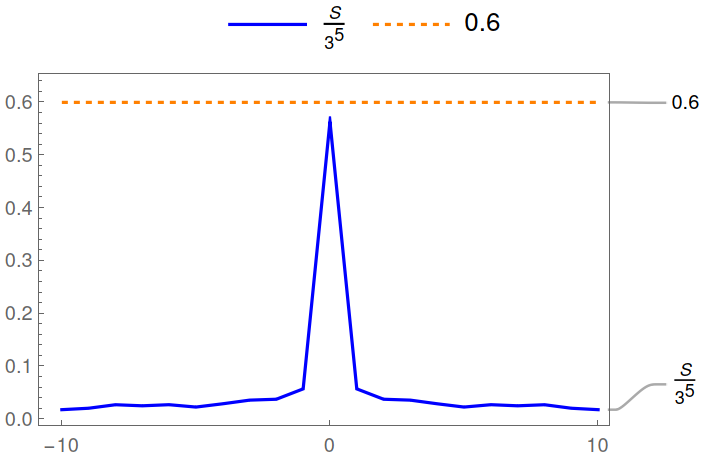}};
\node at (1.5,0) {$p=2,N=5$};
\end{tikzpicture}
,\quad
\begin{tikzpicture}[anchorbase]
\node at (0,0) {\includegraphics[height=4.8cm]{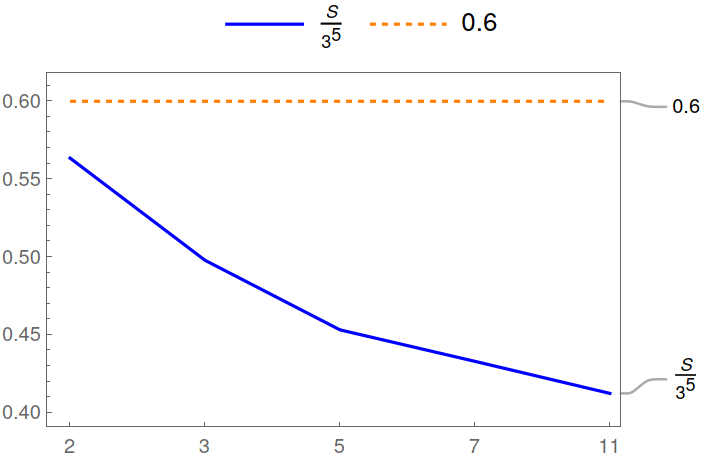}};
\node at (1.5,0) {$n=0,N=5$};
\end{tikzpicture}
.
\end{gather*}
Here we have divided by $3^{N}$ 
to highlight the bound $0.6$ which is quite crude:
\begin{gather*}
\begin{tikzpicture}[anchorbase]
\node at (0,0) {\includegraphics[height=4.8cm]{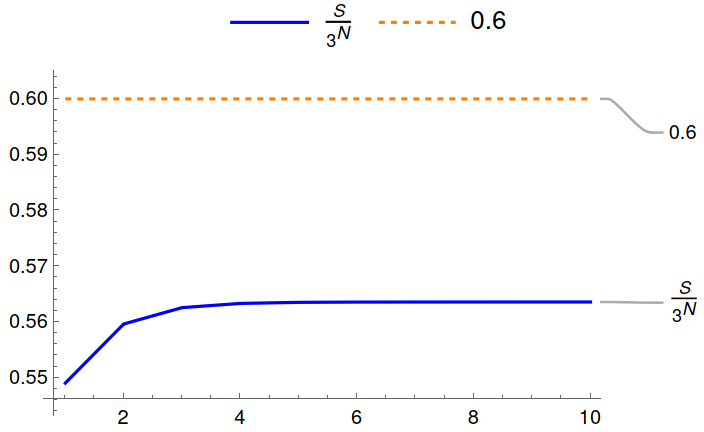}};
\node at (1.5,0) {$p=2,n=0$};
\end{tikzpicture}
.
\end{gather*}

\subsection{Numerical data for odds and ends}\label{SS:ScalarsMainTheorem}

Recall that \autoref{T:MainTheorem} shows that, for some $C_{1},C_{2}\in\R_{>0}$, we have
\begin{gather*}
C_{1}\cdot n^{\alvalue_{p}}\cdot 2^{n}\leq b_{n}\leq C_{2}\cdot n^{\alvalue_{p}}\cdot 2^{n},\qquad n\geq 1.
\end{gather*}
One can probably take $C_{1}=1/4$ and $C_{2}=1$, as illustrated here for $p=2$:
\begin{gather*}
\begin{tikzpicture}[anchorbase]
\node at (0,0) {\includegraphics[height=5.0cm]{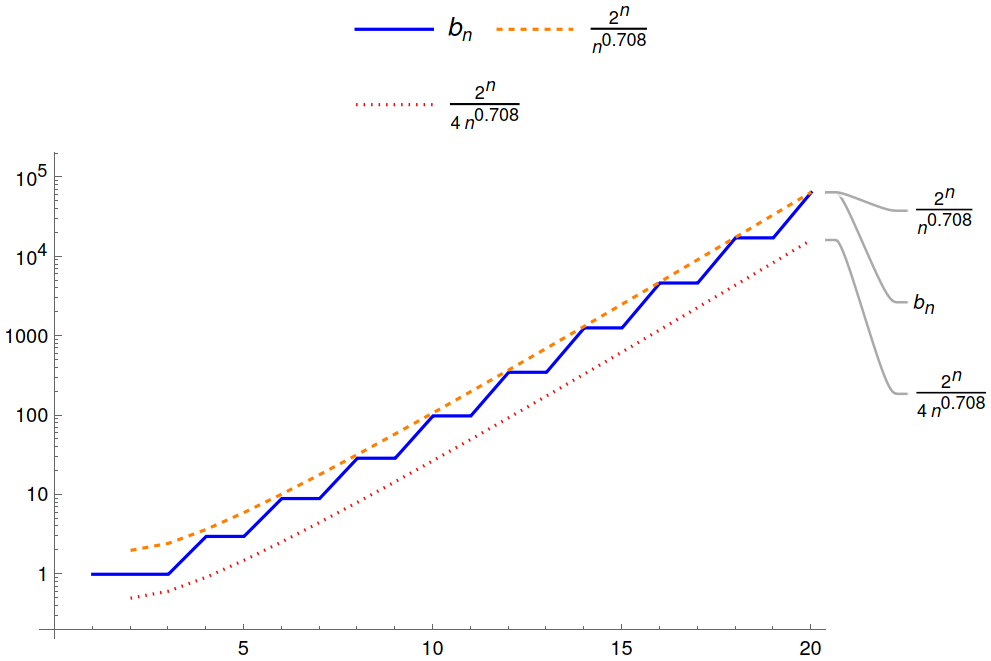}};
\node at (-1.5,0) {$p=2$};
\end{tikzpicture}
,
\begin{tikzpicture}[anchorbase]
\node at (0,0) {\includegraphics[height=5.0cm]{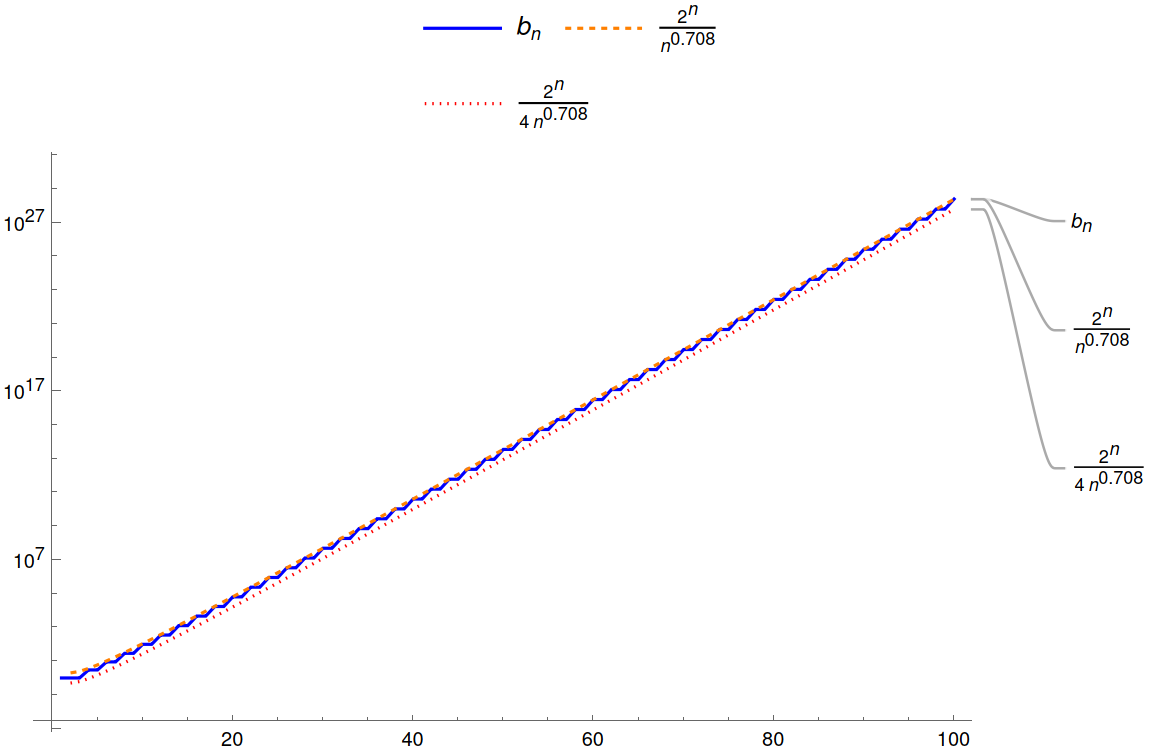}};
\node at (-1.5,0) {$p=2$};
\end{tikzpicture}
.
\end{gather*}
Moreover, \autoref{T:Mono} shows that we have
\begin{gather*}
b_{n+2}\leq 4b_{n},\; n\geq 0.
\end{gather*}
One could conjecture that $\lim_{n\to\infty}b_{n+2}/b_{n}=4$, as illustrated here for $p=2$ where we know this is true by \autoref{P:Mono}:
\begin{gather*}
\begin{tikzpicture}[anchorbase]
\node at (0,0) {\includegraphics[height=4.55cm]{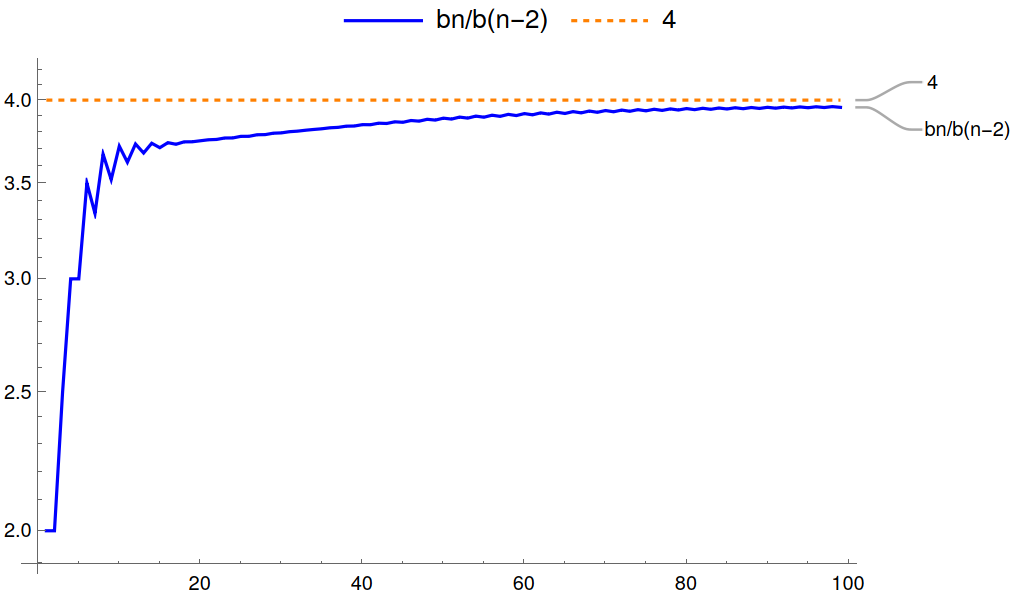}};
\node at (-1.5,0) {$p=2$};
\end{tikzpicture}
\hspace{-0.5cm},
\begin{tikzpicture}[anchorbase]
\node at (0,0) {\includegraphics[height=4.55cm]{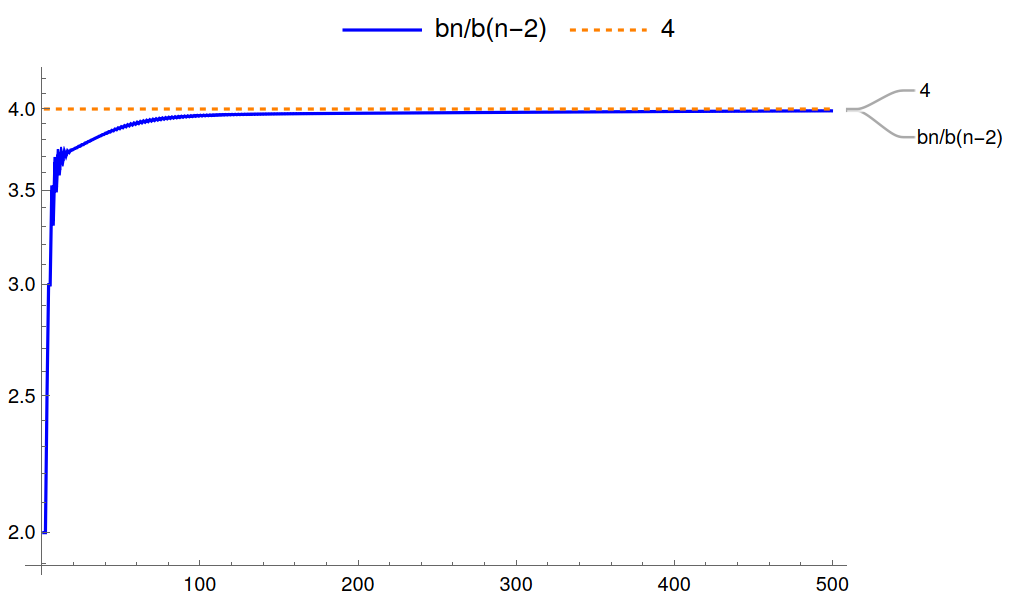}};
\node at (-1.5,0) {$p=2$};
\end{tikzpicture}
.
\end{gather*}
Finally, $W=\mathrm{Sym}^{2}V$ 
with $\dim_{\mathbf{k}}W=3$ is an indecomposable tilting
$SL_{2}(\mathbf{k})$-representation unless $p=2$.
For $p=3$ we get for $b_{n}=b_{n}(W)$ the numbers:
\begin{align*}
(b_{n})_{n\in\N}=(&1,1,2,5,13,35,95,260,715,1976,5486,15301,
\\
&42686,120628,340874,967136,2754455,7872973,\dots).
\end{align*}
One could aim for a statement of the form
\begin{gather*}
D_{1}\cdot n^{\alvalue_{p}}\cdot 3^{n}\leq b_{n}\leq D_{2}\cdot n^{\alvalue_{p}}\cdot 3^{n},\qquad n\geq 1,
\end{gather*}
where $D_{1},D_{2}\in\R_{>0}$. Indeed,
we get the following logplots:
\begin{gather*}
\begin{tikzpicture}[anchorbase]
\node at (0,0) {\includegraphics[height=4.8cm]{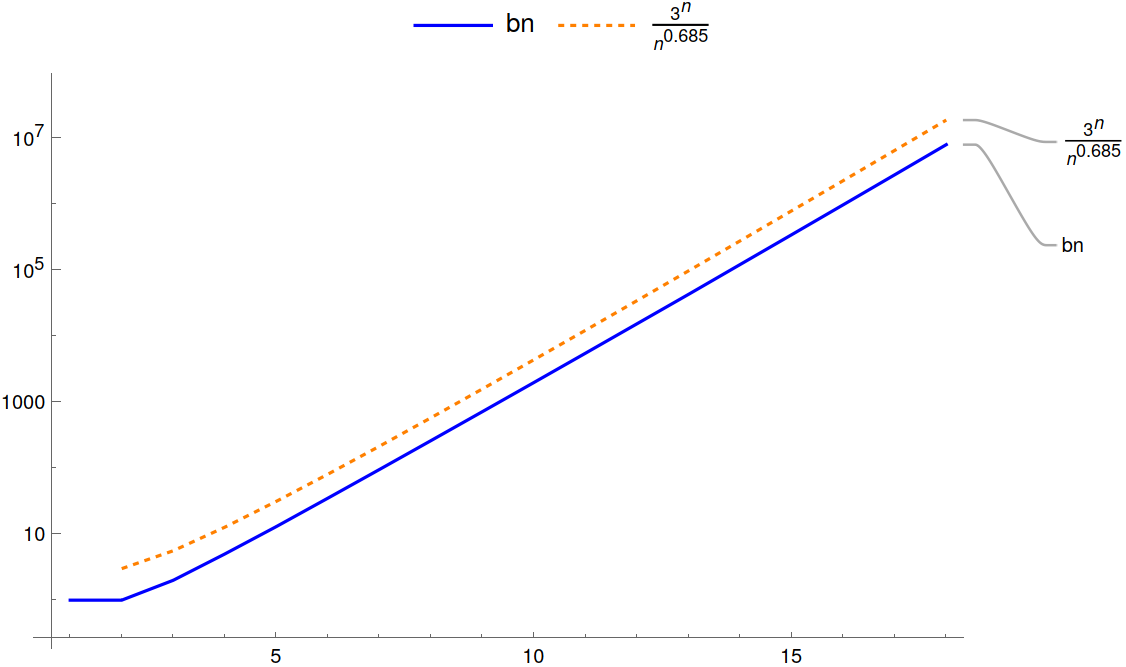}};
\node at (1.5,0) {$p=3$};
\end{tikzpicture}
\hspace{-0.1cm},
\begin{tikzpicture}[anchorbase]
\node at (0,0) {\includegraphics[height=4.8cm]{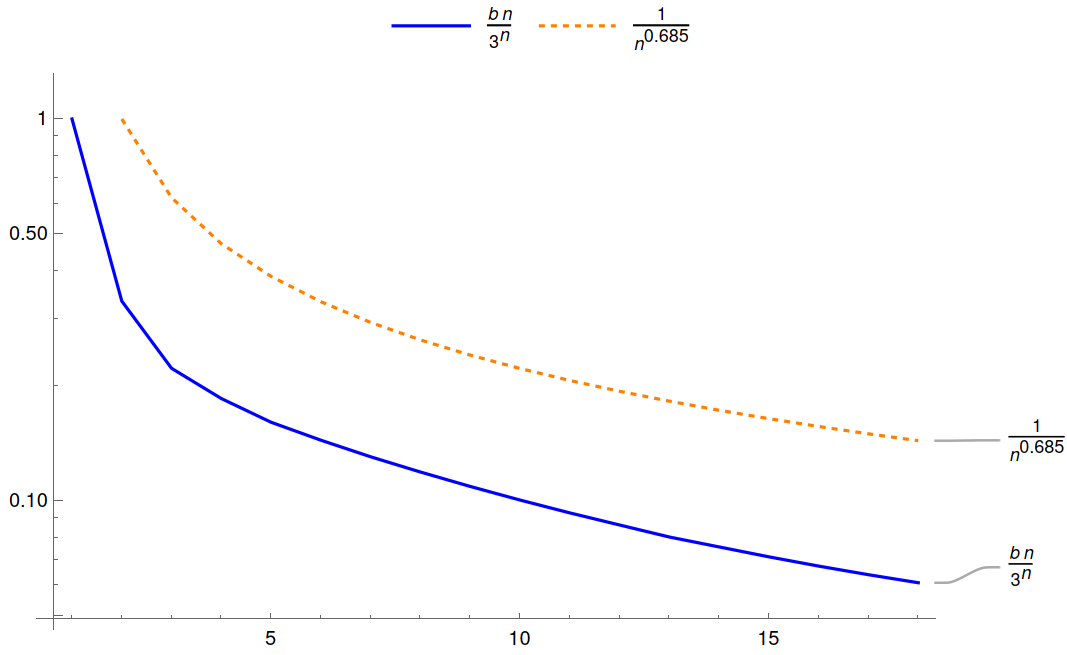}};
\node at (1.5,0) {$p=3$};
\end{tikzpicture}
.
\end{gather*}
As one can see, we expect some nontrivial scalar to make the curves fit better.
Hence, these plots might indicate 
an analog of \autoref{T:MainTheorem} for $W$ instead of $V$.

\appendix

\section{Monoidal subcategories of finite tensor categories}

We use the terminology of tensor categories from \cite{EtGeNiOs-tensor-categories}.
Recall that an object $X$ of a tensor category $\mathbf{C}$ is a \emph{generator}
if any object of $\mathbf{C}$ is isomorphic to a subquotient of a direct sum of tensor products of $X$ and its (iterated) duals. 

For \autoref{R:FiniteTensor}, we need to show that in case $\mathbf{C}$ is \textit{finite} (meaning that $\mathbf{C}$ has enough projective objects and finitely many simple objects), all projective objects appear as direct summands of the tensor powers of $X$. This follows from \autoref{T:AppendixProjTensor}.(d) below and the fact that projective objects in a tensor category are also injective; see \cite[Proposition 6.1.3]{EtGeNiOs-tensor-categories}. Note that the theorem is a generalization of the known case of fusion categories in \cite[Corollary~F.7]{DGNO}. It can also be interpreted as a generalization of the Burnside--Brauer--Steinberg theorem for Hopf algebras; see, e.g. \cite{PaQu-bbs-hopf-algebras}, or a generalization of the main result of \cite{BrKo-tensor-group}.

\begin{Theorem}\label{T:AppendixProjTensor}
The following properties hold for any tensor category $\mathbf{C}$ with finitely many isomorphism classes of simple objects:
\begin{enumerate}

\item Any full additive subcategory that is closed under taking subquotients and internal tensor products is a tensor subcategory (i.e. rigid).

\item For any $X\in\mathbf{C}$, the left and right dual objects $X^{\ast}$ and $^{\ast}X$ appear as subquotients of direct sums of $X^{\otimes d}$ for $d\in\N$.

\item For any simple object $V\in\mathbf{C}$, the left and right dual objects $V^{\ast}$ and $^{\ast}V$ appear as subquotients of $V^{\otimes d}$ for some $d\in\N$.

\item If $X$ is a generator, then every object in $\mathbf{C}$ is a subquotient of direct sums of tensor powers of $X$. 

\end{enumerate}
\end{Theorem}

\begin{proof}
Firstly, we can observe that part (a) implies part (b) and that (b) implies (d). 

Now we prove that (c) actually implies (a). For this, let $\mathbf{D}\subset\mathbf{C}$ be a subcategory as in (a), and we assume that property (c) is true for $\mathbf{C}$.
To show that $\mathbf{D}$ is rigid, we show that every object $Z\in\mathbf{D}$ is rigid. 
We prove this by induction on the length $\ell(Z)$ of $Z\in \mathbf{D}$. 
The base $\ell(Z)=1$ follows from our assumption (c).

To justify the induction step, let $\ell(Z)\geq 2$ and include $Z$ into a short exact sequence
\begin{gather*}
0\to Y\to Z\to A\to 0,
\end{gather*}
where $Y,A$ are nonzero. Note that $A$ and $Y$ are rigid in $\mathbf{D}$
by induction, and $Z$ is defined by an element $h\in\mathrm{Ext}(A,Y):=\mathrm{Ext}_{\mathbf{D}}(A,Y)$.
Now, in general, given an element
$h\in\mathrm{Ext}(A,Y)$ for objects with duals $A$ and $Y$, we have an
element $1\otimes h\otimes 1\in\mathrm{Ext}(Y^{\ast}\otimes A\otimes A^{\ast},
Y^{\ast}\otimes Y\otimes A^{\ast})$ which defines an element $h^{\ast}\in\mathrm{Ext}(Y^{\ast},A^{\ast})$ obtained by composing $1\otimes h\otimes 1$ with evaluations and coevaluations.
The $Z^{\ast}$, the extension of $A^{\ast}$ and $Y^{\ast}$ defined by $h^{\ast}$ is
the left dual of $Z$, and thus belongs to $\mathbf{D}$. The right dual ${}^{\ast}Z$ can be found in a similar way. Thus, $Z$ is rigid.

By the above, we have (c)$\Rightarrow$(a)$\Rightarrow$(b)$\Rightarrow$(d), and it now suffices to prove property (c). This will be done below after some preparation.

The following lemma can be seen as being part of the Frobenius--Perron theorem. 

\begin{Lemma}\label{lee1} 
Let $M$ be a nonzero square matrix with nonnegative entries and $v,w$ be row and column vectors with strictly positive entries such that $vM=\mu_1 v$, $Mw= \mu_2 w$ for some $\mu_1,\mu_2\in \mathbb{R}$. Then $\mu_1=\mu_2>0$ and $M$ is completely reducible (conjugate by a permutation matrix to a direct sum of irreducible matrices). 
\end{Lemma} 

\begin{proof} 
We have $\mu_1 vw=\mu_2vw=vMw>0$, so $\mu_1=\mu_2=\mu>0$. 
It remains to show that if $M$ is reducible, then it is decomposable. 
If $M$ is reducible, then after conjugating $M$ by a permutation matrix, we have $M =\begin{psmallmatrix} M_{11}& M_{12}\\0 & M_{22}\end{psmallmatrix}$. So, if $v=(v_1,v_2)$,
$w=\binom{w_1}{w_2}$ then we have
\begin{gather*}
M_{22}w_2=\mu w_2,\ v_1M_{12}+v_2M_{22}=\mu v_2.
\end{gather*}
Multiplying the second equation by $w_2$, the first equation by $v_2$, and subtracting, we get $v_1M_{12}w_2=0$, which means $M_{12}=0$, so $M$ is decomposable, as claimed. 
\end{proof}

Now, let $A$ be a transitive unital $\mathbb{Z}_+$-ring of finite rank with basis $B=(b_j)$. 
For $X:=\sum x_jb_j$ with $x_j\in \mathbb{R}$, let $M_X$ be the matrix of left multiplication by $X$ on $A\otimes_{\mathbb{Z}}\mathbb{R}$ in the basis $B$. 

\begin{Lemma}\label{lee2} 
\leavevmode
\begin{enumerate}[label=(\roman*)]
\item If $x_j\ge 0$, then $M_X$ is completely reducible.
\item Suppose that $b,b^\ast\in B$ and that the $B$-decomposition of $bb^\ast$ contains $b_0=1$. Then there exists $n\in \N$ such that the $B$-decomposition of $b^n$ contains $b^\ast$.
\end{enumerate}
\end{Lemma} 

\begin{proof} 
\textit{(i).} Let $d\colon A\to \mathbb{R}$ be the Frobenius--Perron dimension. Since $d$ is a character, $d(X)d(b_j)=d(Xb_j)=\sum_k (M_X)_{jk}d(b_k)$, so the column vector $(d(b_k))$ is a right eigenvector of $M_X$. Also, by \cite[Proposition 3.3.6(2)]{EtGeNiOs-tensor-categories}, there exists a left (row) eigenvector $(c_k)$ of
$M_X$ with positive entries. Thus, by \autoref{lee1}, $M_X$ is completely reducible. 

\textit{(ii).} Consider the directed graph $\Gamma$ with vertex set $B$ and edge $b_i\to b_j$ present if and only if
the $B$-decomposition of $bb_i$ contains $b_j$. Let $\Gamma_1$ be the connected component of $\Gamma$ that contains $b_0=1$. Since the $B$-decomposition of $bb^\ast$ contains $1$, 
$b^\ast\in \Gamma_1$. Thus, by (i), there is an oriented path from $1$ to $b^\ast$, i.e., 
there exists $n\in \N$ such that the $B$-decomposition of $b^n$ contains $b^\ast$, as desired. 
\end{proof} 

To complete the proof of \autoref{T:AppendixProjTensor}, it now suffices to observe that \autoref{T:AppendixProjTensor}.(c) follows from applying \autoref{lee2}(ii) to the Grothendieck ring $A=K_0(\mathbf{C})$, with $B$ the basis of simple objects, $b=[V]$ and $b^\ast=[V^\ast]$ or $b^\ast=[{}^\ast V]$.
\end{proof}

\begin{Remark} 
The assumption that $\mathbf{C}$ has finitely many simple objects in \autoref{T:AppendixProjTensor} cannot be dropped (even for categories with enough projective objects). For instance, the statement is not true if 
$\mathbf{C}$ is the representation category of the multiplicative group and $X$ is
a faithful one dimensional representation.
\end{Remark}

\end{document}